\documentclass{birkjour_t2}
\usepackage[utf8]{inputenc}
\usepackage{graphicx}%
\usepackage{multirow}%
\usepackage{amsmath,amssymb,amsfonts}%
\usepackage{amsthm}%
\usepackage{mathrsfs}%
\usepackage[title]{appendix}%
\usepackage{xcolor}%
\usepackage{textcomp}%
\usepackage{manyfoot}%
\usepackage{booktabs}%
\usepackage{algorithm}%
\usepackage{algorithmicx}%
\usepackage{algpseudocode}%
\usepackage{listings}%
\usepackage{subfigure}
\usepackage[backend=bibtex,
            style=nature,
           citestyle=numeric-comp,
           maxnames=99,
           sorting=none]{biblatex}
\usepackage{hyperref} 
\usepackage{multirow}
\usepackage{array}
\usepackage{booktabs}

\newtheorem{theorem}{Theorem}
\newtheorem{defn}{Definition}[section]
\newtheorem{rem}[defn]{Remark}
\definecolor{colorluis}{rgb}{0.8, 0.31, 0.22}

 \definecolor{tiffanyblue}{rgb}{0.04, 0.73, 0.71} 

\providecommand{\testright}[1]{\shortmid\!\xrightarrow{#1}}

\begin{document}

%
%
%
%
%
%
%
%
%

\title[The Riemann problem for three-phase foam flow in porous media]{The Riemann problem for three-phase foam flow in porous media}







\author[L. F. Lozano ]{Luis Fernando Lozano}
\address{%
Laboratory of Applied Mathematics\\
Federal University of Juiz de Fora\\
36036-900\\
Juiz de Fora, MG\\
Brazil}
\email{luisfer99@gmail.com}

\thanks{G. C. and L. L. gratefully acknowledge support from Shell Brasil through the project ``Avan\c{c}ando na modelagem matem\'atica e computacional para apoiar a implementa\c{c}\~ao da tecnologia `Foam-assisted WAG' em reservatórios do Pr\'e-sal'' (ANP 23518-4)  at UFJF and the strategic importance of the support given by ANP through the R\&D levy regulation.
All authors were partly supported by CNPq grant 405366/2021-3.
G. C. was partly supported by CNPq grant 306970/2022-8, and FAPEMIG grant APQ-00206-24.
}
\author[G. Chapiro]{Grigori Chapiro}
\address{%
Laboratory of Applied Mathematics\\
Federal University of Juiz de Fora\\
36036-900\\
Juiz de Fora, MG\\
Brazil}
\email{grigori.chapiro@ufjf.br}
\author[D. Marchesin]{Dan Marchesin}
\address{
National Institute of Pure and Applied Mathematics\\
Estrada Dona Castorina 110\\
Rio de Janeiro, RJ\\
Brazil}
\email{marchesi@impa.br}
\subjclass{
35L65; 
76S05; 
76T30. 
}

\keywords{Riemann problem, three-phase foam flow, porous media, wave curve method}


\begin{abstract}

Gas injection in the context of the three-phase flow in porous media appears in applications such as Enhanced Oil Recovery, aquifer remediation, and carbon capture, utilization, and storage (CCUS). In general, this technique suffers from a difficulty related to excessive gas mobility, which can be circumvented by using foam.
This study addresses the non-linear system of differential equations describing the three-phase foam flow based on Corey relative permeability functions. A major obstacle is an umbilic point, where the characteristic wave velocities for different families coincide, complicating the identification of stable wave structures.

We developed a methodology to solve the Riemann problem describing the three-phase foam displacement in the case when the gas viscosity exceeds that of oil and water. To allow the analysis, we assume foam in local equilibrium (or maximum foam texture), resulting in a constant mobility reduction factor (MRF). These simplifications allowed the classification of possible solutions for the injection of foamed gas and water mixtures under a wide range of initial conditions within the framework of non-classical Conservation Law Theory.
As a relevant industrial application of the proposed solution, we investigate the conditions resulting in oil bank formation.
Besides improving the general physical understanding of foam flow in a porous medium, this analysis can be applied to calibrate numerical simulators and perform uncertainty quantification.
Our analytical estimates were validated through numerical simulations.
\end{abstract}

\maketitle
\section{Introduction}\label{sec1}
Gas injection within three-phase flow in porous media framework appears in several applications, such as $CO_2$ sequestration and underground storage (CCUS), soil remediation, and enhanced oil recovery, see \cite{metz2005ipcc,salimi2012influence,iskandarov2022data}. 
One of the main limitations of gas injection is related to the fingering formation limiting gas storage capacity (for CCUS) and oil displacement (for enhanced oil recovery (EOR) and soil remediation). Foam is a promising technique for controlling this issue. Stable lamellae form in the presence of surfactant, breaking the continuous connection of the gas phase and significantly limiting the gas mobility \cite{Farajzadeh2012,ma2015modeling,hematpur2018foam}.

Numerical simulations of foam flow are challenging and typically do not possess convergence analysis due to the complexity of the partial differential equations involved \cite{chen2010modeling,dePaula2020,DEPAULA2023,de2024numerical}.
This is one reason why analytical solutions to complex phenomena play a key role in improving numerical simulator validity and robustness.
Modeling the interaction between oil, water, and gas is challenging due to intricate wave interactions. Adding foam typically increases the complexity of the problem. 
Modeling and solving analytically the two-phase foam flow in a porous medium is already a challenging task \cite{Ashoori2011a,Zavala2021,Lozano2021,lozano2022,chapiro2022analytical,fritis2024riemann,danelon2024,danelon2025modeling}. Only few studies have addressed this problem for three-phase flow. 
For example, \cite{mayberry2008use} considered foam in equilibrium as a gas with maximum apparent viscosity. They study oil displacement in various foam flooding scenarios, including different viscosities and foam strengths. However, they assumed linear relative permeabilities limiting real applications. 
In \cite{namdar2011method}, the results of \cite{mayberry2008use} were extended to cases where foam strength is a function of water or oil saturation or both. 
In \cite{lee2014three}, the authors also presented analytical and numerical solutions considering linear relative permeabilities but simplifying the foam strength functions. 
Later, in \cite{lee2015moc}, the authors extended these results to multi-layered systems. 
In \cite{tang2019three}, the authors investigated solutions for an implicit-texture foam model with oil using a Corey model with realistic permeability exponents with the foam model implemented in the commercial simulator CMG/STARS, which is widely used in the oil industry. They studied solutions with multiple steady states and explained which one corresponds to stable displacement and is physically correct. Later, \cite{tang2022foam} used analytical solutions for the same foam model with oil, incorporating the Corey model with quadratic permeabilities, to study how the modification of surfactant type (with different oil tolerances) affects the solution of the Riemann problem. Both last two works did not classify all possible Riemann problems due to the complexity of the solution construction. 
In the present work, we assume equilibrium foam as in \cite{mayberry2008use}, and we use convex (quadratic) relative permeabilities as in \cite{tang2022foam} obtaining the solution, which is qualitatively equivalent to one obtained by \cite{tang2019three,lyu2021simulation,tang2022foam}.

In this work, we study the classification of Riemann problem solutions for the
injection of foamed gas and water mixtures under a wide range of initial conditions relevant to industrial applications. Our results show that injecting a two-phase fluid mixture, rather than a single pure phase can enhance oil displacement efficiency due to the formation of an oil bank, which is a crucial phenomenon in petroleum engineering applications.
We follow non-standard Conservation Laws Theory to construct a solution as a sequence of waves (shocks, rarefactions, and composite waves) following \cite{Azevedo2014,aparecido1995,V.2010,Lozano2020pro}.
The specialized software \cite{ELI_web}, developed at the Laboratory of Fluid Dynamics of IMPA, was used to construct these wave curves. The software package enabled us to obtain and analyze integral curves, Hugoniot curves, main bifurcation loci, and phase portraits of dynamical systems and wave curves. These elements are essential for constructing Riemann solutions, as mentioned in \cite{V.2010,Azevedo2014,lozano2018} and other related works. Additionally, our analytical estimates are validated through direct numerical simulations.

This work is organized as follows. 
Section \ref{sec:model} presents a mathematical model describing three-phase foam flow in a porous medium. Section \ref{sec:elementary_waves} presents a brief review of elementary waves, and Section \ref{sec:wave_curve} describes the wave curve method. Section \ref{sec:Classification} presents the construction of the solution for the three-phase foam flow model together with its classification. Section \ref{sec:Results} compares our results with those in the literature for realistic application problems. Finally, Section \ref{sec:discusstion} ends with some discussions and conclusions.


\section{Mathematical Model}
\label{sec:model}

This study focuses on the three-phase flow in a porous medium consisting of water, oil, and foamed gas with constant viscosities. Following \cite{mayberry2008use,mehrabi2020solution,danelon2024, lozano2024a}, the foam is considered to be in local equilibrium, \textit{i.e.,} foam texture remains at its maximum level. We make several assumptions: one-dimensional horizontal flow, incompressible fluids, negligible effects from dispersion, gravity, and capillarity, complete saturation of the rock pores by the fluids, constant temperature, and no mass transfer between the phases. Under these assumptions, the model is described by the following system of conservation equations:
\begin{equation}
\label{eq:Threephase_0}
    \begin{cases} 
	\phi\dfrac{\partial S_w}{\partial t} + \dfrac{\partial u_w}{\partial x} = 0, \\[7pt]
	\phi\dfrac{\partial S_o}{\partial t} + \dfrac{\partial u_o}{\partial x} = 0,
	 \\[7pt]
 \end{cases}  
\end{equation}
where $S_w(x,t)$ and $S_o(x,t)$ represent water and oil saturations; $x\in\mathbb{R}$ and $t\geq 0$. The parameter $\phi$ is the porosity of the medium, while $u_w$ and $u_o$ are the superficial velocities of water and oil.

Using fractional flow theory \cite{lake1989enhanced}, we can express these velocities as:
\begin{equation}
\label{eq:Darcy_Law}
u_i = u\, f_i,  \qquad i = w,o,g,
\end{equation}
where $u = u_w + u_o + u_g$ denotes the total surface velocity. The relative permeability functions are modeled using the Corey approach \cite{Corey1954}:
\begin{equation}
\label{eq:permeabilities}
    k_{ri} = k_{ri}^0 \left(\dfrac{S_i - S_{ri}}{1-S_{rw}-S_{ro}-S_{rg}}\right)^{n_i}, \qquad i = w,o,g,
\end{equation}
where $k_{ri}^0 $ is the endpoint permeability, $S_i$ is the saturation, $S_{ri}$ is the residual saturation, and $n_i$ is the exponent for each phase $i$.
The fractional flow functions are:
\begin{equation}
\label{eq:fi}
    f_{i} = \dfrac{k_{ri}/\mu_i}{k_{rw}/\mu_w + k_{ro}/\mu_o + k_{rg}/(\mu_g \times MRF)}, \qquad i = w,o,g,
\end{equation}
where $\mu_i$ is the viscosity of phase $i$, and $MRF$ is the foamed gas mobility reduction factor, which applies only to foamed gas viscosity \cite{Ashoori2011a, zhou1995applying}. Here, we assume that the foam $MRF$ is constant.

The conservation system \eqref{eq:Threephase_0} can be written in a dimensionless vector form:
\begin{equation}
\label{eq:Threephase_1}
   U_{t_D} + F(U)_{x_D} = 0, \quad F(U) = (f_w(U), f_o(U))^T,
\end{equation}
\begin{equation}
     U(x_{D},t_{D}) = (s_w(x_{D},t_{D}), s_o(x_{D},t_{D}))^T,
\end{equation}
where
\begin{equation}
    x_D = \dfrac{x}{\widehat{x}}, \qquad t_D = \dfrac{t \, u}{\phi \, \widehat{x} (1 - S_{rw} - S_{ro} - S_{rg})}, 
\end{equation}
$\widehat{x}$ refers to the reference length of the porous core and $s_i$ is the normalized saturation given by
\begin{equation}
    s_i = \dfrac{S_i - S_{ri}}{1 - S_{rw} - S_{ro} - S_{rg}}, \quad i = w,o,g.
\end{equation}

\begin{figure}[h!]
	\centering
	\includegraphics[scale=0.45]{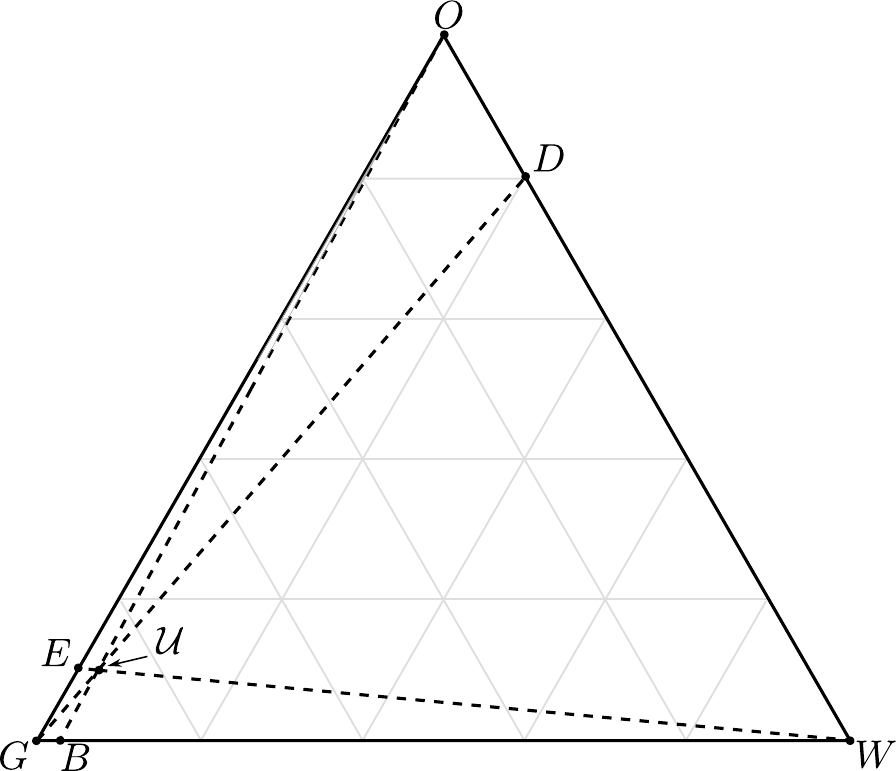}  
	\caption{Saturation triangle $\Omega$. The point $\mathcal{U}$ represents the umbilic point for the parameter values in Table~\ref{tab:table1}. The vertices $W$, $O$, and $G$ correspond to the states with coordinates $(1,0)$, $(0,1)$, and $(0,0)$, respectively.}
	\label{fig:fig1}
\end{figure}

The saturations are defined in a domain $\Omega$, which is contained in a cube bounded by $0 \leq s_w \leq 1$, $0 \leq s_o \leq 1$, and $0 \leq s_g \leq 1$. Given that  $s_w + s_o + s_g = 1$, the domain $\Omega$ is typically visualized as an equilateral triangle using barycentric coordinates, known as the saturation triangle; see Fig.~\ref{fig:fig1}. The vertices $W$, $O$, and $G$ correspond to states with saturations $(s_w,s_o)=(1,0),$ $(s_w,s_o)=(0,1)$, and $(s_w,s_o)=(0,0)$, respectively.

Following  \cite{V.2010, Azevedo2014, Lozano2020pro,Marchesin2001}, we adopt the Corey model with quadratic exponents to describe relative permeabilities ($n_w = n_o = n_g = 2$). This approach simplifies the analysis of mathematical properties while effectively capturing the intricate nonlinear dynamics of wave propagation and interactions in porous media.

The initial condition for system \eqref{eq:Threephase_1} is:
\begin{equation}
\label{eq:RP}
U(x, 0) = \begin{cases}
L = (s_w^L, s_o^L)^T, & \text{if } x < 0, \\
R = (s_w^R, s_o^R)^T, & \text{if } x > 0,
\end{cases}
\end{equation}
where $s_i^L$ and $s_i^R$ ($i = w,o$) are constant. The left state $L$ corresponds to the injection state, while the right state $R$ represents the initial condition in the reservoir. This setup defines a classical Riemann problem.

\subsection{Starting the model analysis}

The Riemann problem associated with equations \eqref{eq:Threephase_1} and \eqref{eq:RP} typically admits solutions which are called waves. The nature of these waves, their interactions, and the analytical approaches used to study them depend on the classification of the system \eqref{eq:Threephase_1} \cite{Bressan2013}. A system of partial differential equations (PDEs) is classified as {\it strictly hyperbolic} if the Jacobian matrix of the flux function, $ DF(U) $, possesses only real and distinct eigenvalues: $ \lambda_1(U) \neq \lambda_2(U) $, for every $ U \in \Omega $. These eigenvalues are known as the {\it characteristic speeds} of the system.

If there are points in $ \Omega $ where $ \lambda_1(U) = \lambda_2(U) $, these points are referred to as {\it umbilic points} \cite{shearer1987} or {\it coincidence points} \cite{mederios1992}, and the system is classified as {\it not strictly hyperbolic} \cite{fritis2024riemann,danelon2024,Marchesin2001}. The model analyzed in this work exhibits an umbilic point (indicated by $\mathcal{U}$ in Fig.~\ref{fig:fig1}) within the interior of $\Omega$ \cite{Marchesin2014, mederios1992, schaeffer1987classification}. 

To improve readability, we follow the literature \cite{Marchesin2001, Azevedo2014} by associating all quantities corresponding to the smaller and larger eigenvalues with the {\it slow family} and {\it fast family}, respectively. Specifically, $ \lambda_s = \lambda_1 $ is the slow-family characteristic speed, and $ \lambda_f = \lambda_2 $ is the fast-family characteristic speed.

As demonstrated in \cite{Azevedo2014}, the characteristic speeds associated with the system \eqref{eq:Threephase_1} are real and positive within the saturation triangle. Along the triangle boundaries, $\lambda_s$ becomes zero because there is a reduction to two-phase flow, which is governed by the Buckley-Leverett equation. Conversely, $\lambda_f$ remains positive everywhere except at the vertices of the triangle, where it also reaches zero. Throughout the closed triangle, the inequality $ \lambda_s < \lambda_f $ holds, except the points $G$, $ W $, $ O $, and $ \mathcal{U} $.

The solution to the Riemann problem depends on the location of the umbilic point, as well as the associated bifurcation and wave curves, as detailed in \cite{isaacson1992multiphase, aparecido1995, V.2010, Azevedo2014, Lozano2020pro}. As shown in \cite{schaeffer1987classification, Marchesin2014, V.2015}, for the system \eqref{eq:Threephase_1}, the umbilic point $ \mathcal{U} $ is unique within $ \Omega $. It is the intersection of the three straight segments $ [G, D] $, $ [W, E] $, and $ [O, B] $ (see Fig.~\ref{fig:fig1}). The coordinates of the umbilic point are given by: 
\begin{equation}
\label{eq:CoordUm}
\mathcal{U} = \left(\mu_w, \mu_o, \mu_g\right) / (\mu_w + \mu_o + \mu_g).
\end{equation}
In this work, we denote the segment between points $A$ and $B$ in $\Omega$ as $[A, B]$; see Fig.~\ref{fig:fig1}. 
Additionally, we use $ (A, B) $ to denote the open segment connecting $ A $ to $B$, excluding both endpoints; $ (A, B] $ and $ [A, B) $ to exclude only one of the endpoints.

In \cite{Andrade2018}, the following viscosity ratios were introduced: \begin{equation}\label{eq:razoes}
\mathcal{R}_{\alpha \beta}^{\pm}=\left(\mu_{\alpha} \pm \mu_{\beta}\right) /
\mu_{\gamma}, \quad \alpha,\,\beta\,,\gamma\in\{w,o,g\}.
\end{equation} 

According to \cite{Andrade2018, lozano2018}, if the viscosity ratios satisfy certain conditions, part of the bifurcation loci lie outside of $\Omega$ so that it is not used, simplifying the classification of the corresponding Riemann problem solutions for $L \in [G, W]$ and right states in the entire region $\Omega$. In \cite{lozano2024b}, Riemann problems were classified considering the viscosity values such that 
\begin{equation} \label{eq:visc_lozano}
\mu_g<\mu_w<\mu_o,\quad \left(\mathcal{R}{g o}^{-}\right)^{2} / \mathcal{R}{g o}^{+} > 8 \quad \text{and} \quad \left(\mathcal{R}{w o}^{-}\right)^{2} / \mathcal{R}{w o}^{+} > 8,
\end{equation}
corresponding to applications dealing with high-viscosity heavy oil.
The umbilic point in \cite{lozano2024b} is close to the corner $O$ in $\Omega$ for these parameters. In \cite{lozano2024b}, the authors present solutions to the Riemann problem considering $L \in [G, W]$ and $R$ in a large portion of the saturation triangle. In the present paper, we investigate foam flow, which results in high apparent gas viscosity displacing the umbilic point to a neighborhood of the point $G$. 
Although the same injection state remains $L \in [G, W]$, the wave sequences in the construction of the solutions are substantially different from those appearing in \cite{lozano2024b}.
In this sense, the present work complements \cite{lozano2024b}.

In the literature, we often find viscosity values satisfying \begin{equation} \label{eq:visc_foam}
\mu_w < \mu_o < \mu_g \times MRF, 
\end{equation} 
when portraying the maximum strength of the foam, with relationships such as
\begin{equation} 
\mu_w / (\mu_g \times MRF) = 10^{-2} \quad \mbox{and}\quad \mu_o / \mu_g = 10^{-1}. 
\end{equation} 
This work considers the foam strength to be $MRF = 1750$, similar to values reported in \cite{gassara2020calibrating} and \cite{tang2022foam}. 
The viscosity parameters are chosen as $\mu_w = 1$, $\mu_o = 4$, and $\mu_g = 0.02$ as summarized in Table \ref{tab:table1}. The foamed gas viscosity becomes $\widehat{\mu_g} =\mu_g \times MRF =  35$ with these values. As a result, the umbilic point lies close to the vertex $G$ in $\Omega$ (Fig.~\ref{fig:fig1}) and satisfies 
\begin{equation}\label{eq:desi_param_foam} 
\left(\mathcal{R}{g w}^{-}\right)^{2} / \mathcal{R}{g w}^{+} > 8 \quad \text{and} \quad \left(\mathcal{R}{g o}^{-}\right)^{2} / \mathcal{R}{g o}^{+} > 8. 
\end{equation}
Strong numerical evidence suggests that the boundary structure of the region where Riemann solutions are considered remains stable for viscosity parameters $ \mu_w $, $ \mu_o $, and $ \widehat{\mu_g} $ satisfying \eqref{eq:desi_param_foam} and $ MRF \geq 1744.43 $.
\begin{table}[h]
\centering
\caption{Parameters for Corey model. Viscosities from \cite{mehrabi2022displacement}.}
\label{tab:table1}%
\begin{tabular}{ccccccccccc}
\hline
$n_w$& $n_o$  & $n_g$ & $k_{rw}^0 $& $k_{ro}^0 $  & $k_{rg}^0 $ &$\mu_w$& $\mu_o$  & $\mu_g$  & $MRF$ &$\widehat{\mu_g} =\mu_g \times MRF$\\
\hline
2    & 2   & 2 &1    & 1   & 1 &1    & 4   & 0.02  &1750 &35 \\
\hline
\end{tabular}
\end{table}

\section{Elementary waves}
\label{sec:elementary_waves}
In general, the solution to a Riemann problem consists of a sequence of wave groups separated by constant states. A wave group is a sequence of consecutive waves with no speed gaps between successive waves, as described in \cite{Marchesin2001, V.2010}. In this paper, we consider a model whose solutions may consist of combinations of three fundamental wave types: rarefaction waves, shock waves, and composite waves. 

\begin{figure}
    \centering
    \subfigure[Rarefaction curves for the slow characteristic speed.]{\includegraphics[width=0.485\linewidth]{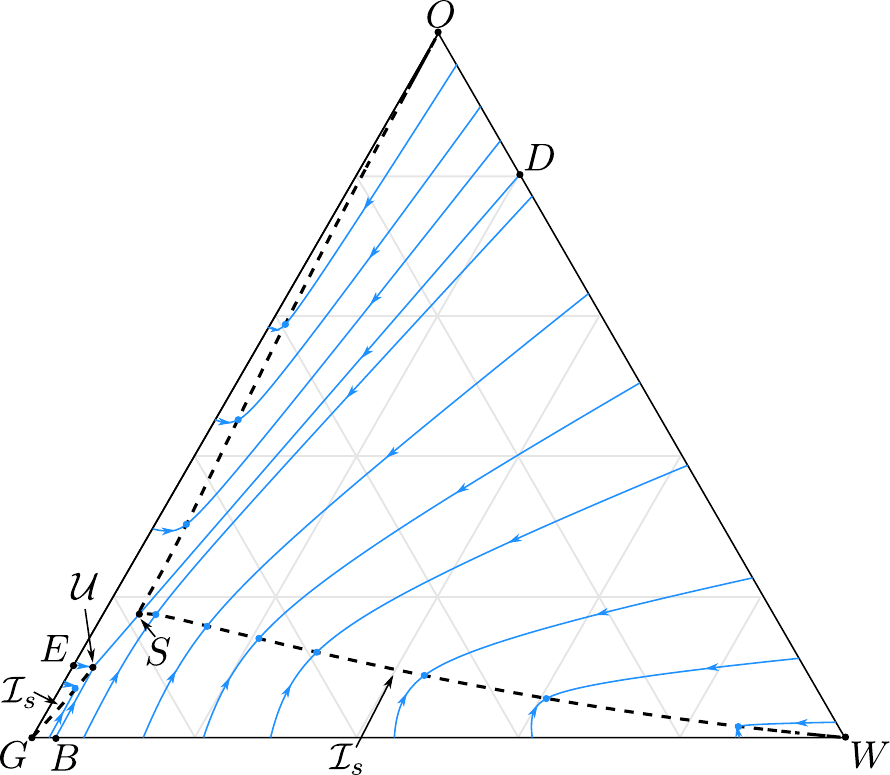}}\hspace{5pt}
      \subfigure[Rarefaction curves for the fast characteristic speed.]{\includegraphics[width=0.485\linewidth]{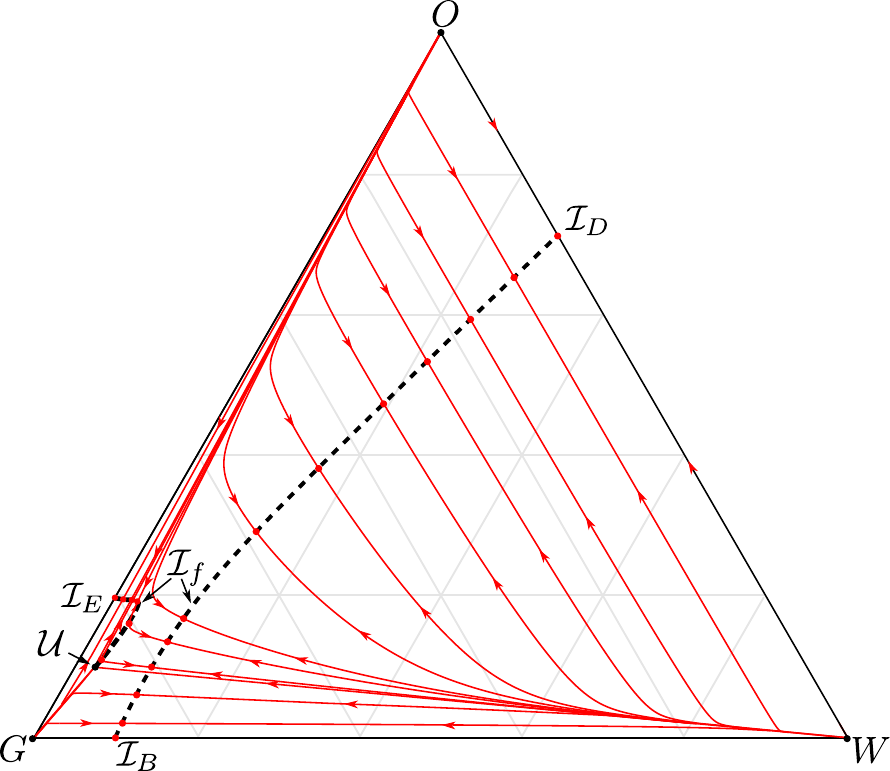}}
    \caption{Rarefaction curves for the parameter values in Table~\ref{tab:table1}. The arrows indicate the direction of increasing characteristic speed. $\mathcal{I}_i$ (dashed lines) represents the $i$-inflection locus of the $i$-family (Definition \ref{def:iinflection}) (a) The state $S$ is the intersection between the $s$-inflection and the segment $[G,D]$.}
    \label{fig:fig2}
\end{figure}

Rarefaction waves are smooth, self-similar solutions to system \eqref{eq:Threephase_1}, expressed as 
\begin{equation}
\label{eq:rarefact} 
U(x,t) = \widehat{U}(\xi),\qquad \xi = x/t. 
\end{equation} 
By substituting \eqref{eq:rarefact} into system \eqref{eq:Threephase_1}, the rarefaction curve is derived from solving the eigenvalue problem \begin{equation}\label{eq:systemEDO_rarefaction} (DF(\widehat{U}) - \lambda_i(\widehat{U}) I)r_i(\widehat{U}) = 0,
\end{equation} 
where $I$ is the identity matrix, and $(\lambda_i(\widehat{U}),r_i(\widehat{U}))$, $i = s, f$, are the eigenpairs associated with the Jacobian matrix of the flux function $DF$. The eigenvector $r_i$ is aligned with $d\widehat{U} /d\xi$, with $\xi = \lambda_i, i = s, f$. 
Figures \ref{fig:fig2}(a) and \ref{fig:fig2}(b) illustrate the rarefaction curves associated with the slow- and fast-family in this model, corresponding to the parameter values in Table \ref{tab:table1}. The black curve $\mathcal{I}_i$ represents the {\it inflection locus} of the $i$-family and is the locus where the characteristic speed reaches an extremum (see Definition \ref{def:iinflection}).

Shock waves, on the other hand, are bounded, discontinuous, self-similar solutions of system \eqref{eq:Threephase_1}, connecting the left state $U^-$ to the right state $U^+$ through a jump (discontinuity), which satisfies the Rankine-Hugoniot condition: \begin{equation}\label{eq:RH1} 
F(U^+)-F(U^-)-\sigma\left( U^+ - U^-\right)=0, \end{equation} 
where $\sigma=\sigma(U^-; U^+)$ represents the propagation speed of the discontinuity. For a fixed $L\in \Omega$, the set of states $U$ such that the pair $(L; U)$ satisfies the Rankine-Hugoniot condition \eqref{eq:RH1} for some $\sigma$ are called the Hugoniot locus for the state $L$ and is denoted by $\mathcal{H}(L)$.

A composite wave, consisting of either a rarefaction followed by a shock or vice versa, is a solution of the Riemann problem \eqref{eq:Threephase_1}, \eqref{eq:RP}, where these two wave types occur in sequence without a constant state segment in between. In this work, we only consider composite waves where both waves belong to the same family.

\subsection{Admissibility criterium}

Since discontinuous solutions in Riemann problems may lead to non-unique solutions, equation \eqref{eq:RH1} must be complemented with an additional criterion to select the unique, physically meaningful solution. These criteria, known as entropy or admissibility conditions, include the ones proposed by Lax \cite{Lax1957}, Liu \cite{Liu1975}, the viscous profile approach \cite{Gelfand1959}, and the vanishing adsorption method \cite{Petrova2022}. It is important to note that different admissibility criteria may fail to be equivalent (see \cite{Schecter1996, Gomes1989}). This study focuses on constructing solutions for shocks that satisfy the viscous profile criterion.


\subsubsection{Viscous profile criterion}\label{sec:visous_profile_Cr}
Considering the parabolic system
\begin{equation}
    \label{eq:parabolic_EPD}
    \dfrac{\partial }{\partial t_D} U + \dfrac{\partial }{\partial x_D}F(U) =\varepsilon \dfrac{\partial}{\partial x_D}\left(D(U)\dfrac{\partial}{\partial x_D} U\right)\quad x_D\in\mathbb{R}\,\,,t_D\geq0\,,U\in\Omega,
\end{equation}
obtained from \eqref{eq:Threephase_1} by considering the viscosity effects (related to the capillary pressure between different phases). 
Notice that, system \eqref{eq:Threephase_1} is an approximation of \eqref{eq:parabolic_EPD} after applying the limit $\varepsilon \xrightarrow[]{}0^+$ \cite{L.1990}.

A shock connecting $U^-$ to $U^+$ satisfies the viscous profile criterion if it is the limit of the {\it traveling wave} solution for the parabolic system \eqref{eq:parabolic_EPD} when the positive parameter $\varepsilon$ approaches zero. The traveling wave solution is given by 
\begin{equation}
    U(x_D, t_D) = \widehat{U}(\xi)\quad\mbox{ with }\quad \xi =(x_D - \sigma \,t_D)/\varepsilon,
\end{equation}
satisfying the boundary conditions
\begin{equation}
    \widehat{U}(-\infty) = U^-\quad \mbox{and}\quad \widehat{U}(+\infty) = U^+,
\end{equation}
where $\sigma = \sigma(U^-; U^+)$ represents the shock speed.  Since the traveling wave $\widehat{U}$ is a smooth function of $\xi$, the system can be reduced to a set of ordinary differential equations (ODEs) (omitting hats):
\begin{equation}
\label{eq:dinamycsyst}
 D(U (\xi)) \dfrac{d U(\xi)}{d \xi} = -\sigma \left(U (\xi) - U^-\right) + F (U (\xi)) - F (U^- ),   
\end{equation}
where $U^-$ and $U^+$ are equilibrium points of \eqref{eq:dinamycsyst}, \textit{i.e.,} the right side of system  \eqref{eq:dinamycsyst} is zero when evaluated at $U^-$ and at $U^+$. A traveling wave solution corresponds to a connection between $U^-$ and $U^+$ and is known as a viscous profile for the shock wave.

In the general case, the matrix $D(U)$ depends on relative permeabilities, viscosities, and capillary pressure function derivatives. In \cite{azevedo2002capillary}, it was proved that $D(U)$ is positive definite in the interior of $\Omega$. In this work, we follow \cite{aparecido1995, isaacson1992multiphase, V.2010, Azevedo2014} and assume that $D(U)$ is a multiple of the identity matrix.
This assumption ensures that the ODE system \eqref{eq:dinamycsyst} remains invariant along the segments $[G, D]$, $[W, E]$, and $[O, B]$. This implies, among other things, that the admissibility of nonlocal shocks depends on the positions of $U^-$ and $U^+$ relative to segments $[G, D]$, $[W, E]$, or $[O, B]$; see \cite{lozano2018}.

\subsubsection{Lax  admissibility criterion}
We can classify the types of discontinuities that satisfy \eqref{eq:RH1} and the Lax admissibility criterion as follows:  
\begin{itemize}
    \item A discontinuity that satisfies  
\begin{equation}
    \sigma <\lambda_s(U^-) \quad \mbox{and} \quad \lambda_s(U^+)<\sigma<\lambda_f(U^+)
\end{equation}
is called a {\it slow shock} (or $s$-shock, or Lax 1-shock \cite{Lax1957}). 
    \item For a {\it fast shock} ($f$-shock or Lax $2$-shock \cite{Lax1957}), the discontinuity must satisfy  
\begin{equation}
   \lambda_s(U^-)<\sigma<\lambda_f(U^-) \quad \mbox{and} \quad \lambda_f(U^+)<\sigma.
\end{equation}
\end{itemize}

For this model, we also have two other relevant discontinuities, known as non-Lax shocks (non-classical shocks) \cite{isaacson1992multiphase, V.2015}:  
\begin{itemize}
    \item In the case of {\it undercompressive shocks} ($u$-shocks or {\it transitional shocks}), we have  
        \begin{equation}
            \lambda_s(U^-)<\sigma<\lambda_f(U^-) \quad \mbox{and} \quad \lambda_s(U^+)<\sigma<\lambda_f(U^+).
        \end{equation}
    \item In the case of {\it overcompressive shocks} ($O$-shocks), we have  
        \begin{equation}
             \lambda_f(U^+)<\sigma<\lambda_s(U^-).
        \end{equation}
\end{itemize}

From the theory of ODEs, there is a relation between the type of shock and the nature of the equilibrium states $L$ and $R$ \cite{V.2015}:
\begin{itemize}
    \item $s$-shock, denoted by $S_s$, $U^-$ is a repeller and $U^+$ is a saddle.
    \item $f$-shock, denoted by $S_f$, $U^-$ is a saddle and $U^+$ is an attractor.
    \item $u$-shock, denoted by $S_u$, $U^-$ and $U^+$ are saddles.
    \item $O$-shock, denoted by $S_O$, $U^-$ is a repeller and $U^+$  is an attractor.
\end{itemize}
In this work, as $ D(U) $ is replaced by the identity matrix, $ u $-shocks are admissible for states $U^-$ and $U^+$ along the segments $[G, D]$, $[W, E]$ and $[O, B]$; see \cite{Marchesin2001, isaacson1992multiphase, lozano2018, Lozano2020pro, lozano2024c, V.2015, Andrade2018}. For the general case of $ D(U) $, see, for example, \cite{lozano2018, Lozano2020pro, lozano2024c}.

 \subsection{Notations}
Let us indicate notation for shocks and rarefactions following \cite{V.2015}. 
A shock from state $A$ to state $B$ may be denoted by $A \xrightarrow{S} B$. 
More generally, a shock of type $X$ (where $X$ can be $s$, $f$, $u$, or $O$) between states $A$ and $B$ is written as $A \xrightarrow{S_X} B$ (if $A$ is a single point in the characteristic $xt$-space) or $A \testright{S_X} B$ (if $A$ corresponds to a positive measure area in the characteristic $xt$-space). 
The same convention applies to rarefactions: an $i$-rarefaction, with $i = s$ or $f$, is denoted by $A \xrightarrow{R_i} B$ or $A \testright{R_i} B$, according to characteristic $xt$-space.

Let $A$, $B$, $C$ be constant states such that $A$ is connected to $B$ by an $a$-wave $ A\xrightarrow{a} B$, and $B$ is connected to $C$ by a $b$-wave $B\xrightarrow{b}C$ (slow, fast, or undercompressive wave group) resulting in the sequence $ A\xrightarrow{a} B\xrightarrow{b}C$. We denote by $v_\text{f}^a$ the final velocity of the $a$-wave, and by $v_\text{i}^b$ the initial velocity of the $b$-wave. The wave sequence is said to be \emph{compatible} if and only if 
\begin{equation}\label{eq:compatib}
    v_\text{f}^a \leq v_\text{i}^b.
\end{equation}
The compatibility is analogous for the larger wave sequences.
 
\subsection{Bifurcation loci}
In this section, we define and characterize various bifurcation loci, which are regions in state space where the topology of solutions undergoes qualitative changes. Understanding these structures is essential for analyzing solution behavior, particularly in cases where the classical method for constructing solutions to the Riemann problem fails; \cite{isaacson1992multiphase, V.2015}.

\begin{defn}\label{def:iinflection}
A state $U$ belongs to the $i$-inflection locus if
\begin{equation}
    \nabla \lambda_i(U)\cdot r_i(U) = 0,
\end{equation}
where $r_i(U)$ is the right eigenvector of the Jacobian matrix $DF(U)$ associated with $\lambda_i(U)$,  $i\in\{s,f\}$.
\end{defn}

\begin{defn}\label{def:seconBif}
A state $U$ belongs to the {\it secondary bifurcation locus} for the family $i\in\{s,f\}$ if there exists a state $U' \neq U$, such that
\begin{equation}
U' \in \mathcal{H}(U) \mbox{  with  } \lambda_i(U') = \sigma(U; U') \mbox{  and  } l_i(U') \cdot (U'-U) = 0,  
\end{equation}
where $l_i(U')$ is the left eigenvector of the Jacobian matrix $DF(U')$ associated with $\lambda_i(U')$.
\end{defn}

In our model, the secondary bifurcation loci are composed of the segments $[G, D]$, $ [W, E]$, and $ [O, B] $,  see Fig. \ref{fig:fig1}. In \cite{Azevedo2014}, it was shown that the segments $ [E, \mathcal{U}] $, $ [B, \mathcal{U}] $, and $ [D, \mathcal{U}] $ remain invariant under the slow-characteristic vector field, whereas $ [G, \mathcal{U}] $, $ [W, \mathcal{U}] $, and $ [O, \mathcal{U}] $ are invariant under the fast-characteristic vector field.



\begin{defn}\label{def:iLRExt}
Let $\mathcal{L}$ be a locus (a curve with possible self-intersections) in state space. For $i \in \{s,f\}$, we define the sets:
 \begin{itemize}
     \item $i$-left-extension of $\mathcal{L}$:
     \begin{equation}
     \{U^+ \in \Omega : \exists\, U^- \in \mathcal{L} \,\text{ such that }\, U^+ \in \mathcal{H}(U^-),\, \sigma(U^-; U^+) = \lambda_i(U^-)\}.
     \end{equation}
     \item $i$-right-extension of $\mathcal{L}$:
     \begin{equation}
     \{U^+ \in \Omega : \exists\, U^- \in \mathcal{L} \,\text{ such that }\, U^+ \in \mathcal{H}(U^-),\, \sigma(U^-; U^+) = \lambda_i(U^+)\}.
     \end{equation}
 \end{itemize}
\end{defn}

\begin{defn}\label{def:extension_umbilic}
A state $U$ belongs to the extension of the umbilic point $\mathcal{U}$ if it satisfies
\begin{equation}\label{eq:ext_um}
    \lambda_s(\mathcal{U}) = \lambda_f(\mathcal{U}) = \sigma(\mathcal{U}; U).
\end{equation}
\end{defn}
In this work, the extension of the umbilic point consists of three states, each located within a secondary bifurcation segment: $[W, E]$, $[G, D]$, and $[O, B]$. These states are denoted by $E_0$, $D_0$, and $B_0$, respectively—one in each segment.




\section{Wave curve Method}
\label{sec:wave_curve}

A wave curve is characterized by its initial state, family, and direction (forward or backward). A forward wave curve for the $i$-family is a trajectory that contains all states $U$ within the $i$-wave group in $\Omega$, parameterized by increasing velocity $\lambda$ or $\sigma$ ($i = {s, f}$), starting from the state $L$. Conversely, the backward wave curve for the $i$-family is a path consisting of all states $U$ in the $i$-wave group, parameterized by decreasing velocity, beginning from the state $R$.
We denote $\mathcal{W}_i(R)$ the backward $i$-wave curve from the state $R$ ($i = {s, f}$). 

Typically, wave curves were constructed using the local continuation algorithm, which was initially developed by Liu in \cite{Liu1975} to extend Ole\u{i}nik approach for scalar Riemann solutions \cite{Oleinik1957}. However, in cases where the system is non-strictly hyperbolic, the traditional approach may fail. This failure occurs because of nonlocal branches of the Hugoniot locus \cite{Marchesin2001,isaacson1992multiphase}, implying that the wave curves in this model may exhibit disconnected or nonlocal segments (we refer to the local branch of $\mathcal{W}_i(N)$ as the branch that contains the initial state $N$). Therefore, in this work, we apply the {\it global continuation algorithm} to build wave curves, following \cite{lozano2018,eschenazi2025solving}. This method constructs wave groups based on a list of compatible wave families (see \cite{Schecter1996, V.2015}), ensuring that each group satisfies the admissibility criterion derived from the existence of viscous profiles. Unlike the local continuation approach, this method enables the construction of complete and admissible wave curves even in degenerate or complex cases. For other wave curve construction methods, see \cite{mehrabi2020solution}.

To solve the Riemann problem \eqref{eq:Threephase_1}-\eqref{eq:RP}, the method involves constructing the forward $s$-wave curve starting from the left state $L$ together with the backward $f$-wave curve starting from the right state $R$. The intersection of these two curves determines the intermediate state $M$. Then, the solution to the Riemann problem is constructed in the forward direction from $L$ to $R$ and consists of the constant states $L$, $M$, and $R$ and the wave groups used to connect these states.
Notice that the only physically admissible saturation values that appear in the displacement are the constant states $L$, $M$, and $R$ and the saturation values along the rarefaction waves. 

\begin{figure}[ht]
    \centering
    \includegraphics[width=0.5\linewidth]{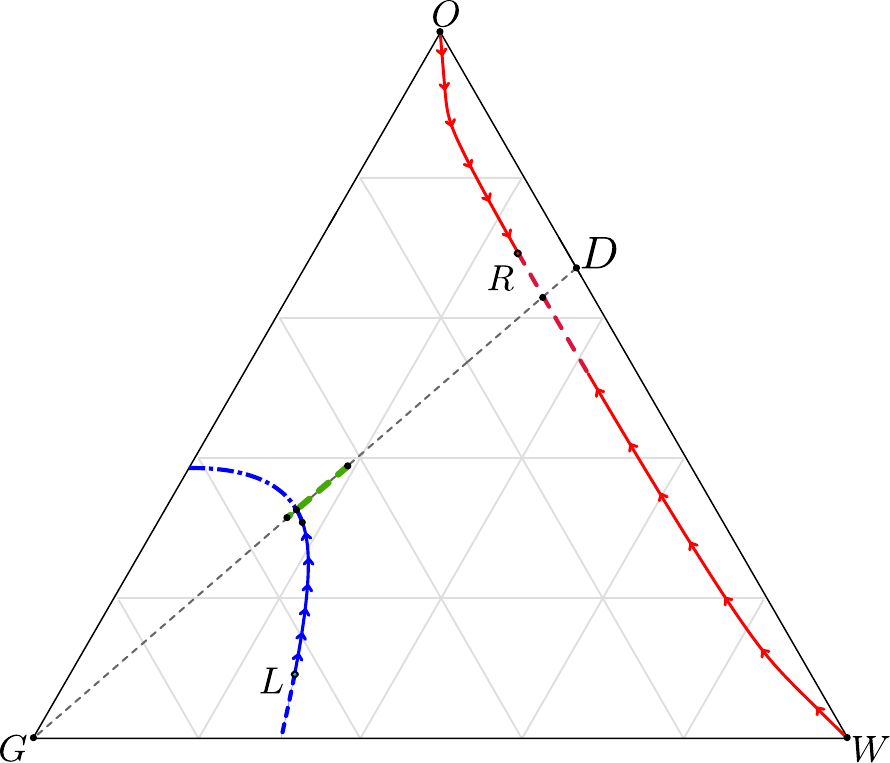}
    \caption{Graphical conventions for depicting wave curves in figures. Dashed curves indicate shocks. The dash-dotted curve corresponds to a composite curve. Continuous curves with arrows represent rarefaction curves. Arrows indicate the direction of increasing characteristic velocity. The colors blue, red, and green represent slow, fast, and undercompressive families.}
    \label{fig:Wave_curve_colors}
\end{figure}

We now present the graphical conventions for the wave curves inside $\Omega$, see Fig.~\ref{fig:Wave_curve_colors}. The components of fast wave curves are shown in red, while the components of slow wave curves are depicted in blue. Colored dashed curves indicate shocks, dash-dotted curves represent composite segments, and continuous curves represent rarefaction segments, with arrows indicating the direction of increasing characteristic velocity.

It has been verified numerically that shocks satisfy the viscous profile criterion. Additionally, only the parts of the wave curve that are used to solve Riemann problems are depicted in $\Omega$. Keeping this in mind, the right state may be isolated, indicating that it can only be reached by an $f$-shock. We plot in green admissible $u$-shock segments.  

\section{Riemann solutions}
\label{sec:Classification}
We aim to solve the Riemann problem with left states $L$ in the segment $[G, W]$ and right states $R$ in the region $\Omega$. A useful approach is to identify regions within the saturation triangle where the solutions exhibit similar properties. To achieve this, we partition  $\Omega$ into $ \mathcal{R} $-regions and the edge $[G, W]$ into $ \mathcal{L} $-segments.

Let us define $ \mathcal{R} $-regions and $ \mathcal{L} $-segments. Let $V$ be an open set of states in $\Omega$. We define $V$ as an $ \mathcal{R} $-region if, for every right state $R\in V$, the backward fast wave curves have the same structure, meaning it consists of the same wave sequence. Similarly, we define $ \mathcal{L}$-segment as a subset of $[G, W]$ where the slow wave curves possess analogous structures. Although $ \mathcal{L} $-segments and $ \mathcal{R} $-regions are constructed independently, their relationship is important to establish the existence and uniqueness of solutions of Riemann problems \cite{V.2015}.

Let us construct the $\mathcal{L}$-segments associated to each $\mathcal{R}$-region to identify Riemann problem solutions with the same structure.
Let $V$ be a backward $ \mathcal{R} $-region. Then for each $R \in V$, we construct the backward fast wave curve $ \mathcal{W}_f(R) $. Similarly, from each of the states $ M \in \mathcal{W}_f(R) $, we generate a backward slow wave curve, which continues until it reaches the edge $ [G, W] $. This edge is then subdivided into several $ \mathcal{L} $-segments. Notice that, by construction, each of the $ \mathcal{L} $-segments depends on $R$. We also construct the $s$-right-extension of $\mathcal{W}_f(R)$, which identifies the admissible $s$-shocks that can reach states over $\mathcal{W}_f(R)$, in order to ensure speed compatibility \eqref{eq:compatib}. Notice that in the procedure described above, certain parts of the edge $[G, W]$ may not be intersected. Consequently, the $s$-wave curves for these left states do not intersect $\mathcal{W}_f(R)$, yielding the Riemann problem with no classical solution. To circumvent this issue, non-classical waves must be introduced, as discussed below.  

In \cite{isaacson1992multiphase,L.1990,Andrade2018, lozano2018}, it was shown that there are two types of non-classical waves in the Riemann problem solutions for the three-phase model with quadratic exponents \eqref{eq:Threephase_0}-\eqref{eq:fi}: undercompressive shocks and transitional rarefactions. These waves appear when there is no direct intersection between the $f$- and $s$-wave curves. Since the viscosity matrix is the identity, these waves only occur along the invariant lines $[G, D]$, $[W, E]$, and $[O, B]$; see Fig.~\ref{fig:Macro_Classif_3foam}. As for our viscosity parameters, the umbilic point lies near $G$, and we are focusing on states $L \in [G, W]$, we describe only the non-classical waves along the invariant line $[G, D]$.


Following \cite{lozano2018, lozano2024c}, we begin the identification of non-classical waves by defining two key points along the segment $[G, D]$, namely $S$ and $D_0$. The point $S$ is the intersection of the $s$-inflection locus $\mathcal{I}_s$ with the segment $[G, D]$, while $D_0$ is an extension state of the umbilic point that satisfies condition \eqref{eq:ext_um}; see Figs.~\ref{fig:fig2} and \ref{fig:Macro_Classif_3foam}, as well as Definition \ref{def:extension_umbilic}. For any state $R \in \Omega$, we construct $\mathcal{W}_f(R)$  and define the point $M =\mathcal{W}_f(R) \cap (\mathcal{U}, D]$. The relative position of $M$ with respect to $S$ and $D_0$ determines the wave structure:
\begin{itemize}
    \item If $M$ lies on $(S, D_0]$, there exists $\widehat{M} \in [\mathcal{U}, S)$ such that any states $N$ in $[G, \mathcal{U})$ connect to $M$ through the wave sequence:
\begin{equation}\label{eq:transRarefactionM3}
N \xrightarrow{R_f} \mathcal{U} \xrightarrow{R_s} \widehat{M}\xrightarrow{S_s} M,
\end{equation}
which consists of a transitional rarefaction $N \xrightarrow{R_f} \mathcal{U} \xrightarrow{R_s}\widehat{
M}$ together with an $s$-shock $\widehat{M}\xrightarrow{S_s} M$.
\item If $ M$ lies on $(D_0, D] $, there exits a segment $[\mathcal{F}, \mathfrak{S}]\subset(G, \mathcal{U}) $ such that any state $ N$ in this segment is connected to $M $ by an undercompressive shock:
\begin{equation}\label{eq:ushock}
    {N}\xrightarrow{S_u} M.
\end{equation}
\end{itemize}
With this procedure, it is possible to explicitly identify $ \mathcal{L}$-segments and $\mathcal{R}$-regions where the solution of the Riemann problem involves each type of non-classical wave. This approach enables the identification of Riemann solutions with the same overall structure. For other works that follow this methodology without non-classical waves, see \cite{Azevedo2014, andrade2016oil, lozano2024b}.

\begin{rem}\label{rem:expl_M}
The state $M$ defines the undercompressive segment $[\mathcal{F}, \mathfrak{S}] \subset (G, \mathcal{U})$, where all $u$-shocks \eqref{eq:ushock} are admissible according to the viscous profile criterion. However, for left states $L \in [G, W]$, only a subset of this segment satisfies the velocity compatibility condition \eqref{eq:compatib} among the $s$-, $u$-, and $f$-shocks (see \cite{Andrade2018,lozano2018,lozano2024c} for more details). To identify this subset, it is necessary to compute the $s$-right extension of $[G, \mathcal{U}] \subset [G, D]$ and the Hugoniot locus $\mathcal{H}(M)$ in order to find the states $Z^* \in [G, \mathcal{U}]$ and $T_{\mathfrak{S}} \in \mathcal{H}(M)$ that satisfy $\sigma(Z^*; M) = \sigma(T_{\mathfrak{S}}; M) = \sigma(Z^*; T_{\mathfrak{S}})$ and $\sigma(T_{\mathfrak{S}}; M) = \lambda_s(T_{\mathfrak{S}})$. Here, the state $Z^*$ represents one of the endpoints of the segment contained in $[\mathcal{F}, \mathfrak{S}]$, while the other endpoint depends on the right state~$R$.
\end{rem}

\subsection{Subdivision in \texorpdfstring{$\mathcal{R}$}{R}-regions}\label{subsec2}
\begin{figure}[ht]
    \centering
    \includegraphics[width=0.7\linewidth]{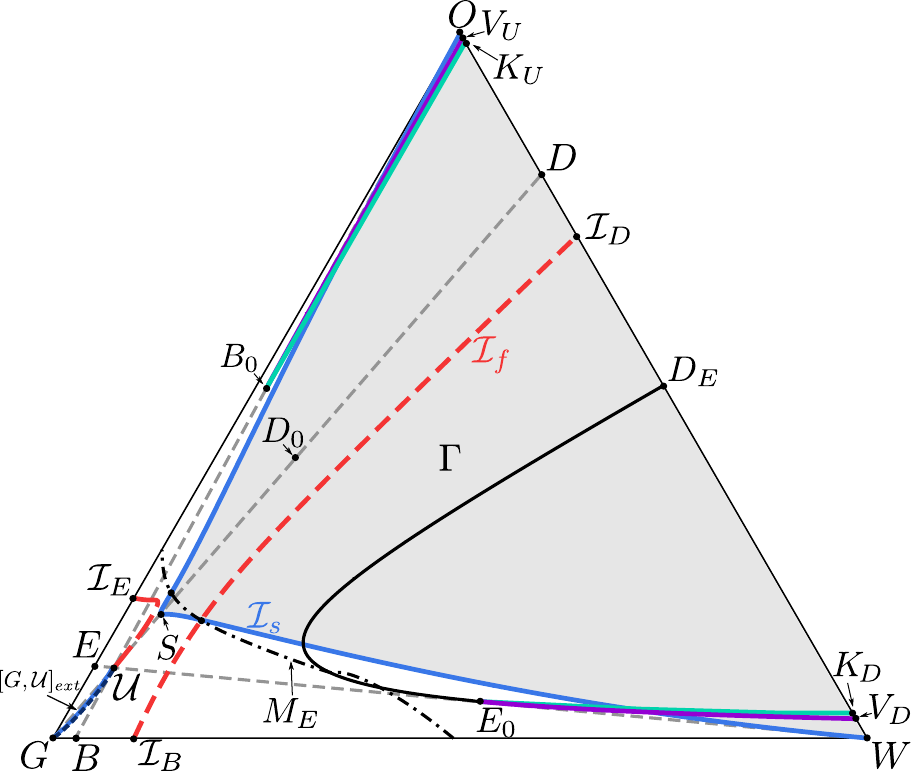}
    \caption{Bifurcation loci in the  $\Omega$ for the viscosity parameters presented in Table~\ref{tab:table1}. The gray region $\Gamma$ is bounded by the $s$-inflection locus $\mathcal{I}_s$ and the edge $[W, O]$. In this work, we consider right states $R$ within the region $\Gamma$. The curve $[G, \mathcal{U}]_{\text{ext}}$ (dashed dark blue) is the $s$-right extension of the segment $[G, \mathcal{U}]$.}
    \label{fig:Macro_Classif_3foam}
\end{figure}
As in \cite{lozano2024b}, we begin by introducing the bifurcation loci within  $\Omega$ that influence significantly the structure of the backward $f$-wave curve for right states $R$ (see Fig.~\ref{fig:Macro_Classif_3foam} for a visual representation). The bifurcation loci were computed numerically using the computer-assisted package ELI \cite{ELI_web}. Numerical experiments indicate that these bifurcation loci maintain their shapes and relative positions under variations in viscosity parameters that satisfy the inequalities \eqref{eq:desi_param_foam}.

Figure~\ref{fig:Macro_Classif_3foam} shows the bifurcation loci. The red curves $[\mathcal{I}_B, \mathcal{I}_D]$ and $[\mathcal{U}, \mathcal{I}_E]$  represent the $f$-inflection loci, while the blue curves $[G, \mathcal{U}]$ and $[W, S]\cup[S, O]$ denote the $s$-inflection loci. The segments $[G, D]$, $[W, E]$, and $[O, B]$ (gray dashed lines) form the secondary bifurcation loci. Along these loci and the boundaries of $\Omega$, system \eqref{eq:Threephase_1} reduces to the scalar Buckley-Leverett equation, see \cite{Azevedo2014}. 

\begin{figure}[]
    \centering
   \subfigure[]{ \includegraphics[width=0.57\linewidth]{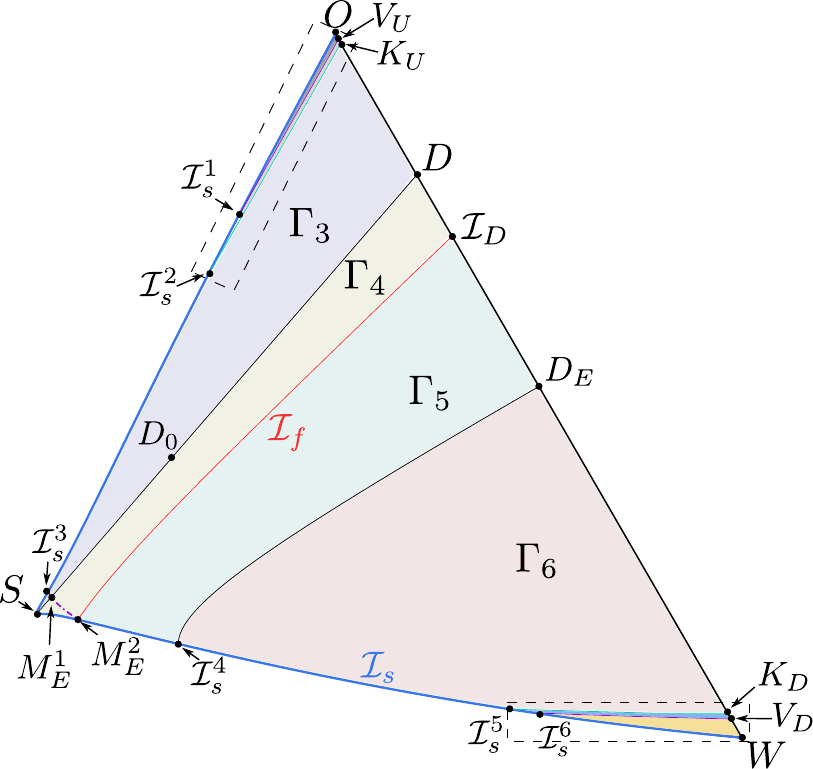}}
   \hspace{10pt}
   \subfigure{ \includegraphics[width=0.35\linewidth]{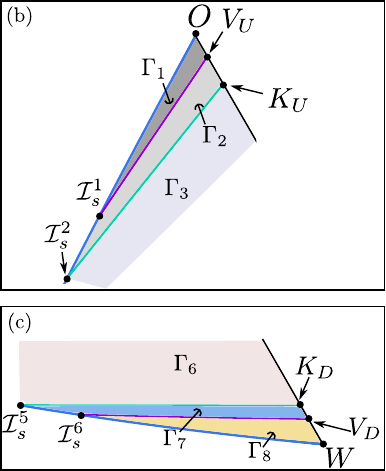}}
    \caption{
Subdivision of the region $\Gamma$ into eight $ \mathcal{R} $-regions $\Gamma_i$, where $i \in \{1, 2, \dots, 8\}$. (a)   Zoom of the region $\Gamma$.
(b) Zoom of the region close to the corner $O$.  
(c) Zoom of the region close to the corner $W$.  }
    \label{fig:Subdivision_Gamma}
\end{figure}

See Definition \ref{def:iLRExt}. The following loci represent extensions of the secondary bifurcation and inflection loci, which identify locations where significant changes occur in the backward $f$-wave curves:  
\begin{itemize}  
    \item The solid black curve $[D_E, E_0]$ in Fig.~\ref{fig:Macro_Classif_3foam} is the $f$-left extension of the segment $[D, \mathcal{U}]$. The state $E_0\in[W, E]$ is an extension state of the umbilic point that satisfies \eqref{eq:ext_um}.  
    \item The green curves $[K_D, E_0]$ and $[K_U, B_0]$ in Fig.~\ref{fig:Macro_Classif_3foam} are the $f$-left extensions of the segment $[G, \mathcal{U}]$. The state $B_0\in[O, B]$ is an extension state of the umbilic point that satisfies \eqref{eq:ext_um}.  
    \item The purple curves $[V_D, E_0]$ and $[V_U, B_0]$ in Fig.~\ref{fig:Macro_Classif_3foam} are double extension curves of the segment $[G, \mathcal{U}]$. This extension process involves first computing the segment $[G, \mathcal{U}]_{\text{ext}}$, the $s$-right extension of $[G, \mathcal{U}]$, followed by its $s$-left extension.  
    \item The dashed black curve $M_E$ in Fig.~\ref{fig:Macro_Classif_3foam} is the $f$-left-extension of the $s$-inflection $\mathcal{I}_s$. 
\end{itemize}
Figure~\ref{fig:Macro_Classif_3foam} highlights the region $\Gamma$ in gray, bounded by the edge $[W, O]$ and the $s$-inflection locus  $\mathcal{I}_s$, given $[O, S]\cup[S, W]$ (blue curve); see also Fig.~\ref{fig:Subdivision_Gamma}. In this work, we focus on the right states $R$ within the subregion $\Gamma \subset \Omega$ for the following reasons:  
\begin{itemize}
    \item Consider viscosity parameters $\mu_w, \mu_o, \widehat{\mu_g}$ satisfying \eqref{eq:desi_param_foam}. Numerical experiments indicate that as the parameter $MRF$ increases, the umbilic point $\mathcal{U}$ moves closer to the corner $G$ along the straight line $[G, D]$, and simultaneously, the state $S$ approaches $G$. During this process, the region $\Gamma$ bounded by the $s$-inflection locus expands to cover nearly the whole saturation triangle $\Omega$. The relative positions of the bifurcation loci observed earlier are preserved.  
    \item Analyzing the Riemann problem solutions for right states $R \in \Gamma$ thus allows to cover a significant portion of  $\Omega$. Moreover, since our study is motivated by applications in oil recovery, the region $\Gamma$ is particularly relevant, as it corresponds to the case with a significant presence of initial oil.  
\end{itemize}
\begin{rem}
The umbilic point in the present work is close to the corner $G$ (here, in the context of EOR applications with the injection of a mixture of foamed gas and water) yielding portions of the backward wave curves with role in the Riemann problem solution significantly different from those in \cite{lozano2024b}.
\end{rem}

\subsection{Subdivision of the \texorpdfstring{$\Gamma$}{T}-region}
Figure \ref{fig:Subdivision_Gamma} shows the subdivision of the region $\Gamma$ into eight $ \mathcal{R} $-regions, denoted as $\Gamma_i$, $i \in \{1, 2, \dots, 8\}$. Figures \ref{fig:Subdivision_Gamma}(b) and \ref{fig:Subdivision_Gamma}(c) provide zoomed-in views near the corners $O$ and $W$, respectively.   

We use the notation $[A, B]_{\text{ext}}$ to indicate the $s$-right extension of a backward $f$-wave curve or a segment of a secondary bifurcation. These extension curves are plotted in black.

\subsubsection{Riemann problem solution for \texorpdfstring{$R\in \Gamma_1$}{RET1} }\label{sec:Gamma1}
\begin{figure}[ht]
      \centering
        \subfigure[Structure of the Riemann solution for a generic $R \in \Gamma_1$, with $R = (0.00628771,0.916593)$.] 
      {\includegraphics[width=0.465\linewidth]{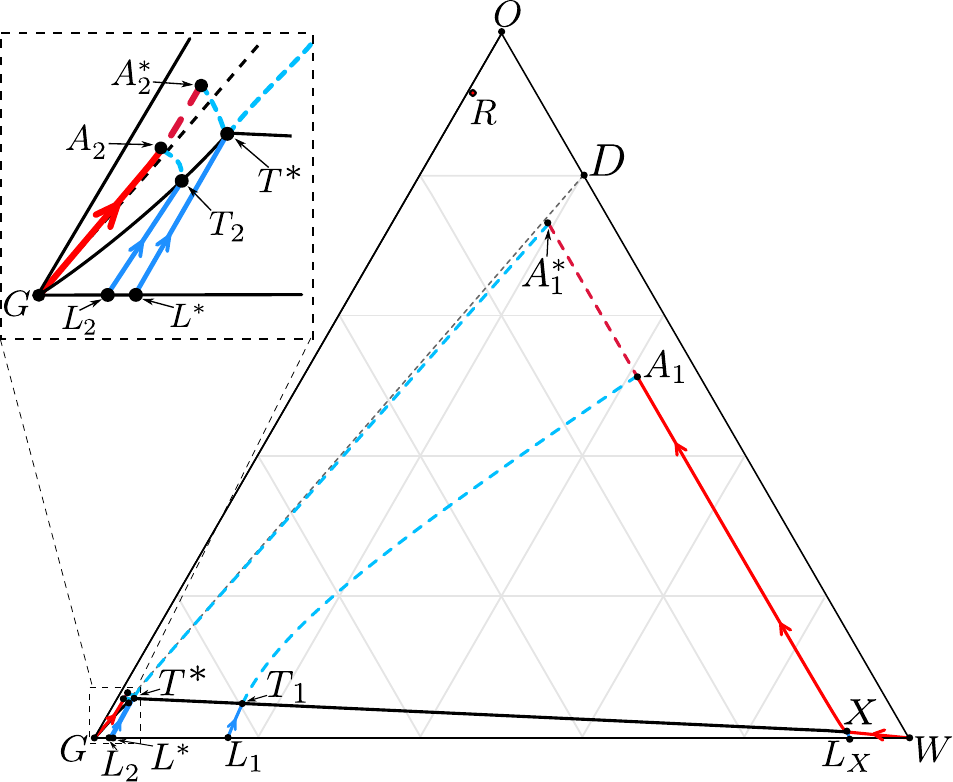}}\hspace{15pt}
       \subfigure[Structure of the Riemann solution for a generic $R \in \Gamma_2$, with $R = (0.00865231, 0.793866)$.]
      {\includegraphics[width=0.465\linewidth]{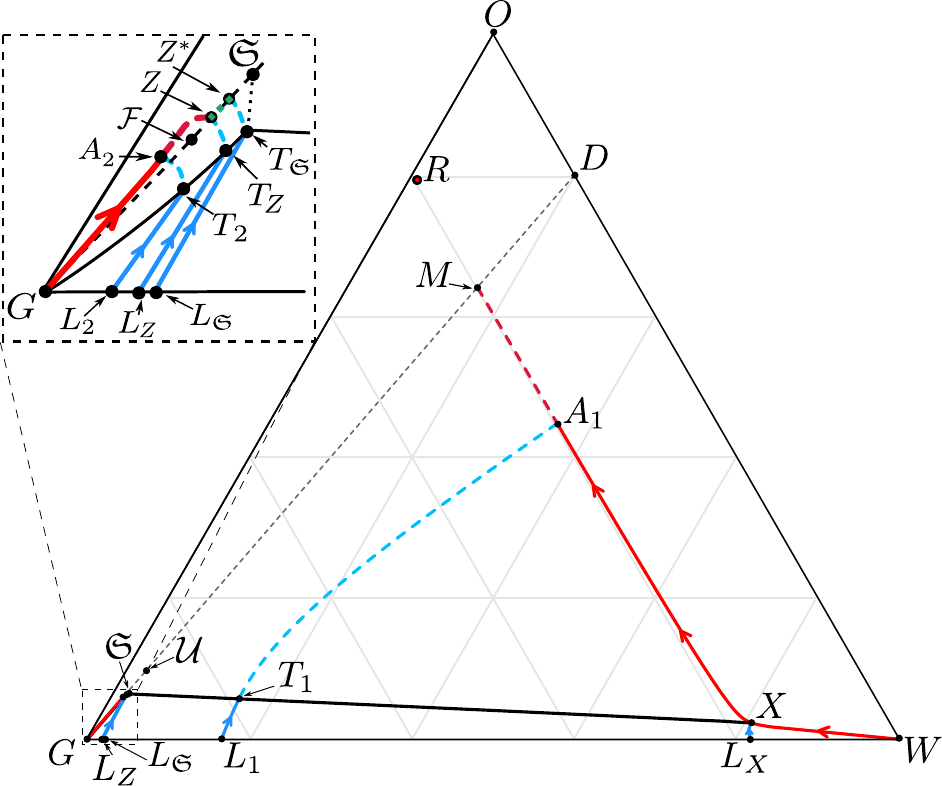}}\hspace{15pt}
      \caption{The blue dashed curves (respectively, red) represent $s$-shock curves (respectively, $f$-shock curves). The blue continuous curve (respectively, red) represents the $s$-rarefaction curves (respectively, the $f$-rarefaction curves). The arrows indicate increasing characteristic velocity. (a) The black curves $[X, T^*]$ and $[G, T^*]$ represent the $s$-right extensions of $\mathcal{W}_f(R)$. (b) The black curves $[X, T_1]_{\text{ext}} \cup [T_1, T_{\mathfrak{S}}]_{\text{ext}} \cup [G, T_Z]_{\text{ext}}$ represent the $s$-right extension of $\mathcal{W}_f(R)$. The green dotted segment $(Z, Z^*] \subset [G, D]$ identifies admissible $u$-shocks with $M$. The black curve $(T_Z, T_{\mathfrak{S}}]_{\text{ext}}$ is the $s$-extension of $(Z, Z^*]$.}
      \label{fig:Gamma_1_e_2}
  \end{figure}

This region corresponds to the dark gray area in Fig.~\ref{fig:Subdivision_Gamma}(b).
The region $\Gamma_1$ is bounded by:
\begin{itemize}  
     \item The blue curve $[O, \mathcal{I}_s^1]$, which is part of the $s$-inflection locus;  
    \item The line segment $[O, V_U]$, which belongs to the boundary of $\Omega$; and  
    \item The purple curve $[V_U, \mathcal{I}_s^1]$, which is part of the double extension curve of $[G, \mathcal{U}]$.  
\end{itemize}  
The purple curve $[V_U, \mathcal{I}_s^1]$ marks the boundary where $u$-shock compatibility changes (see Claim 4.10 in \cite{Andrade2018} and \cite{lozano2018} for more details). For $L \in [G, W]$ and $R \in \Gamma_1$, $ u$-shock velocities are incompatible with the $s$- and $ f$-shock velocities, and thus the Riemann problem solution consists solely of classical wave groups. 
In this case, $\mathcal{W}_f(R)$ consists of a local branch containing the state $R$ and a nonlocal branch containing $G$. Figure \ref{fig:Gamma_1_e_2}(a) presents the components of $\mathcal{W}_f(R)$ used to construct the solution of the Riemann problem for $L \in [G, W]$: 
\begin{itemize}
    \item The state $R$;
    \item  The $f$-rarefaction curves $[W, A_1)$ and $[G, A_2)$ (Fig.~\ref{fig:Gamma_1_e_2}(a)), with $\lambda_f(A_1) = \sigma(A_1; R)$ and $\lambda_f(A_2) = \sigma(A_2; R)$; and
    \item The $f$-shock curves $[A_1, A_1^*]$ and $[A_2, A_2^*]$ (Fig.~\ref{fig:Gamma_1_e_2}(a)), with $\sigma(A_1^*; R)= \sigma(A_2^*; R).$
\end{itemize}  The state $X$ is the intersection of $\mathcal{W}_f(R)$ and the $s$-inflection locus, with $X \in [W, A_1)$. 
 The segments $(G, T_2)_{\text{ext}}$, $[T_2, T^*]_{\text{ext}}$, $(X, T_1)_{\text{ext}}$, and $[T_1, T^*)_{\text{ext}}$ 
denote the  $s$-right-extension of the backward $f$-wave curve segments $(G, A_2)$, $[A_2, A_2^*]$, $(X, A_1)$, and $[A_1, A_1^*)$ respectively.

Let $L_1, L_2, L^*$, and $L_X$ denote the intersection points of the backward $s$-wave curves through the states $A_1, A_2, A_1^*$, and $X$ with the edge $[G, W]$ (see Fig.~\ref{fig:Gamma_1_e_2}(a)). Then, the structure of the solution for the Riemann problem for $L \in [G, W]$ and $R \in \Gamma_1$ (see Fig.~\ref{fig:Subdivision_Gamma}(b)) consists of:

\begin{enumerate}
         \item[(i)] $L\testright{R_f} A_1 \xrightarrow{S_f} R$, for $L=W$;
       
         \item[(ii)] $ L\testright{R_s} Y  \testright{R_f} A_1 \xrightarrow{S_f} R,$ for $L\in(W, L_{X}]$ and
          $Y\in(W,X]\subset \mathcal{W}_f(R)$;
 \item[(iii)]$L\testright{R_s} T  \xrightarrow{S_s} Y \testright{R_f} A_1 \xrightarrow{S_f} R$, for $L\in(L_{X},L_{1})$, $T\in(X, T_1)_{\text{ext}}$ and $Y\in(X,A_{1})\subset \mathcal{W}_f(R)$;
\item[(iv)] $L\testright{R_s} T  \xrightarrow{S_s} Y  \testright{S_f} R$, for $L\in[L_1,L^*)$, $T\in[T_1, T^*)_{\text{ext}}$ and $Y\in[A_1, A_{1}^*)\subset \mathcal{W}_f(R)$;
 \item[(v)] $L\testright{R_f} A_2  \xrightarrow{S_f} R$, for $L = G$;

\item[(vi)] $L\testright{R_s} T  \xrightarrow{S_s} Y \testright{R_f} A_2 \xrightarrow{S_f} R$, for $L\in(G,L_{2})$, $T\in(G, T_2)_{\text{ext}}$ and $Y\in(G,A_{2})\subset \mathcal{W}_f(R)$;
 \item[(vii)] $L\testright{R_s} T  \xrightarrow{S_s} Y  \testright{S_f} R$, for $L\in[L_{2},L^*)$, $T\in[T_{2},T^*)_{\text{ext}}$ and $Y\in[A_{2}, A_2^*)\subset \mathcal{W}_f(R)$.
     \end{enumerate}
     \begin{rem}
       According to the triple shock rule \cite{Azevedo2014}, $(\sigma(T^*; R) =  \sigma(T^*; A_j^* ) = \sigma(A_j^*; R), \, j = 1, 2)$, the solution for $L = L^*$ is represented by three possible wave sequences in state space, one of which does not satisfy the viscous profile criterion. The admissible sequences are
       $$L^*\testright{R_s} T^*  \xrightarrow{S_s} A_1^*\testright{S_f}R\,\,\mbox{ and }\,\, L^*\testright{R_s} T^*  \xrightarrow{S_s} A_2^*\testright{S_f} R.$$
       These sequences correspond to a unique solution in $xt$-space; see \cite{andrade2016oil, Azevedo2014, fritis2024riemann}.
     \end{rem}

\subsubsection{Riemann problem solution for \texorpdfstring{$R\in \Gamma_2$}{RET2}}\label{sec:Gamma2}

This region corresponds to the light gray area in Fig.~\ref{fig:Subdivision_Gamma}(b). The regions $ \Gamma_2 $ and $ \Gamma_1 $ are separated by the purple curve where the compatibility with the $ u $-shocks changes, $ [V_U, \mathcal{I}_s^1] $. The other boundaries of $ \Gamma_2 $ are:  
\begin{itemize}
    \item The line segment $[V_U, K_U]$, which belongs to the boundary of $\Omega$;
    \item The blue curve $[\mathcal{I}_s^1, \mathcal{I}_s^2]$, which is part of $\mathcal{I}_s=[W, S]\cup[S, O]$; and
    \item The green curve $[K_U, \mathcal{I}_s^2]$, which marks the boundary where the admissibility of nonlocal $f$-shocks changes. It corresponds to the $f$-left-extension of $ [G, \mathcal{U}] $ (see \cite{Andrade2018, lozano2018}). 
\end{itemize}
Thus, for $L \in [G, W]$ and $R \in \Gamma_2$, the solution to the Riemann problem comprises both classical and non-classical wave groups. For $R \in \Gamma_2$, $\mathcal{W}_f(R)$ consists of a local branch containing $R$ and a nonlocal branch containing $G$. Figure \ref{fig:Gamma_1_e_2}(b) illustrates the components of $\mathcal{W}_f(R)$ used to construct the solution of the Riemann problem for $L \in [G, W]$:
\begin{itemize}
    \item The state $R$;
    \item The $f$-shock curves $[A_1, M]$ and $[A_2, Z]$ in Fig.~\ref{fig:Gamma_1_e_2}(b); and
    \item The $f$-rarefaction curves $[W, A_1)$ and $[G, A_2)$ in Fig.~\ref{fig:Gamma_1_e_2}(b), with $\lambda_f(A_1) = \sigma(A_1; R)$ and $\lambda_f(A_2) = \sigma(A_2; R)$.
\end{itemize}
The state $X \in [W, A_1)$ corresponds to the intersection of $\mathcal{W}_f(R)$ with the $s$-inflection locus $\mathcal{I}_s$. The state $M$ is defined by the intersection of the local branch of $ \mathcal{W}_f(R) $ and the segment $[G, D]$. Both $Z$ and $M$ lie on $[G, D]$, where $M$ defines the undercompressive shock segment $(Z, Z^*] \subset [\mathcal{F}, \mathfrak{S}]$. This segment satisfies velocity compatibility condition \eqref{eq:compatib} for $s$-, $u$-, and $f$-shocks. The state $Z^* \in \mathcal{H}(M)$ is associated with $T_{\mathfrak{S}} \in \mathcal{H}(M)$, where both satisfy the following conditions (see Remark \ref{rem:expl_M}):
\begin{equation}\label{eq:zetaestrella}
\sigma(Z^*; M) = \sigma(T_{\mathfrak{S}}; M) = \sigma(Z^*; T_{\mathfrak{S}}),\quad \mbox{and}\quad \sigma(T_{\mathfrak{S}}; M) = \lambda_s(T_{\mathfrak{S}}).
\end{equation}

The segments $(G, T_2)_{\text{ext}}$, $[T_2, T_Z]_{\text{ext}}$, $(X, T_1)_{\text{ext}}$, and $[T_1, T_{\mathfrak{S}})_{\text{ext}}$ denote the $s$-right-extensions of the backward $f$-wave curve segments $(G, A_2)$, $[A_2, Z]$, $(X, A_1)$, and $[A_1, M)$. The segment $(T_Z, T_{\mathfrak{S}}]_{\text{ext}}$ denotes the  $s$-right-extension of $(Z, Z^*] \subset [G, \mathcal{U}]$.

Let $L_1, L_2, L_Z, L_{\mathfrak{S}}$, and $L_X$ denote the intersection points of the backward $s$-wave curves through states $A_1$, $A_2$, $Z$, $Z^*$, and $X$ with the edge $[G, W]$ (see Fig.~\ref{fig:Gamma_1_e_2}(b)). Then, the structure of the solution for the Riemann problem for $L \in [G, W]$ and $R \in \Gamma_2$ (see Fig.~\ref{fig:Subdivision_Gamma}(b)) consists of:
\begin{enumerate}
         \item[(i)] $L\testright{R_f} A_1 \xrightarrow{S_f} R$, for $L = W$;
         \item[(ii)]$ L\testright{R_s} Y  \testright{R_f} A_1 \xrightarrow{S_f} R$, for $L\in(W,L_{X}]$ and $Y\in(W,X]\subset \mathcal{W}_f(R)$;
 \item[(iii)] $ L\testright{R_s} T  \xrightarrow{S_s} Y \testright{R_f} A_1 \xrightarrow{S_f} R$, for $L\in(L_{X},L_{1})$, $T\in(X, T_{1})_{\text{ext}}$ and  $Y\in(X,A_{1})\subset \mathcal{W}_f(R)$;
 \item[(iv)] $L\testright{R_s} T  \xrightarrow{S_s} Y  \testright{S_f} R$, for $L\in[L_{1}, L_{\mathfrak{S}})$, $T\in[T_{1}, T_{\mathfrak{S}},)_{\text{ext}}$ and  $Y\in[A_{1},M)\subset \mathcal{W}_f(R)$;
 \item[(v)] $L\testright{R_f} A_2  \xrightarrow{S_f} R$, for $L = G$;
 \item[(vi)] $L\testright{R_s} T  \xrightarrow{S_s} Y \testright{R_f} A_2 \xrightarrow{S_f} R$, for $L\in(G,L_{2})$, $T\in(G, T_2)_{\text{ext}}$ and $Y\in (G,A_2)\subset \mathcal{W}_f(R)$; 
 \item[(vii)] $L\testright{R_s} T  \xrightarrow{S_s} Y  \testright{S_f} R$, for $L\in[L_{2},L_Z]$, $T\in[T_{2},T_Z]_{\text{ext}}$ and $Y\in[A_{2},Z]\subset \mathcal{W}_f(R)$.
\item[(viii)] $L\testright{R_s} T  \xrightarrow{S_s} N  \testright{S_u} M \testright{S_f} R$, $L\in(L_{Z},L_{\mathfrak{S}}]$, $T\in(T_{Z},T_{\mathfrak{S}}]_{\text{ext}}$ and $N\in(Z,Z^*]$.
      \end{enumerate}
  

\subsubsection{Riemann problem solution for \texorpdfstring{$R\in \Gamma_3$}{RET3}}
\label{sec:Gamma3}

This region corresponds to the lavender area in Figs.~\ref{fig:Subdivision_Gamma}(a)-(b). The region $ \Gamma_3 $ is adjacent to $ \Gamma_2 $, they are separated by the green curve $[K_U, \mathcal{I}_s^2]$, where admissibility of nonlocal $f$-shocks changes; see Remark 4.12 in \cite{Andrade2018} and \cite{lozano2018}. The remaining boundaries of $ \Gamma_3 $ are:  
\begin{itemize}
    \item The blue curve $[\mathcal{I}_s^2, S]$, which belongs to the $s$-inflection locus; 
    \item The line segment $[S, D]$, and 
    \item The line segment $[D, K_U]$, which is part of the boundary of $\Omega$.
\end{itemize}
 Therefore, for $L \in [G, W]$ and $R \in \Gamma_3$, the Riemann problem solution involves both classical and non-classical wave groups.

\begin{figure}[ht]
	\centering
     \subfigure[Structure of the Riemann solution for a generic $R \in \Gamma_3$ with $M \in {(S, D_0]}$, located to the left of the curve $M_E$, for $R = {(0.045, 0.2)}$.]%
	 {\includegraphics[width=0.45\linewidth]{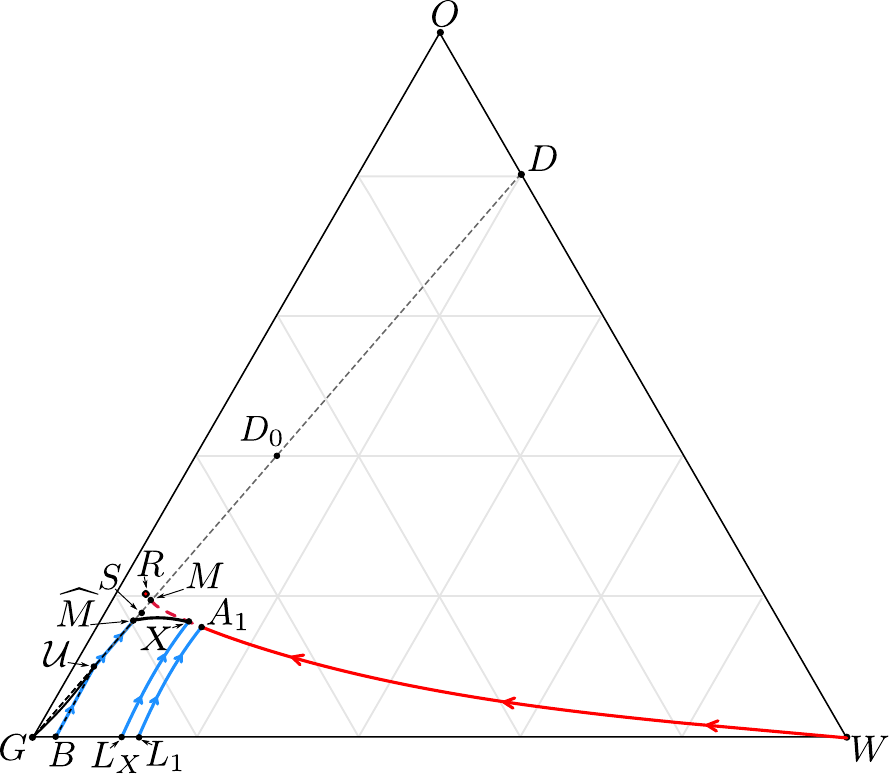}} 
	\hspace{0.8mm}
	\subfigure[Structure of the Riemann solution for a generic $R \in \Gamma_3$ with $M \in {(S, D_0]}$, located to the right of the curve $M_E$, for $R = (0.0540512, 0.322225)$.]
	 {\includegraphics[width=0.45\linewidth]{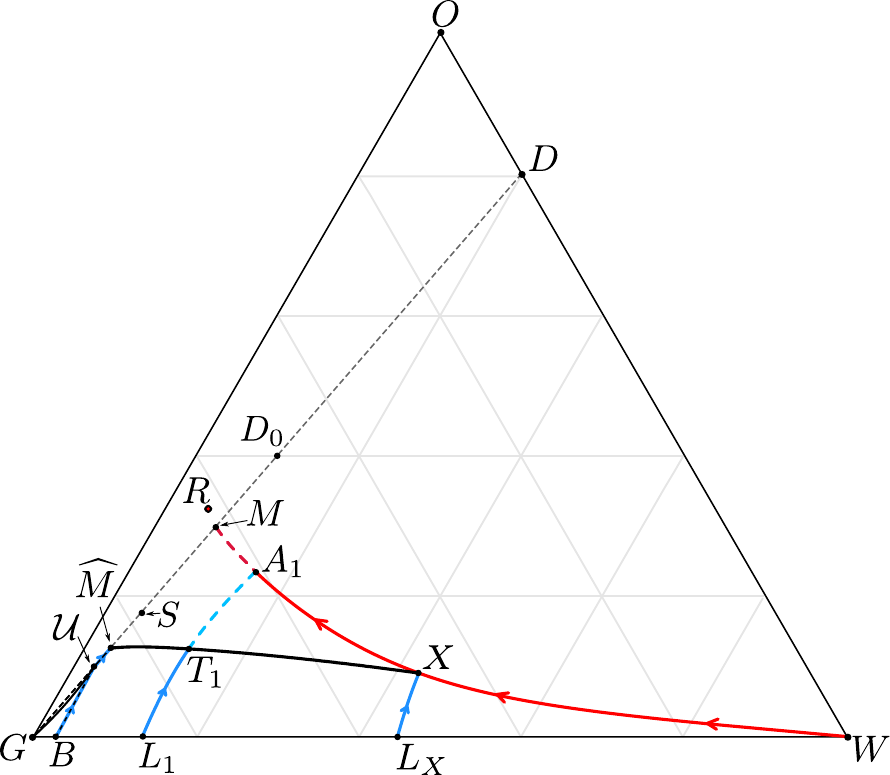}} 
	\hspace{0.8mm}
		\caption{ Structure of the Riemann solution for $R\in\Gamma_3$. The blue dashed curve (respectively red) represents $s$-shock curves (respectively $f$-shock curves). The blue continuous curve (respectively red) represents $s$-rarefaction curves (respectively the $f$-rarefaction curve). The arrows indicate the increasing characteristic velocity. The black curve $[X, \widehat{M}]_{\text{ext}}$ represents the $s$-right extension of $\mathcal{W}_f(R)$. The black curve $[G, \mathcal{U}]_{\text{ext}}$ represents the $s$-right extension of $[G,\mathcal{U}]$. 
	}
	\label{fig:Gamma3_Up_Down}
\end{figure}

\begin{figure}[ht]
	\centering
    \subfigure[Structure of the Riemann solution for generic $R\in\Gamma_3$ with $M\in{(D_0,D]}$ and $R = (0.116779, 0.831619)$. ] 
	{\includegraphics[width=0.45\linewidth]{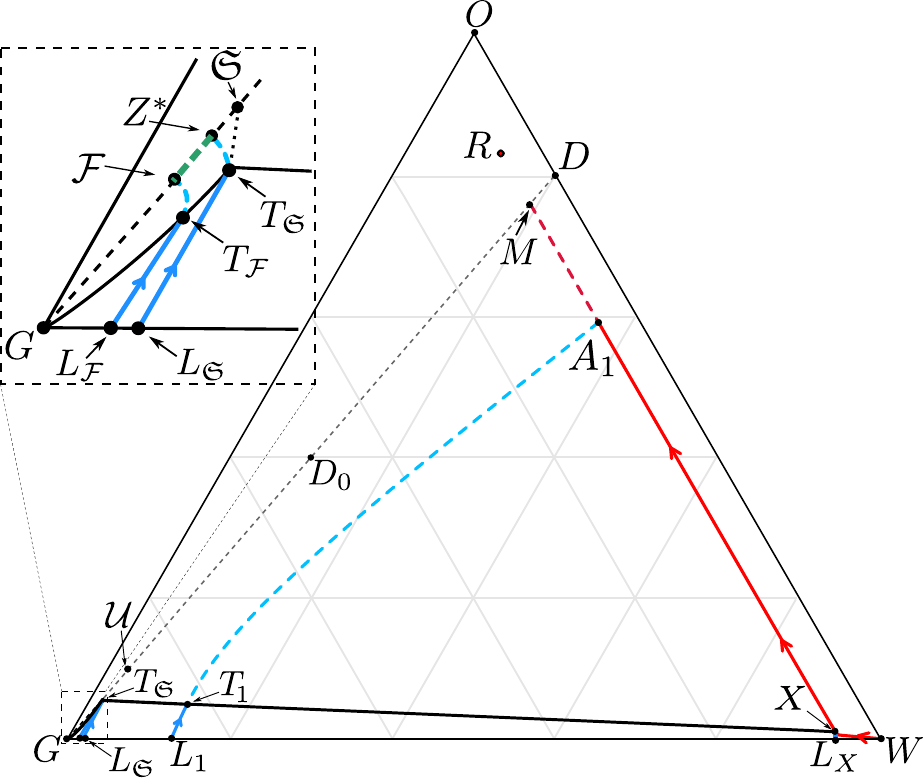}} 
		\caption{ The blue dashed curves (respectively, red) represent $s$-shock curves (respectively, $f$-shock curves). The blue continuous curve (respectively, red) represents the $s$-rarefaction curves (respectively, the $f$-rarefaction curves). The arrows indicate increasing characteristic velocity. The black curve  $[X, T_{\mathfrak{S}}]_{\text{ext}}$ represents the $s$-right extension of $\mathcal{W}_f(R)$. The green dotted segment $[\mathcal{F}, Z^*] \subset$ $[G, D]$ identifies admissible $u$-shocks with $M$.  The black curve $[G, T_{\mathfrak{S}}]_{\text{ext}}$ is the $s$-right extension of $[G, Z^*]$.
		}
	\label{fig:Gamma3_D0_D}
\end{figure}

According to \cite{lozano2018, lozano2024b, Andrade2018}, for $R \in \Gamma_3$, the wave curve $\mathcal{W}_f(R)$ consists exclusively of a local branch connecting the states $O$, $R$ and $W$. Figures \ref{fig:Gamma3_Up_Down} and \ref{fig:Gamma3_D0_D} illustrates the components of $\mathcal{W}_f(R)$ used to construct the Riemann problem solution for $L \in [G, W]$: 
\begin{itemize}
    \item The state $R$;
    \item The $f$-shock curve $[A_1, M]$; and
    \item The $f$-rarefaction curve $[W, A_1)$, with $\lambda_f(A_1) = \sigma(A_1; R)$.
\end{itemize}
 The state $X$ marks the intersection of $\mathcal{W}_f(R)$ and $\mathcal{I}_s$. Its position depends on the location of $R$ relative to the curve $M_E$ (see Figs.~\ref{fig:Macro_Classif_3foam} and \ref{fig:Subdivision_Gamma}):
 \begin{itemize}
      \item If $R$ is on the left of curve $M_E$, then $X$ belongs to the $f$-shock curve $[A_1, M]$, see Fig.~\ref{fig:Gamma3_Up_Down}(a);
     \item If $R$ is on the right of the curve $M_E$, then $X$ belongs to the $f$-rarefaction curve $[W, A_1)$, see Figs.~\ref{fig:Gamma3_Up_Down}(b) and \ref{fig:Gamma3_D0_D}.  
 \end{itemize}
Here, we focus on the case $X \in [W, A_1)$; the alternative case is discussed in Remark \ref{rem:X_ME_Gamma_3}.
The state $M$ is determined by the intersection of $ \mathcal{W}_f(R) $ and the segment $ [G, D] $, while $S$ denotes the intersection of $[G, D]$ and $ \mathcal{I}_s $, see Figs.~\ref{fig:Gamma3_Up_Down} and \ref{fig:Gamma3_D0_D}. Let $ D_0$ represent an extension of the umbilic point. As previously noted, the type of non-classical solution used to solve the Riemann problem depends on whether $ M \in (S, D_0] $ or $ M \in (D_0, D] $:
 \begin{itemize}
    \item For $ M \in (S, D_0] $, there exists a state $ \widehat{M} \in [\mathcal{U}, S) $ that enables $ M $ to be reached via the wave structure \eqref{eq:transRarefactionM3} (Fig.~\ref{fig:Gamma3_Up_Down});
    \item For $ M \in (D_0, D] $, $ M $ defines the undercompressive shock segment $[\mathcal{F}, Z^*] \subset [\mathcal{F}, \mathfrak{S}]$, which satisfies the velocity compatibility condition \eqref{eq:compatib} for $ s $-, $ u $-, and $ f $-shocks, see Fig.~\ref{fig:Gamma3_D0_D}. Here, $ Z^* \in \mathcal{H}(M) $ and $ T_{\mathfrak{S}} \in \mathcal{H}(M) $ both satisfy \eqref{eq:zetaestrella} (Remark \ref{rem:expl_M}).
\end{itemize}

The segment $(G, \mathcal{U})_{\text{ext}}$ represents the $s$-right extension of $(G, \mathcal{U})$ (see Fig.~\ref{fig:Gamma3_Up_Down}). Analogously, the segments $(G, T_{\mathcal{F}}]_{\text{ext}} \cup (T_{\mathcal{F}}, T_{\mathfrak{S}}]_{\text{ext}}$ correspond to the $s$-right extension of $(G, \mathcal{F}] \cup (\mathcal{F}, Z^*] \subset [G, D]$ (see Fig.~\ref{fig:Gamma3_D0_D}).

The extension of $\mathcal{W}_f(R)$ for $R \in \Gamma_3$ is given by:
\begin{itemize}
    \item The segment $[X,\widehat{M}]_{\text{ext}}$, if $R$ lies to the left of the curve $M_E$ and $M \in (S, D_0]$ (Fig.~\ref{fig:Gamma3_Up_Down}(a)).  
    \item The segments $[X,T_1]_{\text{ext}} \cup [T_1, \widehat{M}]_{\text{ext}}$, if $R$ lies to the right of the curve $M_E$ and $M \in (S, D_0]$ (Fig.~\ref{fig:Gamma3_Up_Down}(b)).    
    \item The segments $[X,T_1]_{\text{ext}} \cup [T_1, T_{\mathfrak{S}})_{\text{ext}}$, if $R$ lies to the right of the curve $M_E$ and $M \in (D_0, D]$ (Fig.~\ref{fig:Gamma3_D0_D}).  
\end{itemize}

Let $ L_1 $, $ L_{\mathcal{F}} $, $ L_{\mathfrak{S}} $, and $ L_X $ denote the intersection points of the backward $ s $-wave curves originating from states $ A_1 $, $ \mathcal{F} $, $ Z^* $, and $ X $, respectively, with the edge $ [G, W] $ (Figs.~\ref{fig:Gamma3_Up_Down} and \ref{fig:Gamma3_D0_D}). Then, the structure of the solution for the Riemann problem for $L \in [G, W]$ and $R \in \Gamma_3$ (see Fig.~\ref{fig:Subdivision_Gamma}(a)) consists of:
\begin{enumerate}
     \item[(i)] $L\testright{R_f} A_1 \xrightarrow{S_f} R$, for $L = W$;
     \item[(ii)] $L\testright{R_s}  Y \testright{R_f} A_1 \xrightarrow{S_f} R$, for $L\in(W, L_{X}]$ and $Y\in(W, X]\subset \mathcal{W}_f(R)$;
     \item[(iii)] $L\testright{R_s} T \xrightarrow{S_s} Y \testright{R_f} A_1 \xrightarrow{S_f} R$, for $L\in(L_X,L_{1})$, $T\in(X, T_{1})_{\text{ext}}$ and $Y\in(X, A_1)\subset \mathcal{W}_f(R)$.
     \end{enumerate}
\begin{itemize}
       \item[$\bullet$] Case $M\in (S, D_0]$, see Fig.~\ref{fig:Gamma3_Up_Down}(b):
    \begin{enumerate}
          \item[(iv)] $L\testright{R_f} \mathcal{U}\xrightarrow{R_s} \widehat{M}  \xrightarrow{S_s} M \testright{S_f} R$, for $L = G$;
          \item[(v)] $L\testright{R_s} T \xrightarrow{S_s} N \testright{R_f}\mathcal{U}\xrightarrow{R_s} \widehat{M}  \xrightarrow{S_s} M \testright{S_f} R$, for $L\in(G, B)$, $T\in(G, \mathcal{U})_{\text{ext}}$ and $N\in(G, \mathcal{U})\subset[G, D]$;
           \item[(vi)] $L\testright{R_s} T \xrightarrow{S_s} Y \testright{S_f} R$, for $L\in[L_{1}, B]$, $T\in[T_{1}, \widehat{M}]_{\text{ext}}$ and $Y\in[A_1, M]\subset \mathcal{W}_f(R)$.
      \end{enumerate}
       \item[$\bullet$] Case $M\in (D_0, D]$, see Fig.~\ref{fig:Gamma3_D0_D}:
    \begin{enumerate}
          \item[(iv)] $L\testright{R_f} \mathcal{F}\xrightarrow{S_u} M \testright{S_f} R$, for $L = G$;
           \item[(v)] $L\testright{R_s} T \xrightarrow{S_s} N\testright{R_f} \mathcal{F}  \xrightarrow{S_u} M \testright{S_f} R$, for $L\in(G,L_{\mathcal{F}})$, $T\in(G,T_{\mathcal{F}})_{\text{ext}}$ and $N\in(G,\mathcal{F})\subset [G,D]$;
          \item[(vi)] $L\testright{R_s} T \xrightarrow{S_s} N  \testright{S_u} M \testright{S_f} R$, for $L\in[L_{\mathcal{F}},L_{\mathfrak{S}}]$, $T\in[T_{\mathcal{F}},T_{\mathfrak{S}}]_{\text{ext}}$ and $N\in[\mathcal{F},Z^*]\subset [G,D]$;
          \item[(vii)] $L\testright{R_s} T \xrightarrow{S_s} Y \testright{S_f} R$, for $L\in[L_1, L_{\mathfrak{S}})$, $T\in[T_1,T_{\mathfrak{S}})_{\text{ext}}$ and $Y\in[A_1, M)\subset \mathcal{W}_f(R)$.
     \end{enumerate}
\end{itemize}
\begin{rem}\label{rem:X_ME_Gamma_3}
Consider the case where $ M \in (S, D_0] $ and $ X $ lies on the $ f $-shock segment $ [A_1, M] $, see Fig.~\ref{fig:Gamma3_Up_Down}(a). The segment $ [X, \widehat{M}]_{\text{ext}} $ represents the $ s $-right extension of the $ f $-shock segment $ [X, M] $. The structure of the solution for the Riemann problem for $L \in [G, W]$ and $R\in \Gamma_3$ includes the cases $L = W$ (case (i)), $L = G$ (case (iv)), and $L \in (G, B)$ (case (v)) as presented previously. The remaining cases are:
\begin{enumerate}
   \item[(ii)] $L \testright{R_s} Y \testright{R_f} A_1 \xrightarrow{S_f} R$, for $L \in (W, L_1)$, and $Y \in (W, A_1)\subset \mathcal{W}_f(R)$.
   \item[(iii)] $L \testright{R_s} Y \testright{S_f} R$, for $L \in [L_1, L_X]$ and $Y \in [A_1, X]\subset \mathcal{W}_f(R)$.
   \item[(vi)] $L \testright{R_s} T \xrightarrow{S_s} Y \testright{S_f} R$, for $L \in (L_X, B]$, $T \in (X, \widehat{M}]_{\text{ext}}$ and $Y \in (X, M]\subset \mathcal{W}_f(R)$.
\end{enumerate}
\end{rem}

\subsubsection{Riemann problem solution for \texorpdfstring{$R\in \Gamma_4$}{RET4}}
\label{sec:Gamma4}

This region corresponds to the beige area in Fig.~\ref{fig:Subdivision_Gamma}(a). The boundaries of $\Gamma_4$ consist of: 
\begin{itemize}
    \item The segment $[S, D]$, which is contained within the invariant segment $[G, D]$;
    \item The segment $[D, \mathcal{I}_D]$, which is part of the boundary of $\Omega$;
    \item The red curve $[\mathcal{I}_D, M_E^2]$, belonging to the $f$-inflection locus; and 
    \item The blue curve $[M_E^2,S]$, which is part of the $s$-inflection locus. 
\end{itemize}
Thus, for $L \in [G, W]$ and $R \in \Gamma_4$, the solution of the Riemann problem includes both classical and non-classical wave groups.

\begin{figure}[ht]
	\centering
    \subfigure[Structure of the Riemann solution for a generic $R\in\Gamma_4$ with $M\in{(S,D_0]}$, located to the left of the curve $M_E$, for $R = (0.0612875,0,180323)$.]%
	{\includegraphics[width=0.45\linewidth]{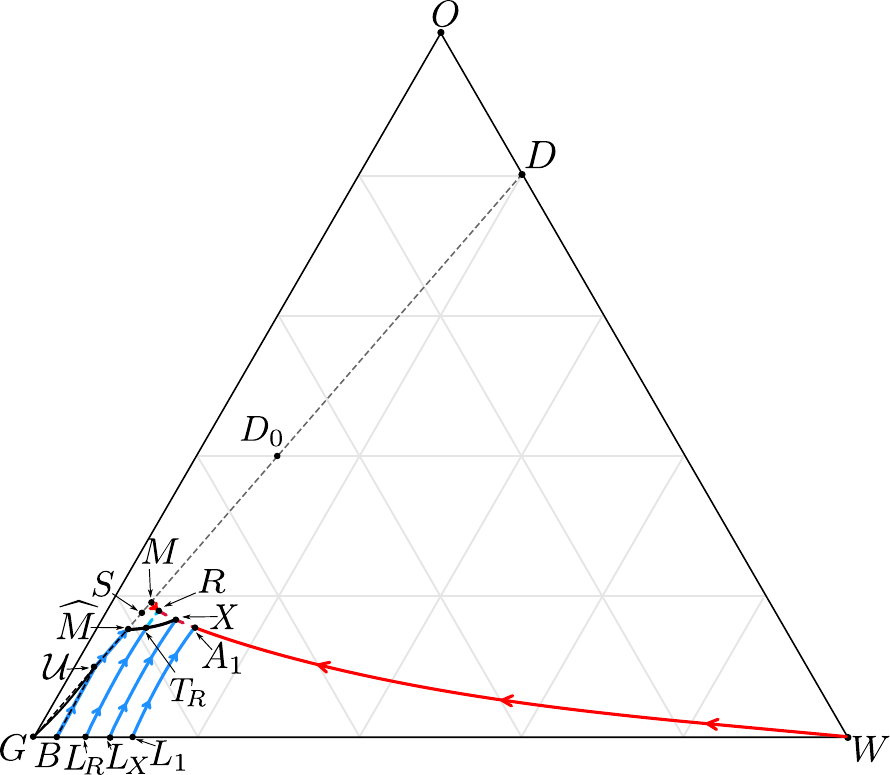}}
  	\hspace{0.8mm}
	 \subfigure[Structure of the Riemann solution for a generic $R\in\Gamma_4$ with $M\in{(D_0,D]}$, located to the right of the curve $M_E$, for $R = (0.0985438, 0.26204)$.]%
	{\includegraphics[width=0.48\linewidth]{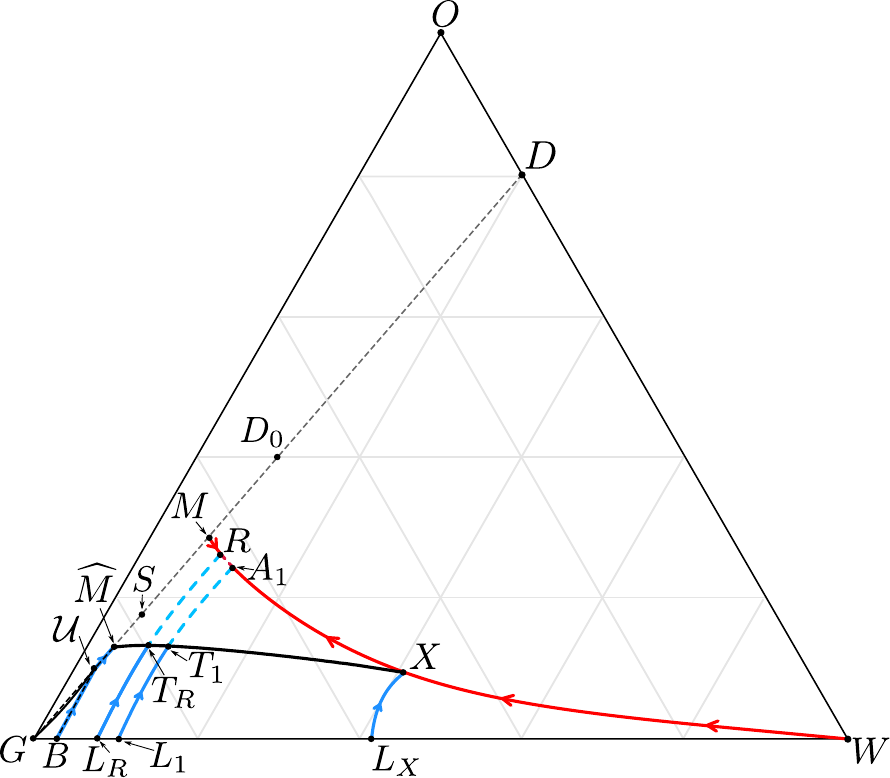}}  
	\caption{The blue dashed curves (respectively, red) represent $s$-shock curves (respectively, $f$-shock curves). The blue continuous curve (respectively, red) represents the $s$-rarefaction curves (respectively, the $f$-rarefaction curves). The arrows indicate increasing characteristic velocity. The black curve $[X, \widehat{M}]_{\text{ext}}$ represents the $s$-right extension of $\mathcal{W}_f(R)$. The black curve $[G, \mathcal{U}]_{\text{ext}}$ represents the $s$-right extension of $[G,\mathcal{U}]$.
		}
	\label{fig:Gamma4}
\end{figure}

\begin{figure}[ht]
	\centering
	 \subfigure[Structure of the Riemann solution for generic $R\in\Gamma_4$ with $M\in{(D_0,D]}$ and $R = (0.217258, 0.745896)$.]
	{\includegraphics[width=0.48\linewidth]{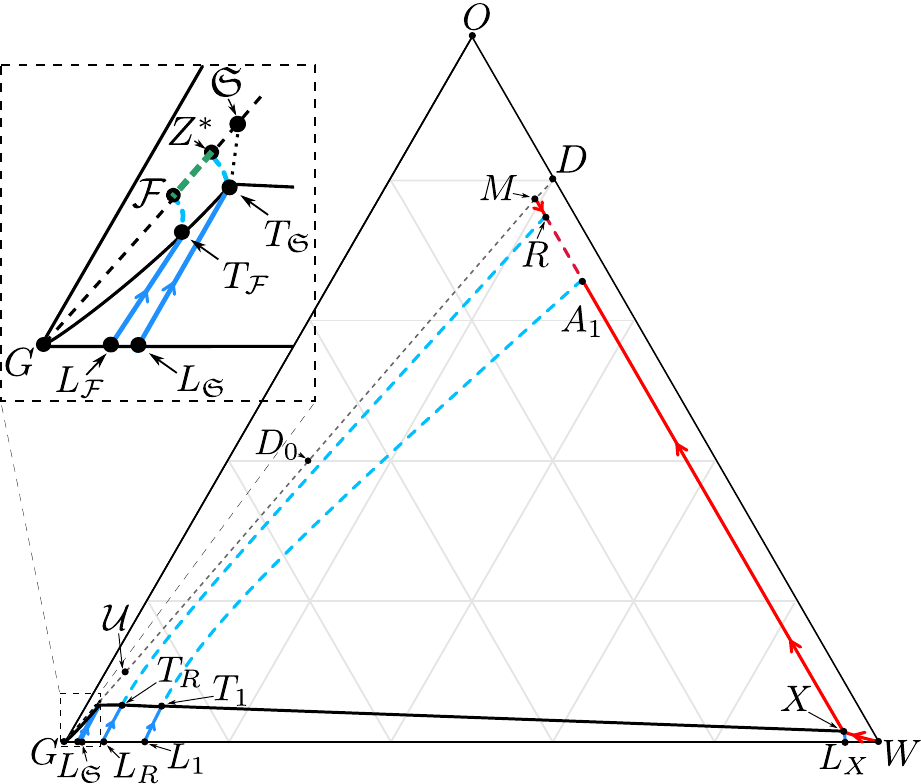}}  
	\caption{The blue dashed curves (respectively, red) represent $s$-shock curves (respectively, $f$-shock curves). The blue continuous curve (respectively, red) represents the $s$-rarefaction curves (respectively, the $f$-rarefaction curves). The arrows indicate increasing characteristic velocity. The black curve  $[X, T_{\mathfrak{S}}]_{\text{ext}}$ represents the $s$-right extension of $\mathcal{W}_f(R)$. The green dotted segment $[\mathcal{F}, Z^*] \subset$ $[G, D]$ identifies admissible $u$-shocks with $M$.  The black curve $[G, T_{\mathfrak{S}}]_{\text{ext}}$ is the $s$-right extension of $[G, Z^*]$.
		}
	\label{fig:Gamma4_D0_D}
\end{figure}

According to \cite{lozano2018, lozano2024b}, $\mathcal{W}_f(R)$ contains exclusively a local branch that connects the states $O$, $R$, and $W$. Figures \ref{fig:Gamma4} and \ref{fig:Gamma4_D0_D} displays the components of $\mathcal{W}_f(R)$ used to construct the solution of the Riemann problem for $L \in [G, W]$: 
\begin{itemize}
    \item The $f$-shock curve $[A_1, R)$; and
    \item The $f$-rarefaction curves $[W, A_1)$, and $[R, M]$ with $\lambda_f(A_1) = \sigma(A_1; R)$.
\end{itemize}
As in $\Gamma_3$, the state $X$ admits two possible configurations:
\begin{itemize}
    \item If $R$ is on the left of the curve $M_E$ (Figs.~\ref{fig:Macro_Classif_3foam} and \ref{fig:Subdivision_Gamma}), then $X$ belongs to the $f$-shock curve $[A_1, R)$ (Fig.~\ref{fig:Gamma4}(a)). 
    \item If $R$ is on the right of the curve $M_E$ (Figs.~\ref{fig:Macro_Classif_3foam} and \ref{fig:Subdivision_Gamma}), then $X$ belongs to the $f$-rarefaction curve $[W, A_1)$ (Figs.~\ref{fig:Gamma4}(b) and \ref{fig:Gamma4_D0_D}).
\end{itemize}
We focus here on $X \in [W, A_1)$; the alternative case appears in Remark \ref{rem:X_ME_Gamma_4}.

The solution classification parallels $\Gamma_3$, but with key differences: since $R$ lies below $[G, D]$, the intersection state $M\in\mathcal{W}_f(R)\cap [G, D]$  now falls within an $f$-rarefaction curve rather than an $f$-shock curve (Figs.~\ref{fig:Gamma4} and \ref{fig:Gamma4_D0_D}). This change alters the elementary wave from $M$ to $R$ to an $f$-rarefaction instead of an $f$-shock (see Section \ref{sec:Gamma3}).

As in the previous section, the type of non-classical solution used to solve the Riemann problem depends on whether $ M \in (S, D_0] $ or $ M \in (D_0, D] $:
 \begin{itemize}
    \item For $ M \in (S, D_0] $, there exists a state $ \widehat{M} \in [\mathcal{U}, S) $ that enables $ M $ to be reached via the wave structure \eqref{eq:transRarefactionM3} (Fig.~\ref{fig:Gamma4});
    \item For $ M \in (D_0, D] $, $ M $ defines the undercompressive shock segment $[\mathcal{F}, Z^*] \subset [\mathcal{F}, \mathfrak{S}]$, which satisfies the velocity compatibility condition \eqref{eq:compatib} for $ s $-, $ u $-, and $ f $-shocks, see Fig.~\ref{fig:Gamma4_D0_D}. Here, $ Z^* \in \mathcal{H}(M) $ and $ T_{\mathfrak{S}} \in \mathcal{H}(M) $ both satisfy \eqref{eq:zetaestrella} (Remark \ref{rem:expl_M}).
\end{itemize}

The segment $(G, \mathcal{U})_{\text{ext}}$ represents the $s$-right extension of $(G, \mathcal{U})$ (see Fig.~\ref{fig:Gamma4}). Analogously, the segments $(G, T_{\mathcal{F}}]_{\text{ext}} \cup (T_{\mathcal{F}}, T_{\mathfrak{S}}]_{\text{ext}}$ correspond to the $s$-right extension of $(G, \mathcal{F}] \cup (\mathcal{F}, Z^*] \subset [G, D]$ (see Fig.~\ref{fig:Gamma4_D0_D}).

The extension of $\mathcal{W}_f(R)$ for $R \in \Gamma_4$ is given by:
\begin{itemize}
    \item The segment $[X,\widehat{M}]_{\text{ext}}$, if $R$ lies to the left of the curve $M_E$ and $M \in (S, D_0]$ (Fig.~\ref{fig:Gamma4}(a)).  
    \item The segments $[X,T_1]_{\text{ext}} \cup [T_1, \widehat{M}]_{\text{ext}}$, if $R$ lies to the right of the curve $M_E$ and $M \in (S, D_0]$ (Fig.~\ref{fig:Gamma4}(b)).    
    \item The segments $[X,T_1]_{\text{ext}} \cup [T_1, T_{\mathfrak{S}})_{\text{ext}}$, if $R$ lies to the right of the curve $M_E$ and $M \in (D_0, D]$ (Fig.~\ref{fig:Gamma4_D0_D}).  
\end{itemize}

Let $L_1, L_{\mathcal{F}}, L_{\mathfrak{S}}, L_R$, and $L_X$ denote the intersection points of the backward $s$-wave curves through the states $A_1$, $\mathcal{F}$, $Z^*$, $R$, and $X$ with the edge $[G, W]$; see Figs.~\ref{fig:Gamma4} and \ref{fig:Gamma4_D0_D}. Then, the structure of the solution for the Riemann problem for $L \in [G, W]$ and $R \in \Gamma_4$ (see Fig.~\ref{fig:Subdivision_Gamma}(a)) consists of:
\begin{enumerate}
     \item[(i)] $L\testright{R_f} A_1 \xrightarrow{S_f} R$, for $L = W$;
     \item[(ii)] $L\testright{R_s}  Y \testright{R_f} A_1 \xrightarrow{S_f} R$, for $L\in(W, L_{X}]$ and $Y\in(W, X]\subset \mathcal{W}_f(R) $;
     \item[(iii)] $L\testright{R_s} T \xrightarrow{S_s} Y \testright{R_f} A_1 \xrightarrow{S_f} R$, for $L\in(L_X,L_{1})$, $T\in(X, T_{1})_{\text{ext}}$ and $Y\in(X, A_1)\subset \mathcal{W}_f(R) $;
     \item[(iv)] $L\testright{R_s} T \xrightarrow{S_s} Y \testright{S_f} R$, for $L\in[L_1,L_{R})$, $T\in[T_1,T_{R})_{\text{ext}}$ and $Y\in[A_1, R)\subset \mathcal{W}_f(R)$;
     \item[(v)] $L\testright{R_s} T \xrightarrow{S_s} R$, for $L=L_R$.
     \end{enumerate}
\begin{itemize}
       \item[$\bullet$] Case $M\in (S, D_0]$, see Fig.~\ref{fig:Gamma4}(b):
    \begin{enumerate}
          \item[(vi)] $L\testright{R_f} \mathcal{U}\xrightarrow{R_s} \widehat{M}  \xrightarrow{S_s} M \testright{R_f} R$, for $L = G$;
          \item[(vii)] $L\testright{R_s} T \xrightarrow{S_s} N \testright{R_f}\mathcal{U}\xrightarrow{R_s} \widehat{M}  \xrightarrow{S_s} M \testright{R_f} R$, for $L\in(G, B)$, $T\in(G, \mathcal{U})_{\text{ext}}$ and $N\in(G, \mathcal{U})\subset[G, D]$.
           \item[(viii)] $L\testright{R_s} T \xrightarrow{S_s} Y \testright{R_f} R,$ for $L\in(L_R,B]$, $T\in(T_R, \widehat{M}]_{\text{ext}}$ and $Y\in(R,M]\subset \mathcal{W}_f(R)$.
      \end{enumerate}
       \item[$\bullet$] Case $M\in (D_0, D]$, see Fig.~\ref{fig:Gamma4_D0_D}:
    \begin{enumerate}
          \item[(vi)] $L\testright{R_f} \mathcal{F}\xrightarrow{S_u} M \testright{R_f} R$, for $L = G$;
           \item[(vii)] $L\testright{R_s} T \xrightarrow{S_s} N\testright{R_f} \mathcal{F}  \xrightarrow{S_u} M \testright{R_f} R$, for $L\in(G,L_{\mathcal{F}})$, $T\in(G,T_{\mathcal{F}})_{\text{ext}}$ and $N\in(G,\mathcal{F})\subset[G, D]$;
          \item[(viii)] $L\testright{R_s} T \xrightarrow{S_s} N  \testright{S_u} M \testright{R_f} R$, for $L\in[L_{\mathcal{F}},L_{\mathfrak{S}}]$, $T\in[T_{\mathcal{F}},T_{\mathfrak{S}}]_{\text{ext}}$ and $N\in[\mathcal{F},Z^*]\subset[G, D]$;
          \item[(ix)] $L\testright{R_s} T \xrightarrow{S_s} Y \testright{R_f} R$, for $L\in(L_R,L_{\mathfrak{S}})$, $T\in(T_R, T_{\mathfrak{S}})_{\text{ext}}$ and $Y\in(R, M)\subset \mathcal{W}_f(R)$.
     \end{enumerate}
\end{itemize}
\begin{rem}\label{rem:X_ME_Gamma_4}
Consider the case where $M \in (S, D_0]$ and $X$ lies on the $f$-shock curve $[A_1, R]$, see Fig.~\ref{fig:Gamma4}(a). The segments $[X, T_R)_{\text{ext}}$ and $[T_R,\widehat{M}]_{\text{ext}}$ represent the $s$-right extensions of the $f$-shock curve $[X, R)$ and the $f$-rarefaction curve $[R, M]$, respectively. The structure of the solution for the Riemann problem for $L \in [G, W]$ and $R\in \Gamma_4$ includes the cases $L = W$ (case (i)), $L = L_R$ (case (v)), $L = G$ (case (vi)), and $L \in (G, B)$ (case (vii)) as presented previously. The remaining cases are:

\begin{enumerate}
   \item[(ii)] $L \testright{R_s} Y \testright{R_f} A_1 \xrightarrow{S_f} R$, for $L \in (W, L_1)$ and $Y \in (W, A_1)\subset \mathcal{W}_f(R)$.
   
   \item[(iii)] $L \testright{R_s} Y \testright{S_f} R$, for $L \in [L_1, L_X]$ and  $Y \in [A_1, X]\subset \mathcal{W}_f(R)$.

     \item[(iv)] $L \testright{R_s} T \xrightarrow{S_s} Y  \testright{S_f} R$, for $L \in (L_X, L_R)$, $T \in (T_R, X)_{\text{ext}}$ and  $Y \in (X, R)\subset \mathcal{W}_f(R)$.
   
   \item[(viii)] $L \testright{R_s} T \xrightarrow{S_s} Y \testright{R_f} R$, for $L \in (L_R, B]$, $T \in (T_R,\widehat{M}]_{\text{ext}}$ and $Y \in (R,M]\subset \mathcal{W}_f(R)$.
\end{enumerate}
\end{rem}

\subsubsection{Riemann problem solution for \texorpdfstring{$R\in \Gamma_5$}{RET5}} \label{sec:Gamma5}

This case corresponds to the green region in Fig.~\ref{fig:Subdivision_Gamma}(a). The boundaries of $\Gamma_5 $ are defined by the following segments:  
\begin{itemize}
\item The $f$-inflection segment $ [M_E^2, \mathcal{I}_D]$ (red curve);  
\item The boundary segment of $ \Omega $, $ [\mathcal{I}_D, D_E] $;  
\item The segment $ [D_E, \mathcal{I}_s^4] $, which is part of the $f$-left-extension of the invariant segment $ [G, D] $; and  
\item The blue curve $[\mathcal{I}_s^4, M_E^2] $, which is part of the $s$-inflection locus $ \mathcal{I}_s $.  
\end{itemize}
For $ L \in [G, W] $ and $ R \in \Gamma_5 $, the solution to the Riemann problem consists of both classical and non-classical wave groups.  

\begin{figure}[ht]
    \centering
      \subfigure[Structure of the Riemann solution for generic $R\in\Gamma_5$ with $M\in{(S, D_0]}$ and $R = (0.206742, 0.287839)$.]
      {\includegraphics[width=0.45\linewidth]{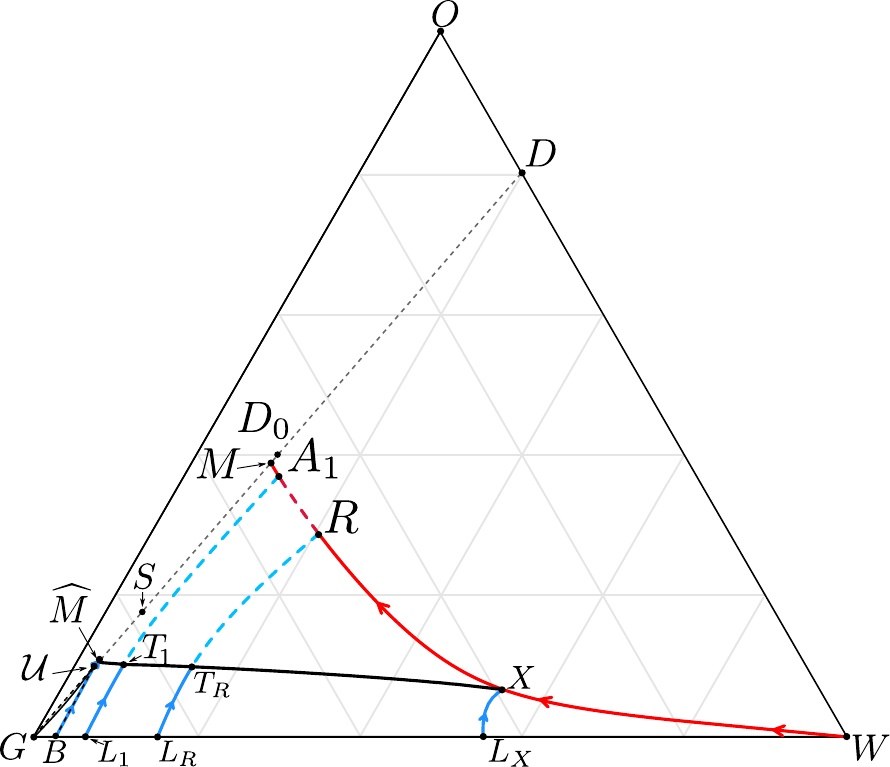}}
    \hspace{0.8mm}
     \subfigure[Structure of the Riemann solution for generic $R\in\Gamma_5$ with $M\in{(D_0,D]}$ and $R = (0.37444, 0.561621)$.]
     {\includegraphics[width=0.48\linewidth]{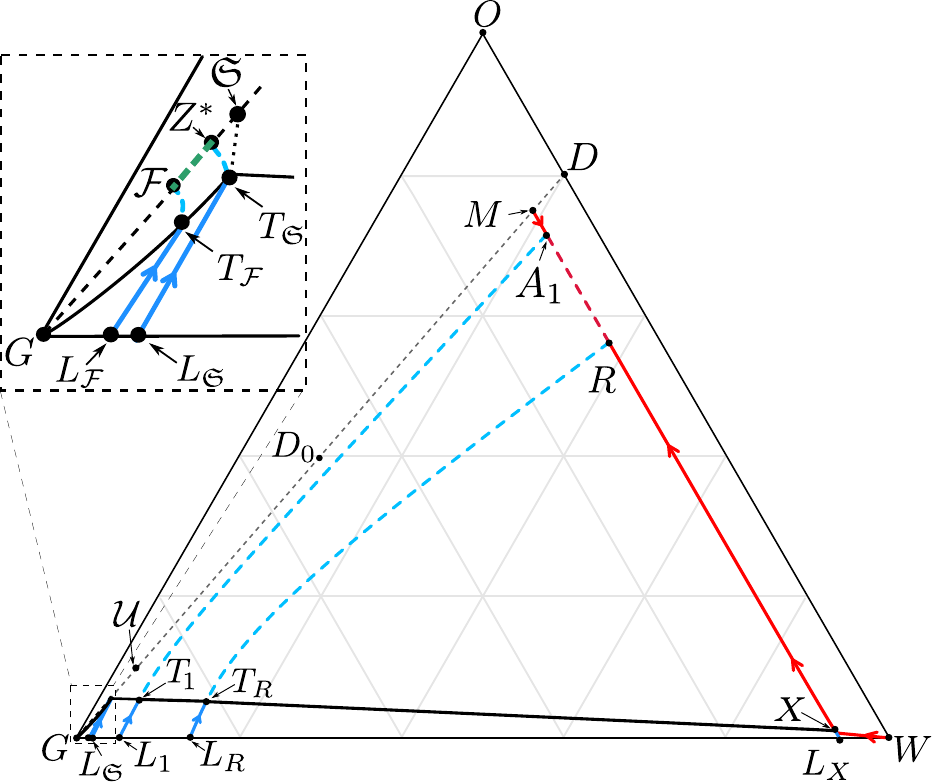}}
       \caption{The blue dashed curve (respectively red) represents $s$-shock curves (respectively $f$-shock curves). The blue continuous curve (respectively red) represents $s$-rarefaction curves (respectively the $f$-rarefaction curve). The arrows indicate the increasing characteristic velocity.  (a) The black curve $[X, \widehat{M}]_{\text{ext}}$ represents the $s$-right extension of $\mathcal{W}_f(R)$. The black curve $[G, \mathcal{U}]_{\text{ext}}$ represents the $s$-right extension of $[G,\mathcal{U}]$.  (b) The black curve $[X, T_{\mathfrak{S}}]_{\text{ext}}$ represents the $s$-right extension of $\mathcal{W}_f(R)$. The green dotted segment $[\mathcal{F}, Z^*] \subset$ $[G, D]$ identifies admissible $u$-shocks with $M$.  The black curve $[G, T_{\mathfrak{S}}]_{\text{ext}}$ represents the $s$-right extension of $[G, Z^*]$.
    }
    \label{fig:Gamma5}
\end{figure}

Following the approach outlined for $ \Gamma_4 $, for $ R\in\Gamma_5 $, $\mathcal{W}_f(R) $ contains only a local branch that connects the states $ O $, $R$, and $ W $. Figure~\ref{fig:Gamma5} shows the components of $ \mathcal{W}_f(R) $ used to construct the solution of the Riemann problem for $ L \in [G, W] $. These include:  
\begin{itemize}
\item The $f$-shock curve $ (R, A_1] $.  
\item The $f$-rarefaction curves $ [M, A_1) $ and $ [W, R) $, with $ \lambda_f(A_1) = \sigma(A_1; R) $.  
\end{itemize}
The primary distinction between $ \Gamma_5 $ and $ \Gamma_4 $ is the relative positioning of states $ A_1 $ and $R$ with respect to the $f$-inflection locus $\mathcal{I}_f$ (see Fig.~\ref{fig:Subdivision_Gamma}(a)), which leads to a different configuration of the Riemann problem solutions. As in region $ \Gamma_4 $, the state $M=\mathcal{W}_f(R)\cap[G, D]$ is located within an $f$-rarefaction curve, while the state $ X $ is defined as the intersection between the $f$-rarefaction curve $[W, R) $ and the $s$-inflection.  
The general solution structure and classification remain similar to those in $ \Gamma_4 $. 

As in the previous section, the type of non-classical solution used to solve the Riemann problem depends on whether $ M \in (S, D_0] $ or $ M \in (D_0, D] $:
 \begin{itemize}
    \item For $ M \in (S, D_0] $, there exists a state $ \widehat{M} \in [\mathcal{U}, S) $ that enables $ M $ to be reached via the wave structure \eqref{eq:transRarefactionM3} (Fig.~\ref{fig:Gamma5}(a));
    \item For $ M \in (D_0, D] $, $ M $ defines the undercompressive shock segment $[\mathcal{F}, Z^*] \subset [\mathcal{F}, \mathfrak{S}]$, which satisfies the velocity compatibility condition \eqref{eq:compatib} for $ s $-, $ u $-, and $ f $-shocks, see Fig.~\ref{fig:Gamma5}(b). Here, $ Z^* \in \mathcal{H}(M) $ and $ T_{\mathfrak{S}} \in \mathcal{H}(M) $ both satisfy \eqref{eq:zetaestrella} (Remark \ref{rem:expl_M}).
\end{itemize}

The segments $(T_R, T_1]_{\text{ext}}$ and $(X, T_R]_{\text{ext}}$ are the $s$-right-extensions of $(R, A_1]$ and $(X, R]$, respectively, which belong to $\mathcal{W}_f(R)$. The $s$-right-extension of $[A_1, M] \subset \mathcal{W}_f(R)$ is given by $[\widehat{M}, T_1]_{\text{ext}}$ when $M \in (S, D_0]$ (Fig.~\ref{fig:Gamma5}(a)) or by $(T_{\mathfrak{S}}, T_1]_{\text{ext}}$ if $M \in (D_0, D]$ (Fig.~\ref{fig:Gamma5}(b)).  

The segment $(G, T_{\mathcal{F}}]_{\text{ext}} \cup (T_{\mathcal{F}}, T_{\mathfrak{S}}]_{\text{ext}}$ represents the $s$-right-extension of $(G, \mathcal{F}] \cup (\mathcal{F}, Z^*] \subset [G, D]$ (Fig.~\ref{fig:Gamma5}(a)), while $(G, \mathcal{U})_{\text{ext}}$ is the $s$-right-extension of $(G, \mathcal{U}) \subset [G, D]$ (Fig.~\ref{fig:Gamma5}(b)).  

Let $L_1$, $L_{\mathcal{F}}$, $L_{\mathfrak{S}}$, $L_R$, and $L_X$ denote the intersection points of the backward $s$-wave curves through states $A_1$, $\mathcal{F}$, $Z^*$, $R$, and $X$ with the edge $[G, W]$ (see Fig.~\ref{fig:Gamma5}). Then, the solution structure for the Riemann problem with $L \in [G, W]$ and $R \in \Gamma_5$ (Fig.~\ref{fig:Subdivision_Gamma}(a)) consists of:  
\begin{enumerate}
     \item[(i)] $L\testright{R_f} R$, for $L = W$;
     \item[(ii)] $L\testright{R_s}  Y \testright{R_f} R$, for $L\in(W, L_{X}]$ and $Y\in(W, X]\subset\mathcal{W}_f(R) $;
     \item[(iii)] $L\testright{R_s} T \xrightarrow{S_s} Y \testright{R_f} R$, for $L\in(L_X,L_{R})$, $T\in(X,T_{R})_{\text{ext}}$ and $Y\in(X, R)\subset\mathcal{W}_f(R)$;
     \item[(iv)] $L\testright{R_s} T_R \xrightarrow{S_s} R$, for $L = L_R$;
     \item[(v)] $L\testright{R_s} T \xrightarrow{S_s} Y \testright{S_f} R$, for $L\in(L_{R}, L_1]$, $T\in (T_R, T_{1}]_{\text{ext}}$ and $Y\in(R, A_1]\subset\mathcal{W}_f(R)$.
     \end{enumerate}
\begin{itemize}
  \item[$\bullet$] Case $M\in (S, D_0]$, Fig.~\ref{fig:Gamma5}(a):
    \begin{enumerate}
          \item[(vi)] $L\testright{R_f} \mathcal{U}\xrightarrow{R_s} \widehat{M}  \xrightarrow{S_s} M  \testright{R_f} A_1 \xrightarrow{S_f} R$, for $L = G$;
          \item[(vii)] $L\testright{R_s} T \xrightarrow{S_s} N \testright{R_f}\mathcal{U}\xrightarrow{R_s} \widehat{M}  \xrightarrow{S_s} M  \testright{R_f} A_1 \xrightarrow{S_f} R$, for $L\in(G, B)$, $T\in(G, \mathcal{U})_{\text{ext}}$ and $N\in(G, \mathcal{U})\subset[G,D]$;
           \item[(viii)] $L\testright{R_s} T \xrightarrow{S_s} Y  \testright{R_f} A_1 \xrightarrow{S_f} R$, for $L\in(L_1,B]$, $T\in(T_1,\widehat{M}]_{\text{ext}}$ and $Y\in(A_1, M]\subset\mathcal{W}_f(R)$.
      \end{enumerate}
  \item[$\bullet$] Case $M\in (D_0, D]$, see Fig.~\ref{fig:Gamma5}(b):
    \begin{enumerate}
          \item[(vi)] $L\testright{R_f} \mathcal{F}\xrightarrow{S_u} M  \testright{R_f} A_1 \xrightarrow{S_f} R$, $L = G$;
           \item[(vii)] $L\testright{R_s} T \xrightarrow{S_s} N\testright{R_f} \mathcal{F}  \xrightarrow{S_u} M \testright{R_f} A_1 \xrightarrow{S_f} R$, for $L\in(G,L_{\mathcal{F}})$, $T\in(G,T_{\mathcal{F}})_{\text{ext}}$ and $N\in(G,\mathcal{F})$;
          \item[(viii)] $L\testright{R_s} T \xrightarrow{S_s} N  \testright{S_u} M \testright{R_f} A_1 \xrightarrow{S_f} R$, for $L\in[L_{\mathcal{F}},L_{\mathfrak{S}}]$, $T\in[T_{\mathcal{F}},T_{\mathfrak{S}}]_{\text{ext}}$ and $N\in[\mathcal{F},Z^*]\subset[G,D]$;
          \item[(ix)] $L\testright{R_s} T \xrightarrow{S_s} Y  \testright{R_f} A_1 \xrightarrow{S_f} R$, for $L\in[L_1, L_{\mathfrak{S}})$, $T\in[T_1, T_{\mathfrak{S}})_{\text{ext}}$ and $Y\in[A_1,M)\subset\mathcal{W}_f(R)$.
     \end{enumerate}
\end{itemize}


\subsubsection{Riemann problem solution for \texorpdfstring{$R\in \Gamma_6$}{RET6}}
\label{sec:Gamma6}

This region corresponds to the light pink area in Fig.~\ref{fig:Subdivision_Gamma}(a). The boundaries of region $ \Gamma_6 $ are defined by the following segments: 
\begin{itemize}
    \item The black curve $[D_E, \mathcal{I}_s^4]$, which is part of the $f$-left-extension of the invariant segment $ [G, D] $;
    \item The blue curve $[\mathcal{I}_s^4, \mathcal{I}_s^5 ]$, which corresponds to part of the $s$-inflection locus; 
    \item The green curve $[\mathcal{I}_s^5, K_D]$, which is the boundary where the admissibility of nonlocal $f$-shocks changes see \cite{Andrade2018, lozano2018}; and
    \item The boundary segment of $ \Omega $, given by $ [K_D, D_E] $. 
\end{itemize} 
In this case, for $L \in [G, W]$ and $R \in \Gamma_6$, the solution to the Riemann problem consists of both classical and non-classical wave groups.

This case is similar to $ \Gamma_3 $, where $ \mathcal{W}_f(R) $ contains only a local branch connecting states $ O $, $R$, and $ W $.
Notice that, in the previous case where $ R \in \Gamma_5 $, the $f$-shock curve $(R, A_1]$ lies below the segment $[G, D]$. When the state $R$ moves from region $ \Gamma_5 $ onto the segment $[D_E, \mathcal{I}_s^4]$ (see Fig.~\ref{fig:Subdivision_Gamma}(a)), the state $ A_1\in\mathcal{W}_f(R) $ reaches the invariant segment $[G, D]$, since $ \lambda_f(A_1) = \sigma(A_1; R) $. When the state $R$ crosses $[D_E, \mathcal{I}_s^4]$, i.e., when $ R \in \Gamma_6 $, the state $ A_1 $ also crosses the invariant segment $[G, D]$, causing the intersection state $M=\mathcal{W}_f(R)\cap[G, D]$ to once again belong to an $f$-shock curve. However, since $A_1$ lies above $[G, D]$, it is no longer relevant for solutions with $ L \in [G, W] $.  

\begin{figure}[ht]
    \centering
    \subfigure[Structure of the Riemann solution for generic $R\in\Gamma_6$ with $M\in{[S,D_0]}$ and $R = (0.337019, 0.175838)$.]
    {\includegraphics[width=0.45\linewidth]{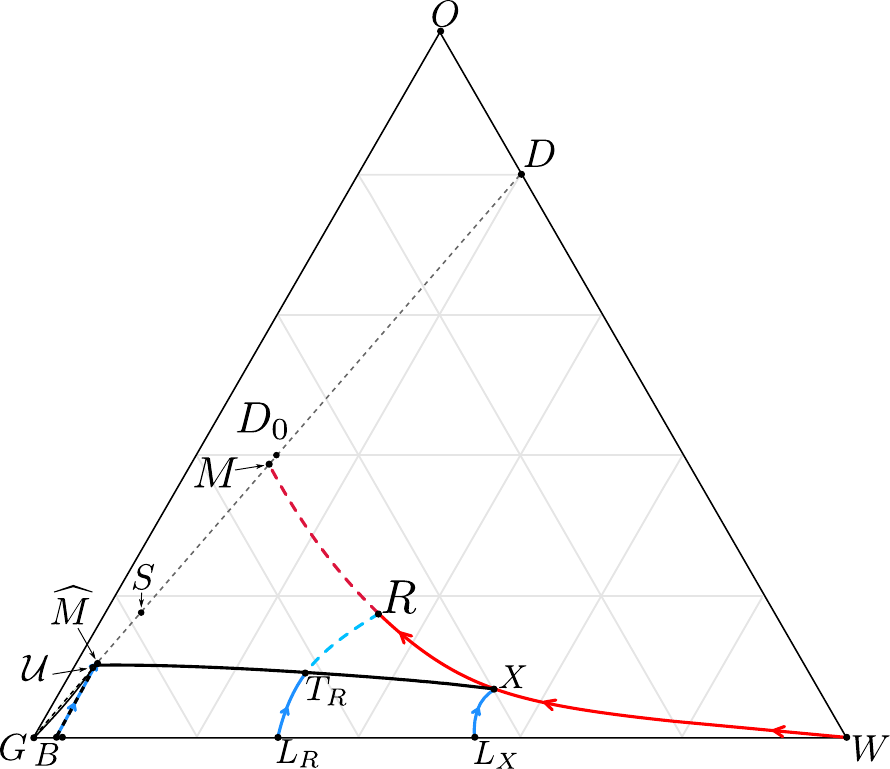}}
    \hspace{0.8mm}
    \subfigure[Structure of the Riemann solution for generic $R\in\Gamma_6$ with $M\in{[D_0,D]}$ and $R = (0.739358, 0.224073)$.]
    {\includegraphics[width=0.48\linewidth]{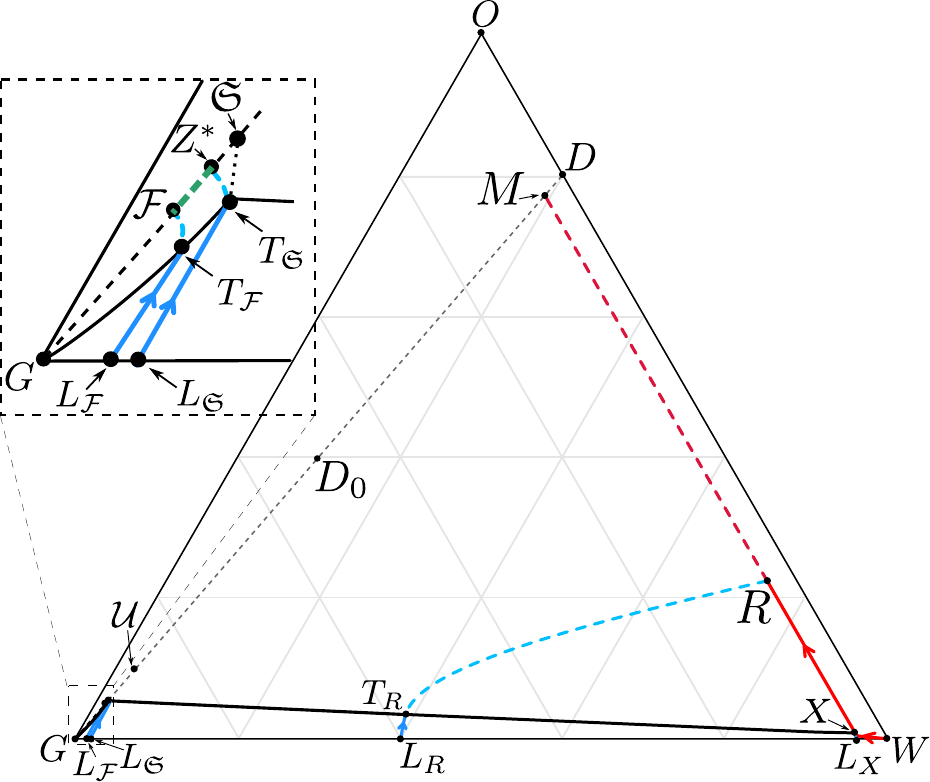}} 
    \caption{The blue dashed curve (respectively red) represents $s$-shock curves (respectively $f$-shock curves). The blue continuous curve (respectively red) represents $s$-rarefaction curves (respectively the $f$-rarefaction curve). The arrows indicate the increasing characteristic velocity.  (a) The black curve $[X, \widehat{M}]_{\text{ext}}$ represents the $s$-right extension of $\mathcal{W}_f(R)$. The black curve $[G, \mathcal{U}]_{\text{ext}}$ represents the $s$-right extension of $[G,\mathcal{U}]$.  (b) The black curve $[X, T_{\mathfrak{S}}]_{\text{ext}}$ represents the $s$-right extension of $\mathcal{W}_f(R)$. The green dotted segment $[\mathcal{F}, Z^*] \subset$ $[G, D]$ identifies admissible $u$-shocks with $M$.  The black curve $[G, T_{\mathfrak{S}}]_{\text{ext}}$ represents the $s$-right extension of $[G, Z^*]$.
    }
    \label{fig:Gamma6}
\end{figure}

Figure \ref{fig:Gamma6} illustrates the relevant portions of $ \mathcal{W}_f(R) $ for constructing the Riemann problem solution with $ L \in [G, W] $: 
\begin{itemize}
    \item The $f$-shock curve $ (R, M] $; and \item The $f$-rarefaction curve $ [W, R) $.  
\end{itemize}
As in regions $ \Gamma_3 $, $ \Gamma_4 $, and $ \Gamma_5 $, the type of non-classical solution used to solve the Riemann problem depends on whether $ M \in (S, D_0] $ or $ M \in (D_0, D] $:
 \begin{itemize}
    \item For $ M \in (S, D_0] $, there exists a state $ \widehat{M} \in [\mathcal{U}, S) $ that enables $ M $ to be reached via the wave structure \eqref{eq:transRarefactionM3} (Fig.~\ref{fig:Gamma6}(a));
    \item For $ M \in (D_0, D] $, $ M $ defines the undercompressive shock segment $[\mathcal{F}, Z^*] \subset [\mathcal{F}, \mathfrak{S}]$, which satisfies the velocity compatibility condition \eqref{eq:compatib} for $ s $-, $ u $-, and $ f $-shocks  (Fig.~\ref{fig:Gamma6}(b)). Here, $ Z^* \in \mathcal{H}(M) $ and $ T_{\mathfrak{S}} \in \mathcal{H}(M) $ both satisfy \eqref{eq:zetaestrella} (Remark \ref{rem:expl_M}).
\end{itemize}

The segment $ (X,T_R]_{\text{ext}} $ is the $s$-right-extension of $ (X, R] \subset \mathcal{W}_f(R) $. The $s$-right-extension of $ [M,R] \subset \mathcal{W}_f(R) $ is given by $ [\widehat{M},T_R]_{\text{ext}} $ when $ M \in (S, D_0] $ (Fig.~\ref{fig:Gamma6}(a)) or by $ (T_{\mathfrak{S}},T_R]_{\text{ext}} $ if $ M \in (D_0, D] $ (Fig.~\ref{fig:Gamma6}(b)).  

The segment $ (G,\mathcal{U})_{\text{ext}} $ is the $s$-right-extension of $ (G,\mathcal{U}) \subset [G, D] $ (Fig.~\ref{fig:Gamma6}(a)), while $ (G,T_{\mathcal{F}}]_{\text{ext}} \cup (T_{\mathcal{F}}, T_{\mathfrak{S}}]_{\text{ext}} $ represents the $s$-right-extension of $ (G,\mathcal{F}] \cup (\mathcal{F},Z^*] \subset [G, D] $ (Fig.~\ref{fig:Gamma6}(b)).

Let $L_{\mathcal{F}}, L_{\mathfrak{S}}, L_R$, and $L_X$ denote the intersection points of the backward $s$-wave curves through states ${\mathcal{F}}, \, Z^*, R$, and $X$ with the edge $[G, W]$; see Fig.~\ref{fig:Gamma6}. Then, the structure of the solution for the Riemann problem for $L \in [G, W]$ and $R \in \Gamma_6$ (see Fig.~\ref{fig:Subdivision_Gamma}(a)) consists of:
\begin{enumerate}
     \item[(i)] $L\testright{R_f} R$, for $L = W$;
     \item[(ii)] $L\testright{R_s}  Y \testright{R_f} R$, for $L\in(W, L_{X}]$, and $Y\in(W, X]\subset\mathcal{W}_f(R)$;
     \item[(iii)] $L\testright{R_s} T \xrightarrow{S_s} Y \testright{R_f} R$, for $L\in(L_X,L_{R})$, $T\in(X, T_{R})_{\text{ext}}$ and $Y\in(X, R)\subset\mathcal{W}_f(R)$;
     \item[(iv)] $L\testright{R_s} T_R \xrightarrow{S_s}R$, for $L = L_R$;
\end{enumerate}
\begin{itemize}
       \item[$\bullet$] Case $M\in (S, D_0]$, Fig.~\ref{fig:Gamma6}(a):
    \begin{enumerate}
          \item[(v)] $L\testright{R_f} \mathcal{U}\xrightarrow{R_s} \widehat{M}  \xrightarrow{S_s} M  \testright{S_f} R$, for $L = G$;
          \item[(vi)] $L\testright{R_s} T \xrightarrow{S_s} N \testright{R_f}\mathcal{U}\xrightarrow{R_s} \widehat{M}  \xrightarrow{S_s} M  \testright{S_f} R$, for $L\in(G, B)$, $T\in(G, \mathcal{U})_{\text{ext}}$ and $N\in(G, \mathcal{U})\subset [G, D]$.
           \item[(vii)] $L\testright{R_s} T \xrightarrow{S_s} Y  \testright{S_f} R$, for $L\in(L_{R}, B]$, $T\in(T_{R},\widehat{M} ]_{\text{ext}}$ and $Y\in(R,M]\subset\mathcal{W}_f(R)$.
      \end{enumerate}
       \item[$\bullet$] Case $M\in (D_0, D]$, see Fig.~\ref{fig:Gamma6}(b):
    \begin{enumerate}
          \item[(v)] $L\testright{R_f} \mathcal{F}\xrightarrow{S_u} M  \testright{S_f} R$, for  $L = G$;
           \item[(vi)] $L\testright{R_s} T \xrightarrow{S_s} N\testright{R_f} \mathcal{F}  \xrightarrow{S_u} M \testright{S_f} R$, for $L\in(G,L_{\mathcal{F}})$, $T\in(G,T_{\mathcal{F}})_{\text{ext}}$ and $N\in(G,\mathcal{F})\subset [G, D]$.
          \item[(vii)] $L\testright{R_s} T \xrightarrow{S_s} N  \testright{S_u} M \testright{S_f} R$, for $L\in[L_{\mathcal{F}},L_{\mathfrak{S}}]$, $T\in[T_{\mathcal{F}},T_{\mathfrak{S}}]_{\text{ext}}$ and $N\in[\mathcal{F},Z^*]\subset [G, D]$.
          \item[(viii)] $L\testright{R_s} T \xrightarrow{S_s} Y  \testright{S_f} R$, for $L\in(L_R, L_{\mathfrak{S}})$, $T\in(T_R, T_{\mathfrak{S}})_{\text{ext}}$ and $Y\in(R, M)\subset\mathcal{W}_f(R)$.
     \end{enumerate}
\end{itemize}


\subsubsection{Riemann problem solution for \texorpdfstring{$R\in \Gamma_7$}{RET7}}
\label{sec:Gamma7}
\begin{figure}[h!]
    \centering
    \subfigure[Structure of the Riemann solution for generic $R\in\Gamma_7$ with $R = (0.79196, 0.0327602)$. Case without intersection between the nonlocal branch of $\mathcal{W}_f(R)$ and $\mathcal{I}_s.$ ]{\includegraphics[width=0.485\linewidth]{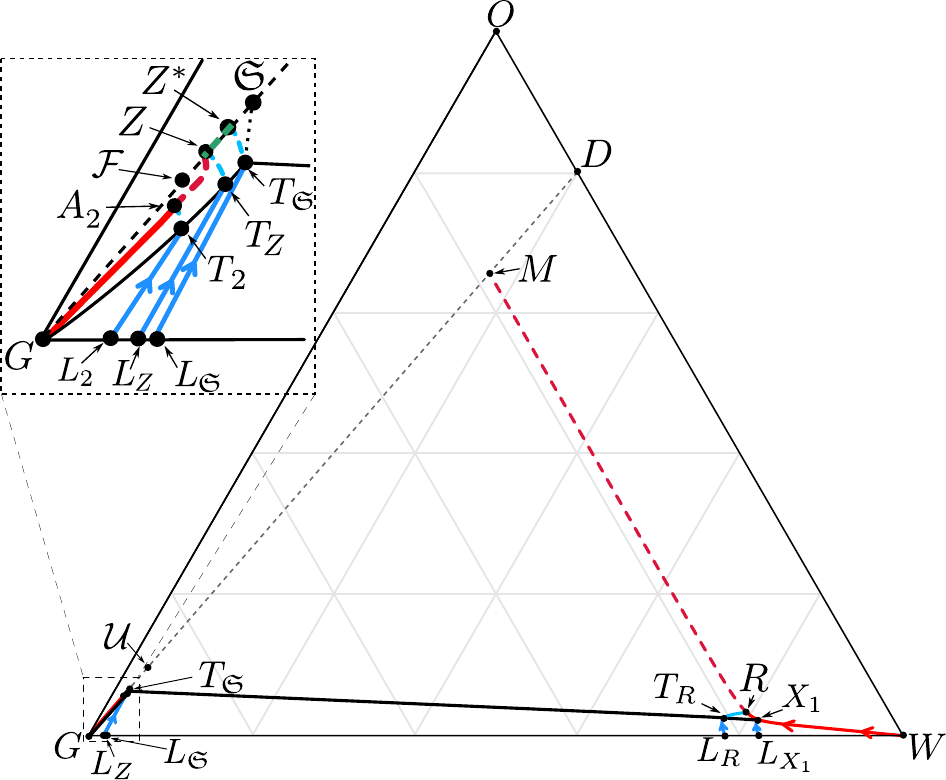}}
    \hspace{0.8mm}
    \subfigure[Structure of the Riemann solution for generic $R\in\Gamma_7$ with $R = (0.924898, 0.0282079)$. Case with $X_2$ and $X_3$ belong to the nonlocal branch of $\mathcal{W}_f(R)$ and $\mathcal{W}_f(R)\cap\mathcal{I}_s.$]{\includegraphics[width=0.485\linewidth]{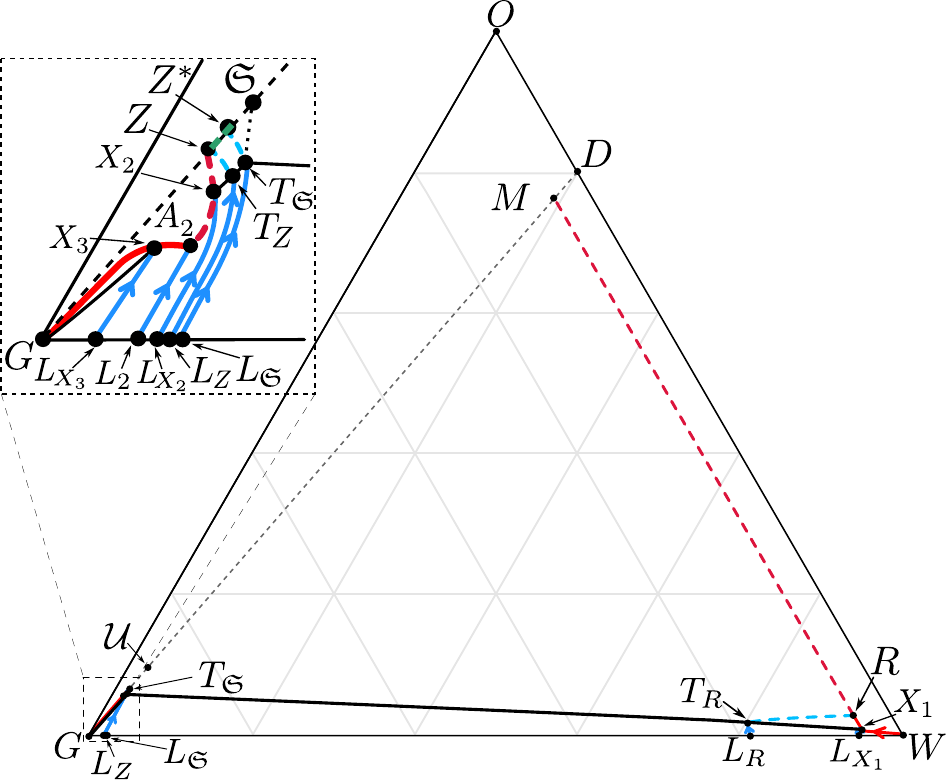}}
    \caption{The blue dashed curve (respectively red) represents $s$-shock curves (respectively $f$-shock curves). The blue continuous curve (respectively red) represents $s$-rarefaction curves (respectively the $f$-rarefaction curve). The arrows indicate the increasing characteristic velocity.
    The green dotted segment $(Z, Z^*] \subset [G, D]$ identifies admissible $u$-shocks with $M$. The black curve $(T_Z, T_{\mathfrak{S}}]_{\text{ext}}$ is the $s$-extension of $(Z, Z^*]$. (a) The black curves $[X_1,T_{\mathfrak{S}}]_{\text{ext}}\cup[G,T_Z]_{\text{ext}}$ represent the $s$-right extension of $\mathcal{W}_f(R)$. (b) The black curves $[X_1,T_{\mathfrak{S}}]_{\text{ext}}\cup[G,X_3]_{\text{ext}}\cup[X_2,T_z]_{\text{ext}}$ represent the $s$-right extension of $\mathcal{W}_f(R)$. }
    \label{fig:Gamma7}
\end{figure}

This region corresponds to the blue area in Figs.~\ref{fig:Subdivision_Gamma}(a)-(c). The boundaries of $ \Gamma_7 $ are defined similarly to those of region $ \Gamma_2 $ (although $ \Gamma_7 $ is close to $ W $ and $ \Gamma_2 $ is close to $O$). The regions $ \Gamma_6 $ and $ \Gamma_7 $ are separated by the green curve where the admissibility of nonlocal $f$-shocks changes, $[\mathcal{I}_s^5, K_D]$; see \cite{Andrade2018,lozano2018}. This implies that when $R$ crosses this boundary, moving from region $ \Gamma_6 $ to region $ \Gamma_7 $, $ \mathcal{W}_f(R) $ admits a nonlocal branch containing the state $G$; see Fig.~\ref{fig:Gamma7}. The remaining boundaries of $ \Gamma_7 $ are:  
\begin{itemize}  
    \item The blue curve $ [\mathcal{I}_s^5, \mathcal{I}_s^6]$, which is part of the $s$-inflection locus;  
    \item The purple curve $ [\mathcal{I}_s^6, V_D] $, which marks the boundary where the admissibility of $u$-shocks changes; and 
    \item The boundary segment of $ \Omega $, $ [K_D, V_D] $.  
\end{itemize}

For $L \in [G, W]$ and $R \in \Gamma_7$, the solution to the Riemann problem consists of both classical and non-classical wave groups. The wave curve $\mathcal{W}_f(R)$ comprises the following components (see Fig.~\ref{fig:Gamma7}):
\begin{itemize}
\item A local branch containing the state $R$, which intersects the $s$-inflection locus at the state $X_1 \in [W, R)$, lying within an $f$-rarefaction curve;
\item A nonlocal branch containing the state $G$, which may intersect the $s$-inflection locus at one or two states.
\end{itemize}
For the intersections of the nonlocal branch with the $s$-inflection locus, we denote a single intersection state as $X_2$, which lies within an $f$-rarefaction curve. When two intersection states exist, we denote them as $X_2$ and $X_3$; they may both lie on an $f$-shock curve, or one may lie on an $f$-rarefaction curve and the other on an $f$-shock curve. In this section, we illustrate the cases where $X_2$ does not exist (Fig.~\ref{fig:Gamma7}(a)), and where $X_2$ and $X_3$ belong to different elementary curves (Fig.~\ref{fig:Gamma7}(b)). The remaining case is similar; see Remark~\ref{rem:shock_hole_gamma7}.

Figure~\ref{fig:Gamma7} shows the components of $ \mathcal{W}_f(R) $ used to construct the solution to the Riemann problem for $ L \in [G, W] $:  
\begin{itemize}  
    \item The $f$-shock curves $(R, M] $ and $ [A_2, Z] $; and  
    \item The $f$-rarefaction curves $ [W, R] $ and $ [G, A_2] $ with $ \lambda_f(A_2) = \sigma(A_2; R)$.  
\end{itemize}  
The state $M$ is defined by the intersection of the local branch of $ \mathcal{W}_f(R) $ and the invariant segment $ [G, D] $. Both $Z$ and $M$ lie on $[G, D]$, where $M$ defines the undercompressive shock segment $(Z, Z^*] \subset [\mathcal{F}, \mathfrak{S}]$. This segment satisfies velocity compatibility condition \eqref{eq:compatib} across the $s$-, $u$-, and $f$-shocks. The state $Z^* \in \mathcal{H}(M)$ is associated with the state $T_{\mathfrak{S}} \in \mathcal{H}(M)$, where both satisfy \eqref{eq:zetaestrella} (see Remark \ref{rem:expl_M}).
The configuration of the $s$-right-extension segments of $\mathcal{W}_f$ depends on the intersections of the nonlocal branch with the $s$-inflection locus:  
\begin{itemize}  
    \item If $ X_2 $ does not exist, the segments $ (G, T_2)_{\text{ext}} $, $ [T_2, T_Z)_{\text{ext}} $, $ (X_1, T_R]_{\text{ext}} $ and $ (T_R, T_{\mathfrak{S}}]_{\text{ext}} $ represent the $s$-right-extensions of the backward $f$-wave curve segments $ (G, A_2) $, $ [A_2, Z] $, $ [X_1, R) $, and $ (R, M]$, see Fig.~\ref{fig:Gamma7}(a).  
    \item If $ X_2 $ and $X_3$ exists, the segments $ (G, X_3)_{\text{ext}} $, $ (X_2, T_Z)_{\text{ext}}$, $ (X_1, T_R]_{\text{ext}} $ and $ (T_R, T_{\mathfrak{S}}]_{\text{ext}}$ represent the $s$-right-extensions of the backward $f$-wave curve segments $ (G, X_3) $, $ [X_2, Z) $, $ [X_1, R) $, and $ (R, M]$, see Fig.~\ref{fig:Gamma7}(b). 
    \end{itemize} 
 The segment $ (T_Z, T_{\mathfrak{S}}]_{\text{ext}} $ is the $s$-right-extension of $ (Z, Z^*] \subset [G, D] $.

Let $L_2, L_Z, L_{\mathfrak{S}}, L_R$, $L_{X_1}$, $L_{X_2}$ and $L_{X_3}$ denote the intersection points of the backward $s$-wave curves through states $A_2$, $Z$, $Z^*$, $R$, $X_1$, $X_2$ and $X_3$ with the edge $[G, W]$ (see Fig.~\ref{fig:Gamma7}). Then, the structure of the solution for the Riemann problem for $L \in [G, W]$ and $R \in \Gamma_7$ (see Fig.~\ref{fig:Subdivision_Gamma}(c)) consists of:
\begin{enumerate}
         \item[(i)] $L\testright{R_f} R$, for $L = W$;
         \item[(ii)] $L\testright{R_s} Y  \testright{R_f} R$, for $L\in(W, L_{X_1}]$ and $Y\in(W,X_1]\subset\mathcal{W}_f(R)$;
 \item[(iii)] $L\testright{R_s} T  \xrightarrow{S_s} Y \testright{R_f} R$, for $L\in(L_{X_1},L_{R})$, $T\in(X_1, T_{R})_{\text{ext}}$ and $Y\in(X_1,R)\subset\mathcal{W}_f(R)$;
 \item[(iv)]  $L\testright{R_s} T_R  \xrightarrow{S_s} R$, for $L=L_{R}$;
 \item[(v)] 
$L\testright{R_s} T  \xrightarrow{S_s} Y  \testright{S_f} R$, for $L\in(L_{\mathfrak{S}},L_R)$,  $T\in(T_{\mathfrak{S}}, T_{R})_{\text{ext}}$ and $Y\in(M, R)\subset\mathcal{W}_f(R)$;
 \item[(vi)]  $  L\testright{R_f} A_2  \xrightarrow{S_f} R$ for $L = G$; 
 \end{enumerate}
\begin{itemize}
    \item[$\bullet$] When $X_2$ does not exist (see Fig.~\ref{fig:Gamma7}(a)):
\begin{enumerate}
 \item[(vii)] $ L\testright{R_s} T  \xrightarrow{S_s} Y \testright{R_f} A_2 \xrightarrow{S_f} R$, for $L\in(G,L_{2})$, $T\in(G,T_{2})_{\text{ext}}$ and $Y\in(G,A_{2})\subset\mathcal{W}_f(R)$;
 \item[(viii)] $L\testright{R_s} T  \xrightarrow{S_s} Y  \testright{S_f} R$, for $L\in[L_{2},L_Z]$, $T\in[T_{2},T_Z]_{\text{ext}}$ and  $Y\in[A_{2},Z]\subset\mathcal{W}_f(R)$;
 \item[(ix)] $L\testright{R_s} T  \xrightarrow{S_s} N  \testright{S_u} M \testright{S_f} R$, for $L\in(L_{Z},L_{\mathfrak{S}}]$, $T\in(T_{Z},T_{\mathfrak{S}}]_{\text{ext}}$ and $N\in(Z,Z^*]\subset[G, D]$.
 \end{enumerate}
 \item[$\bullet$] When $X_2 $ and $X_3$ exists (see Fig.~\ref{fig:Gamma7}(b)):
\begin{enumerate}
 \item[(vii)] $ L\testright{R_s} T  \xrightarrow{S_s} Y \testright{R_f} A_2 \xrightarrow{S_f} R$, for $L\in(G,L_{X_3})$, $T\in(G,X_3)_{\text{ext}}$ and $Y\in(G,X_{3})\subset\mathcal{W}_f(R)$;
 \item[(viii)] $L\testright{R_s} Y \testright{R_f} A_2 \xrightarrow{S_f} R$, for $L\in[L_{X_3},L_2)$ and  $Y\in[X_3,A_{2})\subset\mathcal{W}_f(R)$;
 \item[(ix)] $L\testright{R_s} Y \testright{S_f} R$, for $L\in[L_{2},L_{X_2})$ and  $Y\in[A_{2},X_2)\subset\mathcal{W}_f(R)$;
 \item[(x)] $L\testright{R_s} T  \xrightarrow{S_s} Y  \testright{S_f} R$, for $L\in(L_{X_2},L_Z]$, $T\in(X_{2},T_Z]_{\text{ext}}$ and  $Y\in({X_2},Z]\subset\mathcal{W}_f(R)$;
    \item[(xi)] $L\testright{R_s} T  \xrightarrow{S_s} N  \testright{S_u} M \testright{S_f} R$, for $L\in(L_{Z},L_{\mathfrak{S}}]$, $T\in(T_{Z},T_{\mathfrak{S}}]_{\text{ext}}$ and $N\in(Z,Z^*]\subset[G, D]$;
\end{enumerate}
\end{itemize}

\begin{rem}\label{rem:shock_hole_gamma7}
   As shown in Fig.~\ref{fig:fig2}(a), the $s$-inflection locus, $\mathcal{I}_s$ consists of two distinct branches, the blue curves $[G,\mathcal{U}]$ and $[W, S] \cup [S, O]$. Within the region $\Gamma_7$, the $f$-left-extension of the segment $[G, \mathcal{U}]\subset\mathcal{I}_s$ is present. When the right state $R$ lies on this locus, the intersection state $X_2 \in \mathcal{I}_s\cap\mathcal{W}_f(R)$ satisfy $\lambda_f(X_2) = \sigma(X_2; R)$. Upon crossing this locus, the single state $X_2$ bifurcates into two distinct states, giving rise to an intermediate $O$-shock curve that emerges within the $f$-shock curve $[A_2, Z]$, see Fig.~\ref{fig:Gamma7}(b)). This $O$-shock is not admissible under the viscous profile criterion, resulting in two admissible $f$-shock segments separated by a non-admissible $O$-shock. Under this configuration, only one intersection state $X_2$ remains (denoted $X_3$ in Fig.~~\ref{fig:Gamma7}(b)), lying along the $f$-rarefaction curve $[A_2, G]$.
\end{rem}


\subsubsection{Riemann problem solution for \texorpdfstring{$R\in \Gamma_8$}{RET8} }
\label{sec:Gamma8}

This region corresponds to the yellow area in Figs.~\ref{fig:Subdivision_Gamma}(a)-(c). The region $ \Gamma_8 $ is analogous to region $ \Gamma_1 $, as their boundaries are defined similarly. The purple curve $[\mathcal{I}_s^6, V_D]$ separates $\Gamma_7$ and $\Gamma_8$ and marks the boundary where the admissibility of $u$-shocks changes (see \cite{Andrade2018,lozano2018}). The remaining boundaries of $\Gamma_8$ are:  
\begin{itemize}  
    \item The blue curve $[\mathcal{I}_s^6, W]$, which is part of the $s$-inflection locus; and  
    \item The boundary segment of $ \Omega $, given by $ [W, V_D ]$.  
\end{itemize}  
For $L \in [G, W]$ and $R \in \Gamma_8$, the solution to the Riemann problem consists solely of classical wave groups.

\begin{figure}[ht]
    \centering
    \subfigure[Structure of the Riemann solution for for generic $R\in\Gamma_8$ with $R = (0.794132, 0.0291139)$.]
    {\includegraphics[width=0.55\linewidth]{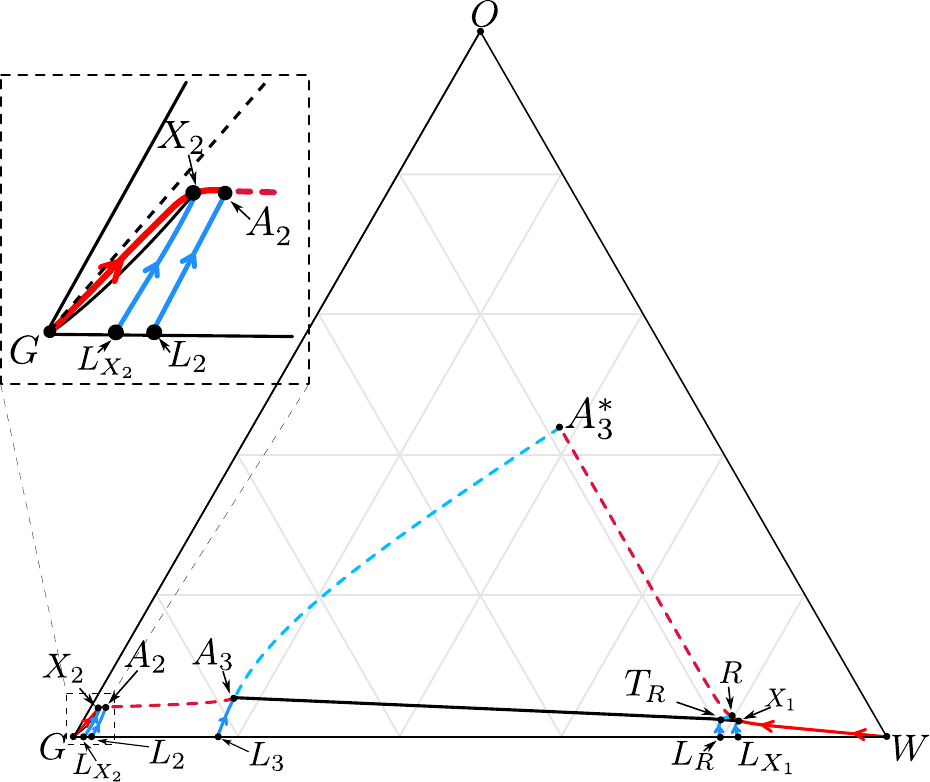}}
    \hspace{0.8mm}
    \caption{Structure of the Riemann solution for $R\in\Gamma_8$. The blue dashed curve (respectively red) represents $s$-shock curves (respectively $f$-shock curves). The blue continuous curve (respectively red) represents $s$-rarefaction curves (respectively the $f$-rarefaction curve). The arrows indicate the increasing characteristic velocity. The black curve $[X,T_R]_{\text{ext}}\cup[T_R,A_3]_{\text{ext}}$ and $[G, X_2]_{\text{ext}}$ represents the $s$-right extension of $\mathcal{W}_f(R)$.
    }
    \label{fig:Gamma8}
\end{figure}

The wave curve $\mathcal{W}_f(R)$ comprises the following components (see Fig.~\ref{fig:Gamma8}):
\begin{itemize}
\item A local branch containing the state $R$, which intersects the $s$-inflection locus at the state $X_1 \in [W, R)$, lying within the a $f$-rarefaction curve;
\item A nonlocal branch containing the state $G$, which intersects the $s$-inflection locus at the state $X_2 \in [G, A_2)$, also within the a $f$-rarefaction curve.
\end{itemize}
Figure~\ref{fig:Gamma8} illustrates the components of $\mathcal{W}_f(R)$ used to construct the solution to the Riemann problem for $L \in [G, W]$:
\begin{itemize}
\item The $f$-rarefaction curves $[W, R)$ and $[G, A_2)$, where $\lambda_f(A_2) = \sigma(A_2; R)$;
\item The $f$-shock curves $(R, A_3^*]$ and $[A_2, A_3]$, where the state $A_3^*$ satisfies $\sigma(R; A_3) = \sigma(R; A_3^*)$.
\end{itemize}

The $s$-right-extension segments consist of $(G, X_2)_{\text{ext}}$, $(X_1, T_R)_{\text{ext}}$, and $[T_R, A_3)_{\text{ext}}$, representing the $s$-right-extensions of the backward $f$-wave curve segments $(G, X_2)$, $[X_1, R)$, and $(R, A_3^*]$.

Let $L_2, L_3, L_R, L_{X_1}$, and $L_{X_2}$ denote the intersection points of the backward $s$-wave curves through states $A_2$, $A_3^*$, $R$, $X_1$, and $X_2$ with the edge $[G, W]$ (see Fig.~\ref{fig:Gamma8}). Then, the structure of the solution for the Riemann problem for $L \in [G, W]$ and $R \in \Gamma_8$ (see Fig.~\ref{fig:Subdivision_Gamma}(b)) consists of:
\begin{enumerate}
    \item[(i)] $L \testright{R_f} R$, for $L = W$;
    \item[(ii)] $L \testright{R_s} Y \testright{R_f} R$, for $L \in(W, L_{X_1}]$ and $Y \in (W, X_1]\subset\mathcal{W}_f(R)$;
    \item[(iii)] $L \testright{R_s} T \xrightarrow{S_s} Y \testright{R_f} R$, for $L \in (L_{X_1}, L_R)$, $T \in (X_1, T_R)_{\text{ext}}$ and $Y \in (X_1, R)\subset\mathcal{W}_f(R)$;
     \item[(iv)] $L \testright{R_s} T \xrightarrow{S_s}R$, for $L = L_R$;
    \item[(v)] $L \testright{R_f} A_2 \xrightarrow{S_f} R$, for $L = G$;
    \item[(vi)] $L \testright{R_s} T \xrightarrow{S_s} Y \testright{R_f} A_2 \xrightarrow{S_f} R$, for $L \in (G, L_{X_2})$, $T \in (G, X_2)_{\text{ext}}$ and $Y \in (G, X_2)\subset\mathcal{W}_f(R)$;
    \item[(vii)] $L \testright{R_s} Y \testright{R_f} A_2 \testright{S_f} R$, for $L \in [L_{X_2}, L_2)$, and  $Y \in [X_2, A_2)\subset\mathcal{W}_f(R)$;
    \item[(viii)] $L \testright{R_s} Y \testright{S_f} R$, for $L \in [L_2, L_3)$ and $Y \in [A_2, A_3)\subset\mathcal{W}_f(R)$;
    \item[(ix)] $L \testright{R_s} T \xrightarrow{S_s} Y \testright{S_f} R$, for $L \in (L_R, L_3)$, $T \in (T_R, A_3)_{\text{ext}}$ and $Y \in (R, A_3^*)\subset\mathcal{W}_f(R)$.
\end{enumerate}
 \begin{rem}
       According to the triple shock rule \cite{Azevedo2014}, $(\sigma(A_3; R) =  \sigma(A_3^*; R) = \sigma(A_3^*; A_3))$, the solution structure for $L = L_3$ (see Fig.~\ref{fig:Gamma8}) is represented by two possible wave sequences in state space:
       $$L_3\testright{R_s} A_3\testright{S_f}R\,\,\mbox{ and }\,\, L_3\testright{R_s} A_3  \xrightarrow{S_s} A_3^*\testright{S_f} R.$$
       These sequences correspond to a unique solution in $xt$-space; see \cite{andrade2016oil, Azevedo2014,fritis2024riemann}.
            \label{rem:gamma_8_with_X_2}
     \end{rem}


\section{Results}\label{sec:Results}

Let us show some applications of the solution classification presented in the previous sections.

\subsection{Oil bank formation}\label{sec:Oil_bank}
The oil bank is a temporary effect characterized by higher oil production and plays an important role in petroleum engineering \cite{kyte_oil_bank1956}. Mathematically, it is represented by a state with oil saturation higher than neighborhood states presented in a solution of the Riemann problem.

\begin{figure}[ht]
    \centering
    \subfigure[Oil saturation profile exhibiting oil bank formation. Saturation values satisfy $\beta>\alpha$ and $\beta> \gamma$]
    {\includegraphics[width=0.35\linewidth]{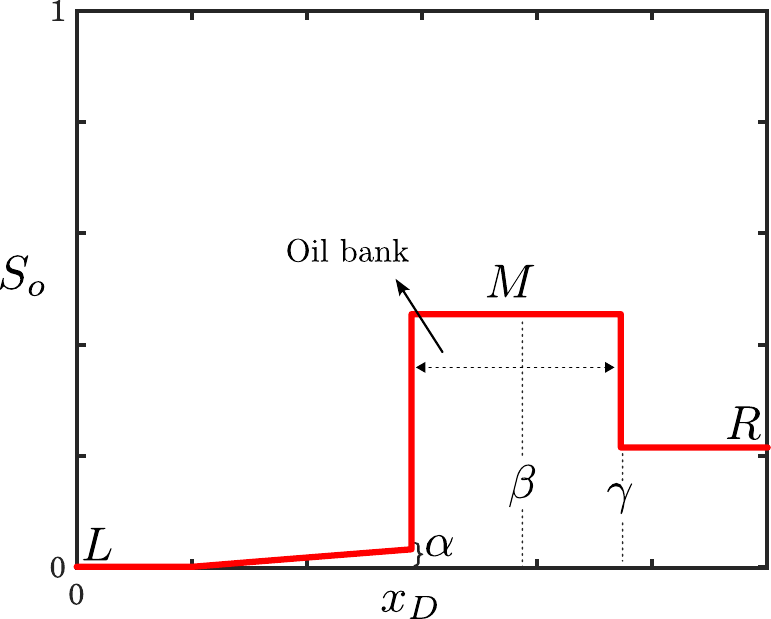}}
    \hspace{5mm}
    \subfigure[Graphical representation of Theorem~\ref{Teo:Teorema_1}. The pink region corresponds to the set of right states that satisfy the theorem, while the blue region does not.]
    {\includegraphics[width=0.5\linewidth]{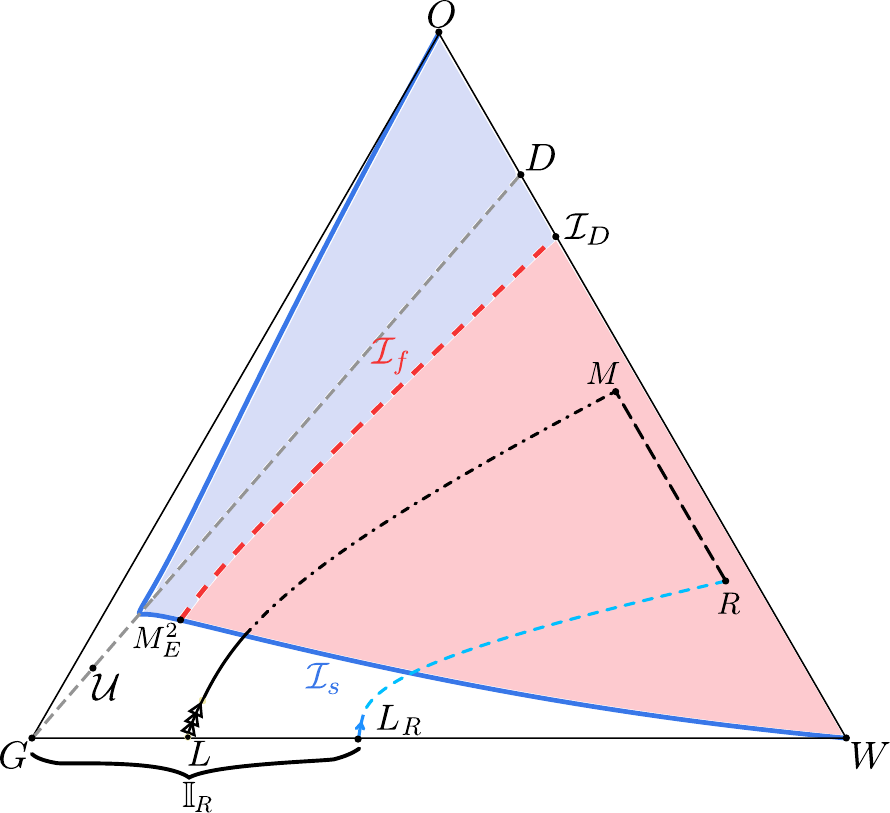}}
    \caption{Oil bank formation. (a) A solution of the Riemann problem with two shocks enclosing a constant state of elevated oil saturation. (b) Validity region of Theorem~\ref{Teo:Teorema_1}, highlighting the right states $R$ and its interval $\mathbb{I}_R$ for which oil bank formation occurs. The black curves represent the composition path for the Riemann problem solution for the left state $L$.}
    \label{fig:Oil_Bank_Theorem}
\end{figure}

In this work, we regard the solution to the Riemann problem that exhibits oil bank formation as a part of the solution satisfying both the following conditions:
\begin{itemize}
    \item[(A)] Two Shock Fronts: the solution consists of a shock followed by a constant state, which is, in turn, followed by another shock.
    \item[(B)] Oil Saturation Condition: the constant state must have an oil saturation higher than the states before and after the shocks.
\end{itemize}
Naturally, both shock fronts must satisfy the velocity compatibility criterion specified in \eqref{eq:compatib} (wave velocities are organized from slow to fast). Figure~\ref{fig:Oil_Bank_Theorem}(a) illustrates an example of an oil saturation profile that satisfies the above conditions. Other cases (oil bank related to rarefaction and shock, for example) can happen for different parameter choices and will be addressed in the future.

Conditions (A) and (B) automatically exclude the regions $\Gamma_1$, $\Gamma_2$, and $\Gamma_3$, which are located above the segment $[G, D]$, since in these regions the oil saturation profile increases monotonically from the left state $L$ to the right state $R$-a behavior incompatible with oil bank formation. Furthermore, when $R \in \Gamma_4$, there exists no corresponding $L \in [G, W]$ that satisfies the required conditions. Specifically:
\begin{itemize}
    \item If $L \in [G, L_R)$, the wave that reaches the right state $R$ in the structure of the Riemann problem solution is an $f$-rarefaction (see Section \ref{sec:Gamma4}, cases (vi), (vii), (viii) and (ix)); hence, condition (A) is not satisfied.
    \item If $L \in [L_R, W]$, the oil saturation profile increases until it reaches the state $R$; thus, condition (B) is not satisfied.
\end{itemize}

For the regions $\Gamma_j$, $j=5,6,7$, and $8$, if $L \in [G, L_R)$, the wave that reaches the right state $R$ in the structure of the Riemann problem solution is an $f$-shock. See, for example:
\begin{itemize}
    \item Cases (v), (vi), (vii), (viii), and (ix) in Section \ref{sec:Gamma5};
    \item Cases (v), (vi), (vii), and (viii) in Section \ref{sec:Gamma6};
    \item Cases (v), (vi), (vii), (viii), and (ix) in Section \ref{sec:Gamma7};
    \item Cases (v), (vi), (vii), (viii), and (ix) in Section \ref{sec:Gamma8}.
\end{itemize}
Since, in these regions, the right state $R$ lies below the $f$-inflection locus $\mathcal{I}_f$, the $f$-shock curve (belonging to the local branch of $\mathcal{W}_f(R)$) necessarily maintains oil saturations greater than those of $R$. Consequently, conditions (A) and (B) are automatically satisfied, yielding the following result:

\begin{theorem}\label{Teo:Teorema_1}
Consider viscosity parameters satisfying \eqref{eq:desi_param_foam} and the shock admissibility criterion given by viscous profile described in Section \ref{sec:visous_profile_Cr}. Let the Riemann problem \eqref{eq:Threephase_1}–\eqref{eq:RP} have a left state $L \in [G, W]$ and a right state $R \in \Gamma$, where $R$ lies below the $f$-inflection locus $\mathcal{I}_f$ (the set $\Gamma$ is defined in Section \ref{subsec2}; see Fig.~\ref{fig:Macro_Classif_3foam}). Let $L_R$ denote the intersection point of the backward $s$-wave curve through $R$ with the edge $[G, W]$. If the left state satisfies $L \in [G, L_R)$, then the solution to the Riemann problem contains an oil bank.
\end{theorem}

Figure~\ref{fig:Oil_Bank_Theorem}(b) illustrates the region where the theorem is valid (highlighted in pink). For the state $R$, the interval $\mathbb{I}_R = [G, L_{R})$ defines the set of left states for which the Riemann problem possesses an oil bank.

\subsection{Comparison of Riemann Problem Solutions}

We present solutions to the Riemann problem for four cases under full-strength foam conditions. For Cases 1 and 2, we analyze foamed wet injection (gas-water co-injection) and compare our solutions with realistic three-phase foam models described in \cite{lyu2021simulation,tang2022foam} obtained using classical waves. For Cases 3 and 4, we study dry foam injection using the Riemann data from \cite{namdar2011method,mayberry2008use}, where the solution construction was based on non-classical waves.

For each case, we validate the analytical estimates with numerical solutions at $t_D = 1$. Direct numerical simulations of System \eqref{eq:Threephase_1} were performed using an implicit finite-difference scheme (FDS), which provides second-order accuracy in both space and time. This was combined with Newton’s method; further details on the implementation can be found in \cite{lambert2020mathematics}.

\subsubsection{Case 1: \texorpdfstring{$R \in \Gamma_6$}{RET6}  and \texorpdfstring{$L \in [L_{\mathfrak{S}}, L_R]$}{LELGLR} }
\label{sec:Case_1}
In this case, we consider a scenario representing a water-wet reservoir. Foam stability increases at high water saturation because the presence of water supports lamellae formation, reducing gas mobility. The injection conditions ($L$) and initial conditions ($R$) are taken from \cite{lyu2021simulation,tang2022foam} (Case 2 of Scenario 3). From Theorem \ref{Teo:Teorema_1}, for these states, the solution to the Riemann problem possesses an oil bank.

\begin{figure}[ht]
	\centering
	\subfigure[Riemann problem solution for $ L =(0.3125 , 0.0001) $ and $ R = (0.775,0.225) $.]
	{\includegraphics[width=0.465\linewidth]{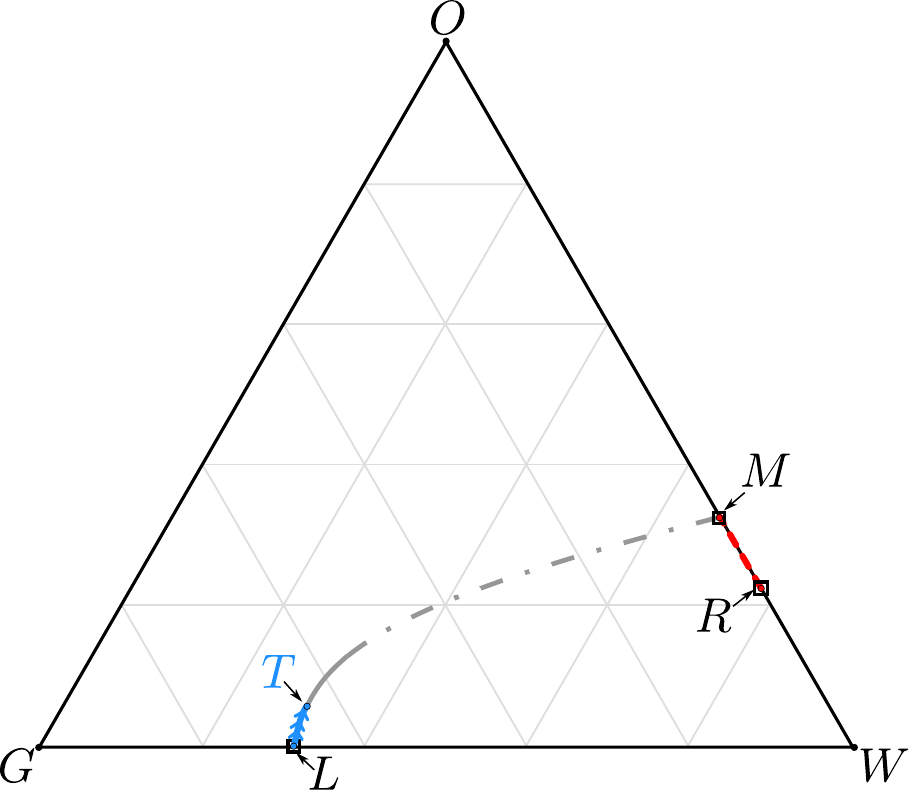}} 
	\hspace{1.5mm}	
 \subfigure[Saturation profiles for water, oil, and gas at $ t_D=1 $.]
	{\includegraphics[width=0.5\linewidth]{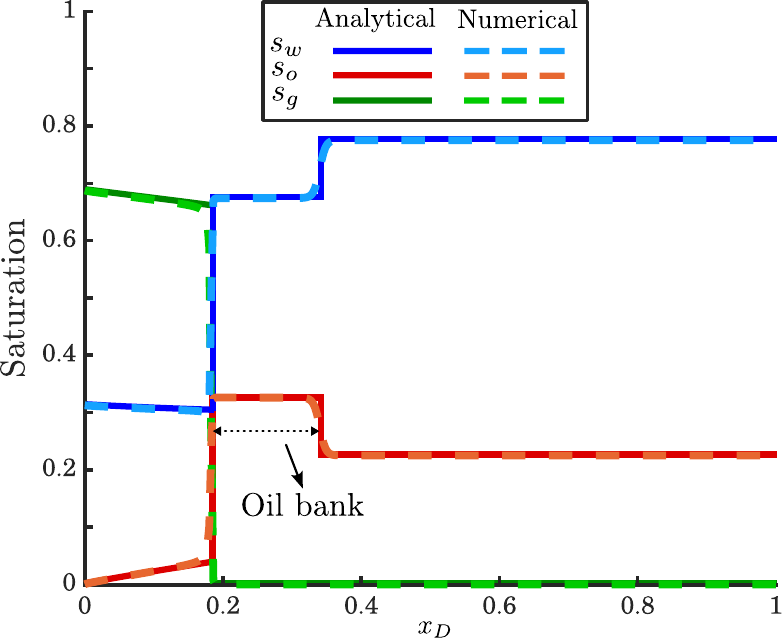}} 
 	\caption{Analytical solution for case 1 with $R\in \Gamma_6$ and $L\in{[L_{\mathfrak{S}}, L_R]}$. (a) The composition path in the saturation triangle that consists of an $s$-rarefaction from $L$ to $T$ follows for an $s$-shock from $T$ to $M$ and an $f$-shock connecting $M$ to $R$.  (b) Analytical profiles (solid curves) compared to numerical simulations (dashed curves).
	}
	\label{fig:Solution_with_foam_Tang2002_Scenario3_caso1}
\end{figure}

Figure~\ref{fig:Solution_with_foam_Tang2002_Scenario3_caso1} illustrates the Riemann problem solution for $ L = (0.3125, 0.0001) $ and $ R = (0.775, 0.225) $. According to the classification presented in Section \ref{sec:Classification}, this case corresponds to $ R \in \Gamma_6 $ and $ L \in [L_{\mathfrak{S}}, L_R] $. The solution consists of an $s$-composite wave from $L$ to $M$ ($ L \testright{R_s} T \xrightarrow{S_s} M $), followed by an $f$-shock wave from $M$ to $R$ ($ M \testright{S_f} R $).  

Figure~\ref{fig:Solution_with_foam_Tang2002_Scenario3_caso1}(a) presents the solution path in the saturation triangle: the blue curve represents $s$-rarefaction from $L$ to $T$, and the red dashed line represents the $f$-shock from $M$ to $R$. The only physically admissible saturation states involved in the displacement are the constant states $L$, $M$, and $R$ (represented by tiny squares in Fig.~\ref{fig:Solution_with_foam_Tang2002_Scenario3_caso1}(a)), as well as the saturations along the rarefaction waves (blue curve in the same figure). The state $T$ (blue), which lies on the $s$-composite wave curve, marks the junction point between rarefaction curve and shock curve, satisfying $ \sigma(T; M) = \lambda_s(T) $.  

Figure~\ref{fig:Solution_with_foam_Tang2002_Scenario3_caso1}(b) compares the water, oil, and gas saturation profiles of the analytical solution (solid curves) with the numerical results (dashed curves) at $t_D = 1$. Initially, the gas profile (green curve) shows a moderate saturation decrease, while, closer to the injection well, a wavefront of gas and water (green and blue curves) displaces oil (red curve). An oil bank forms ahead of this wavefront, producing a solution qualitatively similar to that reported by \cite{lyu2021simulation,tang2022foam} (Figs.~12 and 7) for a more complex three-phase foam flow model.

\subsubsection{Case 2: \texorpdfstring{$R \in \Gamma_3$}{RET3}  and \texorpdfstring{$L \in [L_{1}, L_X]$}{LEL1LX}  }
\label{sec:Case_2}
In this case, we consider a scenario representing a high-oil-content reservoir, reflecting an early-stage enhanced oil recovery (EOR) process. The injection conditions ($L$) and initial conditions ($R$) are taken from \cite{lyu2021simulation,tang2022foam} (Case 2 of Scenario 4). As shown in Section \ref{sec:Oil_bank}, for $R\in\Gamma_3$ the solution to the Riemann problem does not exhibit oil bank formation.   

\begin{figure}[h!]
	\centering
	\subfigure[Riemann problem solution for $ L =(0.3125 , 0.0001) $ and $ R = (0.1875,0.8125) $.]
	{\includegraphics[width=0.465\linewidth]{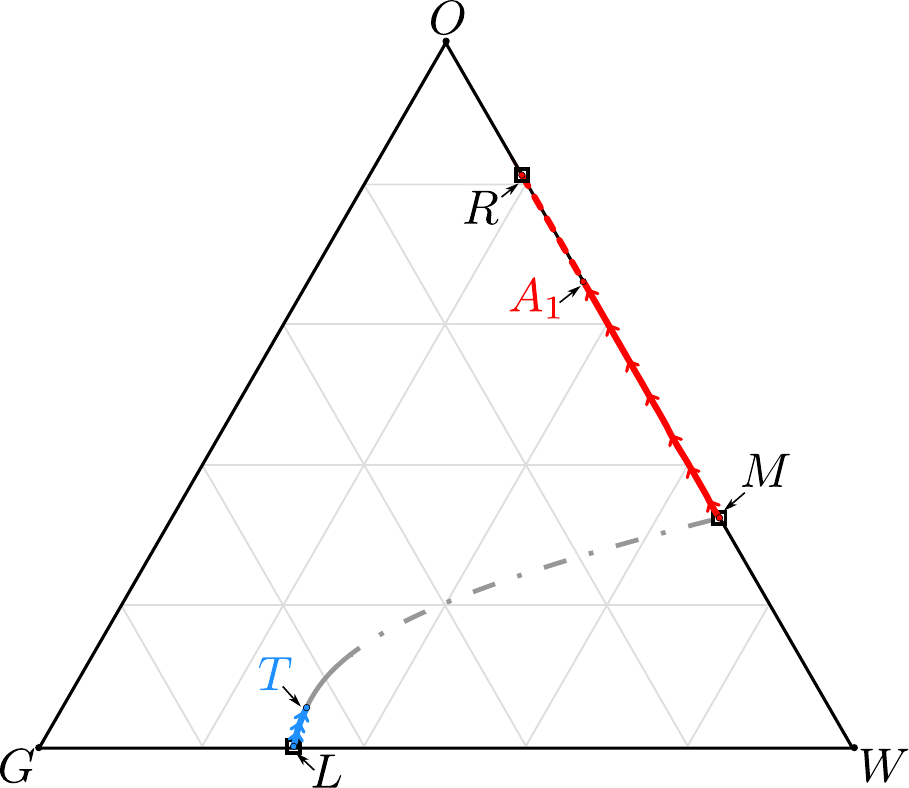}} 
	\hspace{1.5mm}	
 \subfigure[Saturation profiles for water, oil, and gas at $ t_D=1 $.]
	{\includegraphics[width=0.5\linewidth]{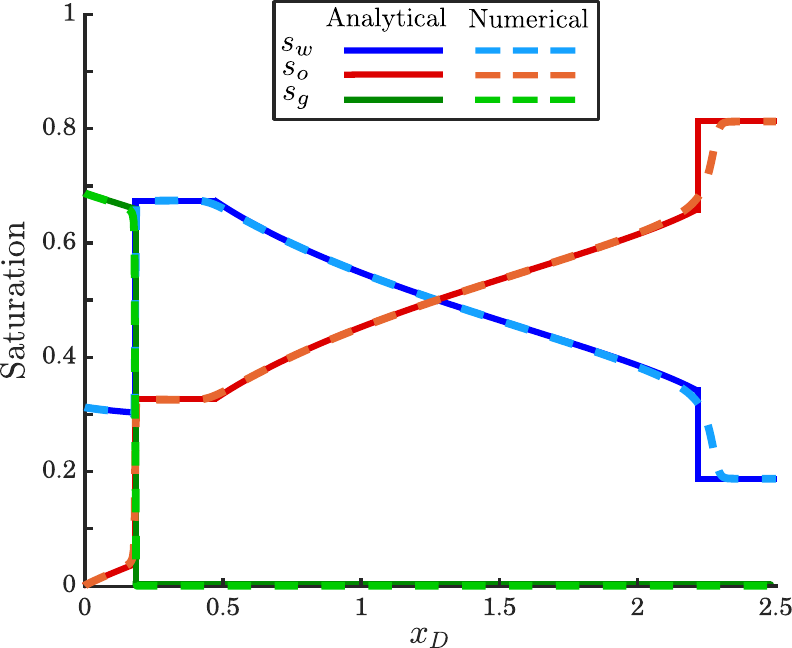}} 
 	\caption{Analytical solution for case 2 with $R\in \Gamma_3$ and $L\in{[L_1,L_{X}]}$. (a) The composition path in the saturation triangle that consists of an $s$-rarefaction from $L$ to $T$ follows for an $s$-shock from $T$ to $M$ follows for an $f$-rarefaction from $M$ to $A_1$ and an $f$-shock connecting $A_1$ to $R$.  (b) Analytical profiles (solid curves) compared to numerical simulations (dashed curves).
	}
	\label{fig:Solution_with_foam_Tang2002_Scenario4_caso1}
\end{figure}

Figure~\ref{fig:Solution_with_foam_Tang2002_Scenario4_caso1} shows the Riemann problem solution for the same state $L = (0.3125, 0.0001)$, but for $R = (0.1875, 0.8125)$. According to the classification presented in Section \ref{sec:Classification}, this case corresponds to $R \in \Gamma_3$ and $L \in [L_1, L_X]$. The solution features an $s$-composite wave from $L$ to $M$ ($ L \testright{R_s} T \xrightarrow{S_s} M $), followed by an $f$-composite wave from $M$ to $R$ ($ M \testright{R_f} A_1 \xrightarrow{S_f} R $). 

Figure~\ref{fig:Solution_with_foam_Tang2002_Scenario4_caso1}(a) presents the solution path in the saturation triangle: the blue curve indicate the $s$-rarefaction from $L$ to $T$, the red curve represent the $f$-rarefaction from $M$ to $A_1$ and the red dashed line represents the $f$-shock from $A_1$ to $R$. The physically admissible saturation states involved in the displacement are the constant states $L$, $M$, and $R$ (represented by tiny squares in Fig.~\ref{fig:Solution_with_foam_Tang2002_Scenario4_caso1}(a)) and the saturations along the rarefaction waves (blue and red curves in Fig.~\ref{fig:Solution_with_foam_Tang2002_Scenario4_caso1}(a)). The states $ {T} $ (blue) and $ {A_1} $ (red), which lie on their respective composite wave curves, mark the junction points between rarefaction waves and shock waves. These points satisfy the conditions $ \sigma({T}; {M}) = \lambda_s({T}) $ and $ \sigma({A_1}; {R}) = \lambda_f({A_1}) $.

Figure~\ref{fig:Solution_with_foam_Tang2002_Scenario4_caso1}(b) compares the water, oil, and gas saturation profiles of the analytical solution (solid curves) with the numerical results (dashed curves) at $t_D = 1$. As in Case 1, the gas profile shows a controlled decrease in saturation (green curves). Then, near the injection well, there is a gas-water wavefront (gas, green curves, and water, blue curves) that displaces the oil (red curves). In addition, a stable foam bank appears ahead of the wavefront. This solution is qualitatively the same as one presented by \cite{lyu2021simulation,tang2022foam} (Figs.~14 and 8) for a more complex three-phase foam flow model.

\subsubsection{Case 3: \texorpdfstring{$ R \in \Gamma_6 $}{RET6}  and \texorpdfstring{$ L = G $}{LG}  }\label{sec:Case_3}
In this case, we consider dry foam injection scenarios (e.g., $CO_2$ or nitrogen injection) in high-water-saturation environments. This scenario enables the study of $CO_2$ sequestration models, particularly focusing on how gas displaces water in nearly saturated formations. 
The injection conditions ($L$) and initial conditions ($R$) are taken from \cite{namdar2011method} (Section 5.2.3, Fig. 3), where the authors employ the Corey model with linear relative permeabilities, in contrast to the quadratic model used in this work. From Theorem \ref{Teo:Teorema_1}, for these $L$ and $R$ states, the solution to the Riemann problem exhibits an oil bank formation.

\begin{figure}[h!]
	\centering
	\subfigure[Riemann problem solution for $ L =(0 , 0) $ and $ R = (0.9123, 0.0875) $.]
	{\includegraphics[width=0.465\linewidth]{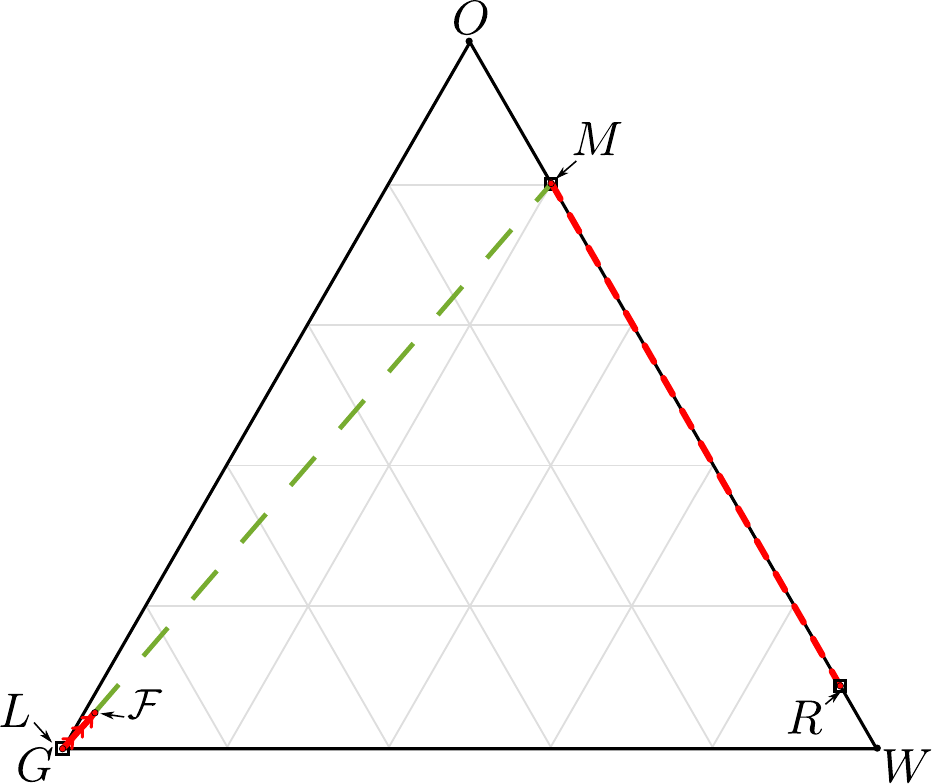}} 
	\hspace{1.5mm}	
 \subfigure[Saturation profiles for water, oil, and gas at $ t_D=1 $.]
	{\includegraphics[width=0.45\linewidth]{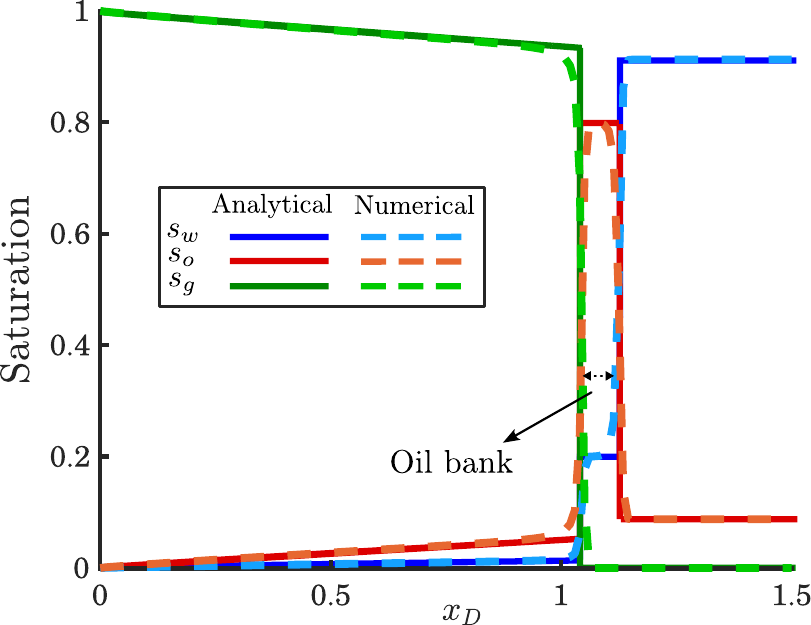}} 
 	\caption{Analytical solution for Case 3 with $R\in \Gamma_6$ and $L=G$. (a) The composition path in the saturation triangle that consists of an $f$-rarefaction from $G$ to $\mathcal{F}$ follows for a $u$-shock from $\mathcal{F}$ to $M$ and an $f$-shock connecting $M$ to $R$.  (b) Analytical profiles (solid curves) compared to numerical simulations (dashed curves).
	}
	\label{fig:Solution_Namdar_I_2_J_1}
\end{figure}

Figure~\ref{fig:Solution_Namdar_I_2_J_1} illustrates the Riemann problem solution for $ L = (0, 0) $ and $ R = (0.9123, 0.0875) $. According to the classification presented in Section \ref{sec:Classification}, this case corresponds to $ R \in \Gamma_6 $ and $ L = G $. The solution features a $ u $-composite wave from $L$ to $M$ ($ L \testright{R_f} \mathcal{F} \xrightarrow{S_u} M $), followed by an $f$-shock wave from $M$ to $R$  ($ M \testright{S_f} R $). 

Figure~\ref{fig:Solution_Namdar_I_2_J_1}(a) presents the solution path in the saturation triangle: the red curve indicate the $f$-rarefaction from $L$ to $ \mathcal{F} $, the green dashed line represents the $ u $-shock from $ \mathcal{F} $ to $M$, and the red dashed line represents the $f$-shock from $M$ to $R$. The constant states $L$, $M$, and $R$ (represented by tiny squares in Fig.~\ref{fig:Solution_Namdar_I_2_J_1}(a)) and the saturations along the rarefaction waves (red curve in Fig.~\ref{fig:Solution_Namdar_I_2_J_1}(a)) are the physically admissible saturation states in the displacement. The state $ \mathcal{F} $, which lies on the $ u $-composite wave curves, marks the junction point between the rarefaction curve and the shock curve, satisfying $ \sigma(\mathcal{F}; M) = \lambda_f(\mathcal{F}) $.  

Figure~\ref{fig:Solution_Namdar_I_2_J_1}(b) compares the water, oil, and gas saturation profiles of the analytical solution (solid curves) with the numerical results (dashed curves) at $ t_D = 1 $. The displacement behavior is similar to that observed in earlier cases: an oil bank (red curve) forms ahead of a wavefront of gas and water (green and blue curves), with no gas production until all recoverable oil is displaced. However, as expected, the wavefront of gas is faster in these cases (dry gas injection) than in cases addressed in sections \ref{sec:Case_1} and \ref{sec:Case_2} (gas + water injection).

\subsubsection{Case 4: \texorpdfstring{ $ R \in \Gamma_3 $ }{RET3}  and \texorpdfstring{ $ L = [G, L_{\mathcal{F}}] $}{LEGLF}}
\label{sec:Case_4}

This case presents a scenario where the injection state is nearly pure gas with a small amount of water, and the reservoir initially contains mostly oil with some water. This setup can be understood as part of the FAWAG (Foam-Assisted Water-Alternating-Gas) procedure, specifically during the stage where foamed gas is injected. The injection conditions ($L$) and initial conditions ($R$) are taken from \cite{mayberry2008use} (Fig. 13), where the authors employ the Corey model with linear relative permeabilities, in contrast to the quadratic model used in the current work. As shown in Section \ref{sec:Oil_bank}, for $R\in\Gamma_3$ the solution to the Riemann problem does not possess an oil bank.

\begin{figure}[h!]
	\centering
	\subfigure[Riemann problem solution for $ L =(0.0125,0.0001) $ and $ R = (0.125, 0.875) $.]
	{\includegraphics[width=0.465\linewidth]{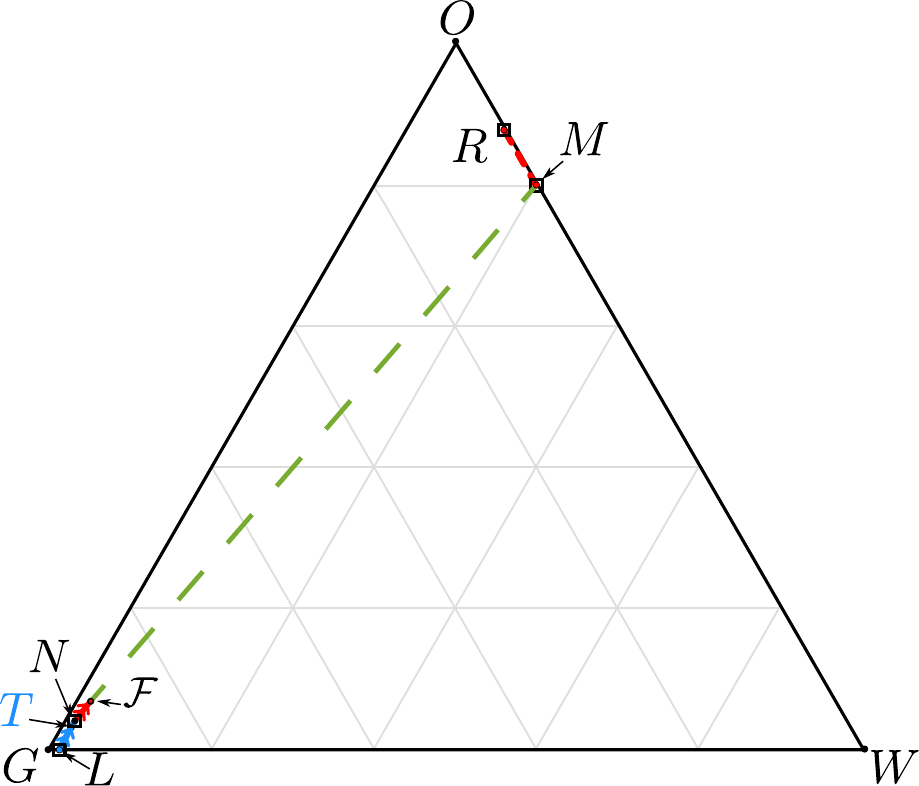}} 
	\hspace{1.5mm}	
 \subfigure[Saturation profiles for water, oil, and gas at $ t_D=1 $.]
	{\includegraphics[width=0.45\linewidth]{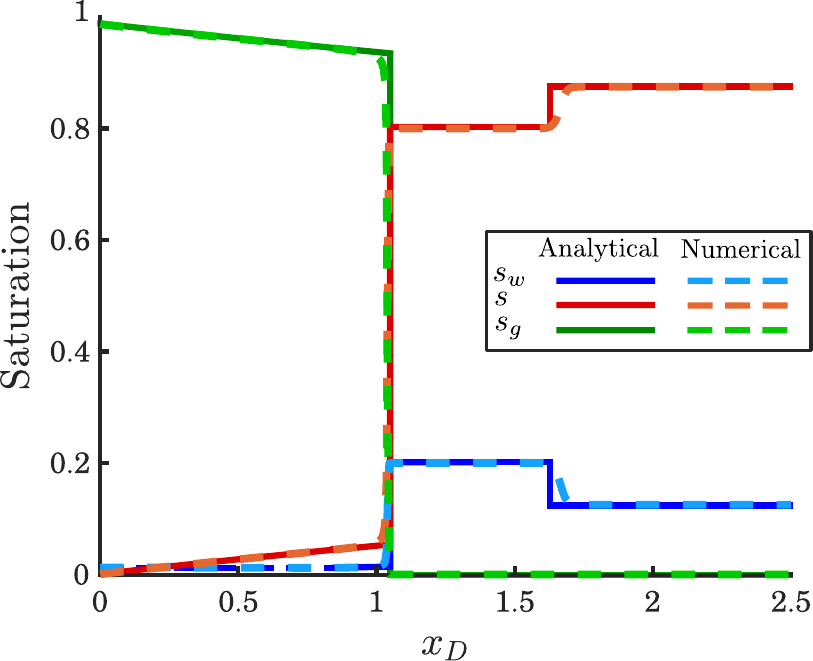}} 
 	\caption{Analytical solution for case 4 with $R\in \Gamma_3$ and $L\in{[G,L_{\mathcal{F}}]}$. (a) The composition path in the saturation triangle that consists of an $s$-rarefaction from $L$ to $T$ follows for an $s$-shock from $T$ to $N$ follows for an $f$-rarefaction from $N$ to $\mathcal{F}$ follows for a $u$-shock from $\mathcal{F}$ to $M$ and an $f$-shock connecting $M$ to $R$.  (b) Analytical profiles (solid curves) compared to numerical simulations (dashed curves).
	}
	\label{fig:Solution__Mayberry_I1_J4}
\end{figure}

Figure \ref{fig:Solution__Mayberry_I1_J4} illustrates the Riemann problem solution for $ L =(0.0125,0.0001) $ and $ R = (0.125, 0.875) $. According to the classification presented in Section \ref{sec:Classification}, this case corresponds to $ R \in \Gamma_3 $ and $ L = [G, L_{\mathcal{F}}] $. The solution features an $s$-composite wave from $L$ to $ N $ ($ L \testright{R_s} T \xrightarrow{S_s} N $), followed by an $u$-composite wave from $ N $ to $M$ ($ N \testright{R_f} \mathcal{F} \xrightarrow{S_u} M $), followed by an $f$-shock wave from $M$ to $R$ $(M \testright{S_f} R)$. 

Figure~\ref{fig:Solution__Mayberry_I1_J4}(a) presents the solution path in the saturation triangle: The blue and red curves represent the $ s$-rarefaction from $ L$ to $ T$ and the $ f$-rarefaction from $ N$ to $ \mathcal{F}$, respectively. The green dashed line corresponds to the $ u$-shock from $ \mathcal{F}$ to $ M$, while the red dashed line represents the $ f$-shock from $ M$ to $ R$. In this case, there are four constant states $L$, $ N$, $ M$, and $ R$ denoted by tiny squares in Fig.~\ref{fig:Solution__Mayberry_I1_J4}(a). The saturations along the slow and fast rarefaction waves (blue and red curves, respectively) are the physically admissible saturation states involved in the displacement. However, the shock amplitude between $T$ (blue) and $ N$ is very small and imperceptible in the simulation shown in Fig.~\ref{fig:Solution__Mayberry_I1_J4}(b). The states $ T$ (blue) and $ \mathcal{F}$, belonging to their respective composite wave curves, mark the junction points between the rarefaction curve and the shock curve, satisfying $ \sigma(T; N) = \lambda_s(T)$ and $ \sigma(\mathcal{F}; M) = \lambda_f(\mathcal{F})$.

Figure~\ref{fig:Solution__Mayberry_I1_J4}(b) compares the water, oil, and gas saturation profiles of the analytical solution (solid curves) with the numerical results (dashed curves) at $ t_D = 1$. The displacement behavior is consistent with that observed in earlier cases: a stable wavefront of gas and water (green and blue curves) displaces the oil, with no gas production until all recoverable oil is displaced. However, as expected, the wavefront of gas is faster in these cases (dry gas injection) than in cases addressed in sections \ref{sec:Case_1} and \ref{sec:Case_2} (gas + water injection).

\section{Discussion and Conclusions}
\label{sec:discusstion}

In this work, we analyze a model for three-phase foam flow in a porous medium with nonlinear relative permeabilities. The analysis leads to the solution of the Riemann problem for the injection of foamed gas and water mixtures under a wide range of initial conditions relevant to industrial applications. Our solutions were validated in two ways: comparison with literature results and numerical simulation. Comparison between selected Riemann data cases shows qualitative agreement between our solutions and those of other foam models \cite{tang2019three,lyu2021simulation}.

Our solution is valid for realistic foam flow models and can be employed in calibrating numerical simulators, performing uncertainty quantification, and analyzing the impact of variations in physical parameters.

The results hold significant industrial relevance. 
Notably, we show that wet foam (Section \ref{sec:Case_2}) displaces faster than dry foam (Section \ref{sec:Case_4}), consistent with findings in \cite{lyu2021simulation,tang2022foam} for more complex models. We identify specific conditions for the formation of oil banks, which occurs only for right states $R$ below the segment $[G, D]$ and left states $L$ within $[G, L_R)$ (Fig.~\ref{fig:Oil_Bank_Theorem}). Oil bank formation is one of the key interests in petroleum engineering. Our analytical estimates for oil bank formation allow determining the breakthrough time, which is important for real-world applications.

\begin{figure}
    \centering
    \includegraphics[width=0.5\linewidth]{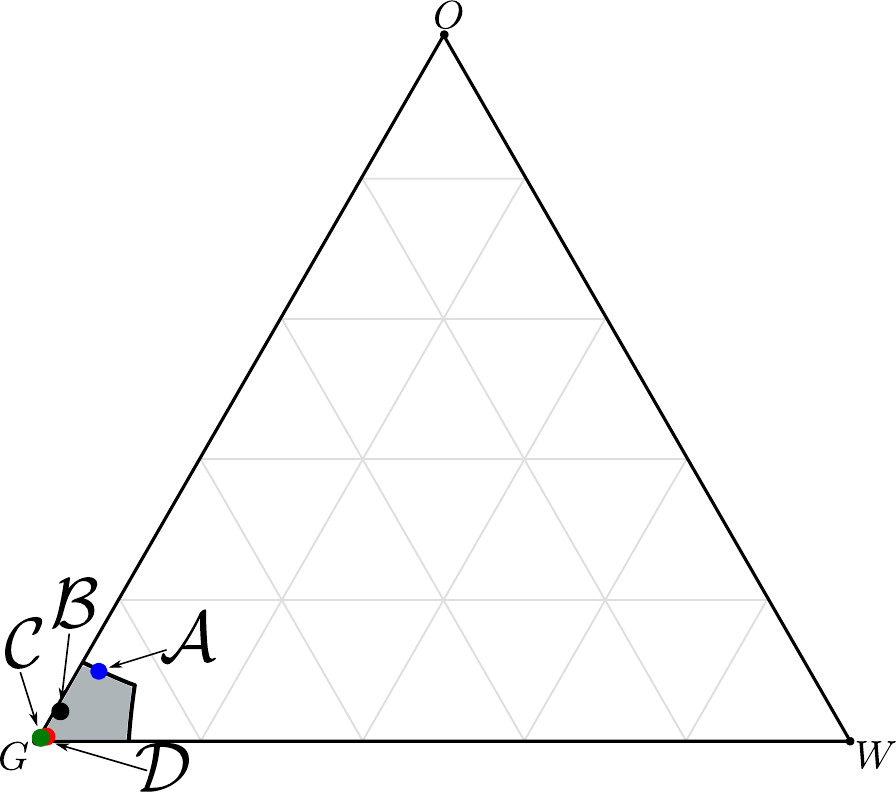}
    \caption{The colored region indicates the area within the saturation triangle where the umbilic point $\mathcal{U}$ is located to satisfy the inequalities in \eqref{eq:desi_param_foam}. 
    Points correspond to: $\mathcal{A}$(blue) - Table \ref{tab:table1}; 
    $\mathcal{B}$ (black) -  \cite{izadi2021investigation}; 
    $\mathcal{C}$ (green) - \cite{namdar2011method};
    $\mathcal{D}$ (red) - \cite{lyu2021simulation,tang2022foam}.}
    \label{fig:Region_parameters_adm}
\end{figure}

The mathematical approach utilized here is non-classical conservation law theory. We show that the solution configurations remain stable under variations in physics parameters. In particular, this is so for changes in viscosity parameters such that the umbilic point remains within the region given by \eqref{eq:desi_param_foam}. In Fig.~\ref{fig:Region_parameters_adm}, this region is indicated with a gray color. Also, we plot a few umbilic points represented by $ \mathcal{A} $ (blue) for synthetic viscosities (Table \ref{tab:table1}), 
$\mathcal{B}$ (black) for viscosities from the Lisama field \cite{izadi2021investigation} , and for experimental data, $\mathcal{C}$ (green) from  \cite{namdar2011method} and $\mathcal{D}$ (red) from \cite{lyu2021simulation,tang2022foam}.
Thus, our analysis encompasses applications for realistic viscosity parameters.


\subsection*{Acknowledgment}

The authors thank Prof. Bradley Plohr for assistance with ELI interactive solver \cite{ELI_web}.

G.C. and L.L. gratefully acknowledge support from Shell Brasil through the projects ``Avançando na modelagem matemática e computacional para apoiar a implementação da tecnologia `Foam-assisted WAG' em reservatórios do Pré-sal'' (ANP 23518-4) at UFJF, and the strategic importance of the support given by ANP through the R\&D levy regulation.

G.C. was partly supported by CNPq grants 306970/2022-8, 405366/2021-3, and FAPEMIG grant APQ-00206-24.

D.M. was partly supported by CAPES grant 88881.156518/2017-01, by CNPq under grants 405366/2021-3, 306566/2019-2, and by FAPERJ under grants E-26/210.738/2014, E-26/202.764/2017, E-26/201.159/2021.

\subsection*{CRediT authorship contribution statement}
\textbf{L. F. Lozano:} 
    Conceptualization, 
    Methodology, 
    Software, 
    Formal analysis, 
    Investigation, 
    Writing – original draft, 
    Writing – review \& editing. 
\textbf{G. Chapiro:} 
    Conceptualization, 
    Methodology, 
    Supervision, 
    Investigation, 
    Project administration, 
    Funding acquisition,
    Writing – original draft, 
    Writing – review \& editing.
\textbf{D. Marchesin}
    Conceptualization, 
    Methodology, 
    Supervision,
    Investigation, 
    Writing – original draft, 
    Writing – review \& editing. 
    
\subsection*{Data Availability Statement} The data that support the findings of this study are available within the
 article.

\subsection*{Declarations}

\subsection*{Conflict of interest statement}  The author has no conflicts to disclose.



\printbibliography 
\end{document}